\documentclass[a4paper,12pt]{amsart}

\usepackage{float}
\usepackage{euscript,eufrak,verbatim}
\usepackage{graphicx}
\usepackage[usenames]{color}
\usepackage[colorlinks,linkcolor=red,anchorcolor=blue,citecolor=blue]{hyperref}
\usepackage{amsmath}
\usepackage{amsthm}
\usepackage[all]{xy}
\usepackage{graphicx} 
\usepackage{amssymb}

\usepackage{mathrsfs}
\usepackage{amscd}

\makeatletter
\newcommand{\svdots}{%
  \vbox{\fontsize{\sf@size}{\sf@size pt}\linespread{0.3}\selectfont
    \kern0.2\baselineskip
    \hbox{.}\hbox{.}\hbox{.}%
    \kern0.1\baselineskip
  }%
}
\makeatother

\theoremstyle{plain}
\newtheorem{main theorem}{Main Theorem}
\newtheorem{theorem}{Theorem}[section]
\newtheorem{lemma}[theorem]{Lemma}
\newtheorem{conjecture}[theorem]{Conjecture}
\newtheorem{corollary}[theorem]{Corollary}
\newtheorem{proposition}[theorem]{Proposition}
\newtheorem{claim}[theorem]{Claim}

\newtheorem{lemma-definition}[theorem]{Lemma-Definition}
\theoremstyle{definition}
\newtheorem{definition}[theorem]{Definition}
\newtheorem{remark}[theorem]{Remark}

\makeatletter 
\@addtoreset{equation}{section}
\numberwithin{equation}{section}

\newcommand{\norm}[1]{\left\lVert#1\right\rVert}

\newcommand{\diam}{\mathrm{Diam}}

\newcommand{\mdim}{\mathrm{mdim}}

\newcommand{\widim}{\mathrm{Widim}}

\newcommand{\umdimm}{\overline{\mathrm{mdim}}_{\mathrm{M}}}
\newcommand{\lmdimm}{\underline{\mathrm{mdim}}_{\mathrm{M}}}

\newcommand{\umdimh}{\overline{\mathrm{mdim}}_{\mathrm{H}}}
\newcommand{\lmdimh}{\underline{\mathrm{mdim}}_{\mathrm{H}}}
\newcommand{\urdim}{\overline{\mathrm{rdim}}}
\newcommand{\lrdim}{\underline{\mathrm{rdim}}}
\newcommand{\dimh}{\dim_{\mathrm{H}}}
\newcommand{\umdimhl}{\overline{\mathrm{mdim}}_{\mathrm{H}, L^1}}
\newcommand{\lmdimhl}{\underline{\mathrm{mdim}}_{\mathrm{H},L^1}}

\title[Mean dimension with potential of $\mathbb{R}^d$-actions: I]
{Variational principle for mean dimension with potential of $\mathbb{R}^d$-actions: I}

\author{Masaki Tsukamoto}

\address
{Department of Mathematics, Kyoto University, Kitashirakawa Oiwake-cho, Sakyo-ku, Kyoto 606-8502, Japan}

\email{tsukamoto@math.kyoto-u.ac.jp}

\begin{document}



\thanks{The author was supported by JSPS KAKENHI JP21K03227.}

\maketitle

\begin{abstract}
We develop a variational principle for mean dimension with potential of $\mathbb{R}^d$-actions.
We prove that mean dimension with potential is bounded from above by the supremum of 
the sum of rate distortion dimension and a potential term.
A basic strategy of the proof is the same as the case of $\mathbb{Z}$-actions.
However measure theoretic details are more involved because $\mathbb{R}^d$ is a continuous group.
We also establish several basic properties of metric mean dimension with potential and 
mean Hausdorff dimension with potential for $\mathbb{R}^d$-actions.
\end{abstract}

\section{Introduction} \label{section: introduction}

\subsection{Background: the case of $\mathbb{Z}$-actions}   \label{subsection: background}

The purpose of this paper is to develop a theory of variational principle for 
mean dimension with potential of $\mathbb{R}^d$-actions.
First we review the theory already established in the case of $\mathbb{Z}$-actions.

Mean dimension is a topological invariant of dynamical systems introduced by Gromov 
\cite{Gromov} at the end of the last century.
It is the number of parameters \textit{per unit time} for describing given dynamical systems.
Mean dimension has several applications to topological dynamics, 
most notably in the embedding problem of dynamical systems 
\cite{Lindenstrauss--Weiss, Lindenstrauss, Gutman--Tsukamoto embedding, 
Gutman--Qiao--Tsukamoto}.

Lindenstrauss and the author \cite{Lindenstrauss--Tsukamoto IEEE, Lindenstrauss--Tsukamoto double, Tsukamoto potential}
began to develop the variational principle in mean dimension theory.
Let $\mathcal{X}$ be a compact metrizable space, 
and let $T\colon \mathcal{X}\to  \mathcal{X}$ be a homeomorphism of $\mathcal{X}$.
The classical variational principle \cite{Goodwyn, Dinaburg, Goodman} states that 
the topological entropy $h_{\mathrm{top}}(T)$ is equal to the supremum of the Kolmogorov--Sinai entropy $h_\mu(T)$
over all invariant probability measures $\mu$:
\begin{equation}    \label{eq: variational principle for entropy}
    h_{\mathrm{top}}(T) = \sup_{\mu\in \mathscr{M}^T(\mathcal{X})} h_{\mu}(T), 
\end{equation}
where $\mathscr{M}^T(\mathcal{X})$ denotes the set of all $T$-invariant Borel probability measures on $\mathcal{X}$.
Ruelle \cite{Ruelle} and then Walters \cite{Walters} generalized (\ref{eq: variational principle for entropy}) to pressure:
Let $\varphi \colon \mathcal{X} \to \mathbb{R}$ be a continuous function, 
and we denote by $P_T(\varphi)$ the topological pressure of $(X, T, \varphi)$.
Then 
\begin{equation}  \label{eq: variational principle for pressure}
   P_T(\varphi) = \sup_{\mu\in \mathscr{M}^T(\mathcal{X})}\left(h_\mu(T) + \int_{\mathcal{X}} \varphi \, d\mu\right).
\end{equation}

In the classical variational principles (\ref{eq: variational principle for entropy}) and (\ref{eq: variational principle for pressure}),
the quantities $h_{\mathrm{top}}(T)$ and $P_T(\varphi)$
in the left-hand sides are topological invariants of dynamical systems. 
The Kolomogorov--Sinai entropy in the right-hand side is an information theoretic quantity.
Therefore (\ref{eq: variational principle for entropy}) and (\ref{eq: variational principle for pressure}) connect topological dynamics to 
information theory.
Lindenstrauss and the author tried to find an analogous structure in mean dimension theory.
(See also the paper of Gutman--Śpiewak \cite{Gutman--Spiewak} 
for a connection between mean dimension and information theory.)
In the papers \cite{Lindenstrauss--Tsukamoto IEEE, Lindenstrauss--Tsukamoto double}
they found that \textit{rate distortion theory} provides a fruitful framework for the problem.
This is a branch of information theory studying \textit{lossy} data compression method under a distortion constraint.

Let $T:\mathcal{X}\to \mathcal{X}$ be a homeomorphism on a compact metrizable space $\mathcal{X}$ as in the above.
We denote the mean dimension of $(\mathcal{X}, T)$ by $\mdim(\mathcal{X}, T)$.
We would like to connect it to some information theoretic quantity.
We define $\mathscr{D}(\mathcal{X})$ as the set of all metrics (distance functions) on $\mathcal{X}$ compatible with the 
given topology. 
Let $\mathbf{d}\in \mathscr{D}(\mathcal{X})$ and $\mu\in \mathcal{M}^T(\mathcal{X})$.
We randomly choose a point $x\in \mathcal{X}$ according to the distribution $\mu$ and consider the orbit $\{T^n x\}_{n\in \mathbb{Z}}$.
For $\varepsilon>0$, we define the \textit{rate distortion function} $R(\mathbf{d}, \mu, \varepsilon)$ 
as the minimum number of bits per unit time for describing $\{T^n x\}_{n\in \mathbb{Z}}$ with average distortion 
bounded by $\varepsilon$ with respect to $\mathbf{d}$.
See \S \ref{subsection: rate distortion theory} for the precise definition of $R(\mathbf{d},\mu,\varepsilon)$ in the case of 
$\mathbb{R}^d$-actions.

We define the \textbf{upper and lower rate distortion dimensions} by 
\[  \urdim(\mathcal{X}, T, \mathbf{d}, \mu) = \limsup_{\varepsilon\to 0} \frac{R(\mathbf{d}, \mu, \varepsilon)}{\log(1/\varepsilon)},\quad 
    \lrdim(\mathcal{X}, T, \mathbf{d}, \mu) = \liminf_{\varepsilon\to 0} \frac{R(\mathbf{d}, \mu, \varepsilon)}{\log(1/\varepsilon)}. \]
Rate distortion dimension was first introduced by Kawabata--Dembo \cite{Kawabata--Dembo}.

Lindenstrauss and the author \cite[Corollary 3.13]{Lindenstrauss--Tsukamoto double} proved that 
\begin{equation}  \label{eq: variational principle for mean dimension}
   \mdim(\mathcal{X}, T) \leq \sup_{\mu \in \mathcal{M}^T(\mathcal{X})} \lrdim(\mathcal{X}, T, \mathbf{d}, \mu)
\end{equation}
for any metric $\mathbf{d}$ on $\mathcal{X}$ compatible with the given topology.
Moreover they proved that if $(\mathcal{X}, T)$ is a free minimal dynamical system then \cite[Theorem 1.1]{Lindenstrauss--Tsukamoto double}
\begin{equation}  \label{eq: double variational principle for mean dimension} 
   \begin{split}
    \mdim(\mathcal{X}, T) 
    &= \min_{\mathbf{d}\in \mathscr{D}(\mathcal{X})}\left(\sup_{\mu \in \mathcal{M}^T(\mathcal{X})} \lrdim(\mathcal{X}, T, \mathbf{d}, \mu)\right) \\
    &= \min_{\mathbf{d}\in \mathscr{D}(\mathcal{X})}\left(\sup_{\mu \in \mathcal{M}^T(\mathcal{X})} \urdim(\mathcal{X}, T, \mathbf{d}, \mu)\right).
   \end{split} 
\end{equation}
They called this “double variational principle” because it involves a minimax problem with respect to the \textit{two} 
variables $\mathbf{d}$ and $\mu$. We conjecture that (\ref{eq: double variational principle for mean dimension}) holds for all dynamical systems
without any additional assumption.

The author \cite{Tsukamoto potential} generalized (\ref{eq: variational principle for mean dimension}) and 
(\ref{eq: double variational principle for mean dimension}) to \textit{mean dimension with potential}, which is 
a mean dimension analogue of topological pressure.
Let $\varphi\colon \mathcal{X}\to \mathbb{R}$ be a continuous function. 
The paper \cite{Tsukamoto potential} introduced mean dimension with potential (denoted by $\mdim(\mathcal{X}, T, \varphi)$) and 
proved that \cite[Corollary 1.7]{Tsukamoto potential}
\begin{equation} \label{eq: variational principle for mean dimension with potential}
    \mdim(\mathcal{X}, T, \varphi) \leq \sup_{\mu \in \mathscr{M}^T(\mathcal{X})} 
    \left(\lrdim(\mathcal{X}, T, \mathbf{d}, \mu) + \int_\mathcal{X} \varphi\, d\mu\right). 
\end{equation}    
Moreover, if $(\mathcal{X}, T)$ is a free minimal dynamical system then \cite[Theorem 1.1]{Tsukamoto potential}
\begin{equation}  \label{eq: double variational principle for mean dimension with potential}
     \begin{split}
     \mdim(\mathcal{X}, T, \varphi) 
    &= \min_{\mathbf{d}\in \mathscr{D}(\mathcal{X})}
    \left\{\sup_{\mu \in \mathscr{M}^T(\mathcal{X})} 
    \left(\lrdim(\mathcal{X}, T, \mathbf{d}, \mu) + \int_\mathcal{X} \varphi\, d\mu\right)\right\} \\
    &= \min_{\mathbf{d}\in \mathscr{D}(\mathcal{X})}
    \left\{\sup_{\mu \in \mathscr{M}^T(\mathcal{X})} 
    \left(\urdim(\mathcal{X}, T, \mathbf{d}, \mu) + \int_\mathcal{X} \varphi\, d\mu\right)\right\}.     
     \end{split}
\end{equation}
We also conjecture that this holds for all dynamical systems.

The main purpose of this paper is to generalize the above 
(\ref{eq: variational principle for mean dimension with potential}) to $\mathbb{R}^d$-actions.
We think that we can also generalize the \textit{double variational principle} 
(\ref{eq: double variational principle for mean dimension with potential}) to free minimal $\mathbb{R}^d$-actions.
However it requires a technically heavy work.
We postpone it to Part II of this series of papers.
In this paper we concentrate on the inequality (\ref{eq: variational principle for mean dimension with potential}).

The motivation to generalize (\ref{eq: variational principle for mean dimension with potential}) and 
(\ref{eq: double variational principle for mean dimension with potential}) to $\mathbb{R}^d$-actions
comes from the fact that many natural examples of mean dimension theory are rooted in
geometric analysis \cite{Gromov, Matsuo--Tsukamoto Brody curves, Tsukamoto Brody curves}.
In geometric analysis we usually consider actions of groups more complicated than $\mathbb{Z}$.
Maybe $\mathbb{R}^d$-actions are the most basic case.
We plan to apply the results of this paper to geometric examples of 
\cite{Gromov, Matsuo--Tsukamoto Brody curves, Tsukamoto Brody curves} in a future paper.

Since $\mathbb{R}^d$ is a continuous group, several new technical difficulties appear.
Especially measure theoretic details are more complicated in the case of $\mathbb{R}^d$-actions 
than in the case of $\mathbb{Z}$-actions.
A main task of this paper is to establish such details.

We would like to mention the paper of Huo--Yuan \cite{Huo--Yuan}.
They develop the variational principle for mean dimension of $\mathbb{Z}^d$-actions.
In \S \ref{section: mean dimension of Z^d-actions} and \S \ref{section: proof of mdim is bounded by mdimh} 
we also touch the case of $\mathbb{Z}^d$-actions.
Some results in these sections were already studied in \cite{Huo--Yuan}.

\subsection{Mean dimension with potential of $\mathbb{R}^d$-actions}   
\label{subsection: mean dimension with potential of R^d-actions}

\begin{figure}[h] 
    \centering
    \includegraphics[width=3.0in]{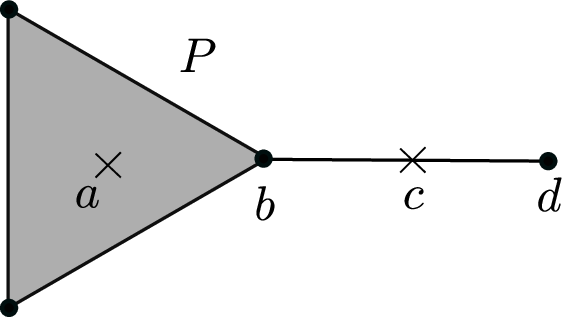}
    \caption{Here $P$ has four vertexes (denoted by dots), four $1$-dimensional simplexes and one $2$-dimensional simplex. 
    The points $b$ and $d$ are vertexes of $P$ wheres $a$ and $c$ are not. 
    We have $\dim_a P = \dim_b P =2$ and $\dim_c P = \dim_d P =1$.}    \label{figure: local dimension}
\end{figure}

In this subsection we introduce mean dimension with potential for $\mathbb{R}^d$-actions.
Let $P$ be a finite simplicial complex. (Here “finite” means that the number of faces is finite.
In this paper we do not consider infinite simplicial complexes. Simplicial complexes are always finite.)
For a point $a \in P$ we define $\dim_a P$ as the maximum of $\dim \Delta$ 
where $\Delta$ runs over all simplices of $P$ containing $a$.
We call $\dim_a P$  the \textbf{local dimension} of $P$ at $a$.
See Figure \ref{figure: local dimension}.
(This is the same as \cite[Fig. 1]{Tsukamoto potential}.)

Let $(\mathcal{X}, \mathbf{d})$ be a compact metric space.
Let $\mathcal{Y}$ be a topological space and $f\colon \mathcal{X}\to \mathcal{Y}$ a continuous map.
For a positive number $\varepsilon$ we call $f$ an \textbf{$\varepsilon$-embedding} if we have 
$\diam f^{-1}(y) < \varepsilon$ for all $y\in \mathcal{Y}$.
Let $\varphi: \mathcal{X}\to \mathbb{R}$ be a continuous function.
We define the \textbf{$\varepsilon$-width dimension with potential} by
\begin{equation} \label{eq: widim with potential}
   \begin{split}
     & \widim_\varepsilon(\mathcal{X}, \mathbf{d}, \varphi)  \\ 
     & =  \inf\left\{\max_{x\in \mathcal{X}} \left(\dim_{f(x)} P + \varphi(x)\right) \middle|
    \parbox{3in}{\centering $P$ is a finite simplicial complex and $f:\mathcal{X}\to P$ is an $\varepsilon$-embedding}\right\}.  
   \end{split}
\end{equation}

Let $d$ be a natural number. 
We consider that $\mathbb{R}^d$ is equipped with the Euclidean topology and 
standard additive group structure.
We denote the standard Lebesgue measure on $\mathbb{R}^d$ by $\mathbf{m}$.
Let $T\colon \mathbb{R}^d\times \mathcal{X}\to \mathcal{X}$ be a continuous action of $\mathbb{R}^d$ on 
a compact metrizable space $\mathcal{X}$.
Let $\mathbf{d}$ be a metric on $\mathcal{X}$ compatible with the topology, and 
let $\varphi\colon \mathcal{X}\to \mathbb{R}$ be a continuous function.
For a bounded Borel subset $A \subset \mathbb{R}^d$ we define a new metric $\mathbf{d}_A$ and 
a new function $\varphi_A \colon \mathcal{X}\to \mathbb{R}$ by 
\[  \mathbf{d}_A(x, y) = \sup_{u\in A} \mathbf{d}(T^u x, T^u y), \quad 
     \varphi_A(x) = \int_{A} \varphi(T^u x)\, d\mathbf{m}(u). \]

If $\varphi(x) \geq 0$ for all $x\in \mathcal{X}$ then we have:
\begin{enumerate}
    \item  \textbf{Subadditivity:} For bounded Borel subsets $A, B \subset \mathbb{R}^d$ 
    \[     \widim_\varepsilon\left(\mathcal{X}, \mathbf{d}_{A \cup B}, \varphi_{A \cup B}\right)
            \leq  \widim_\varepsilon\left(\mathcal{X}, \mathbf{d}_{A}, \varphi_{A}\right)
                   + \widim_\varepsilon\left(\mathcal{X}, \mathbf{d}_{B}, \varphi_{B}\right). \]
     \item  \textbf{Monotonicity:} If $A \subset  B$ then 
    \[   0\leq  \widim_\varepsilon\left(\mathcal{X}, \mathbf{d}_{A}, \varphi_{A}\right) \leq 
         \widim_\varepsilon\left(\mathcal{X}, \mathbf{d}_{B}, \varphi_{B}\right).  \]              
    \item  \textbf{Invariance:} For $a \in \mathbb{R}^d$ and a bounded Borel subset $A \subset \mathbb{R}^d$
    \[   \widim_\varepsilon\left(\mathcal{X}, \mathbf{d}_{a+A}, \varphi_{a+A}\right) 
         =  \widim_\varepsilon\left(\mathcal{X}, \mathbf{d}_{A}, \varphi_{A}\right),    \] 
     where $a+A  = \{a+u\mid u\in A \}$.    
\end{enumerate}
Notice that we need to assume the nonnegativity of $\varphi$ for the properties (1) and (2).

For a positive number $L$ we denote $\mathbf{d}_{[0,L)^d}$ and $\varphi_{[0, L)^d}$ by 
$\mathbf{d}_L$ and $\varphi_L$ respectively for simplicity.
We define the \textbf{mean dimension with potential} of $\left(\mathcal{X}, T, \varphi\right)$ by 
\begin{equation}  \label{eq: definition of mean dimension with potential}
    \mdim\left(\mathcal{X}, T, \varphi\right) = \lim_{\varepsilon \to 0} 
      \left(\lim_{L\to \infty} \frac{\widim_\varepsilon(\mathcal{X}, \mathbf{d}_L, \varphi_L)}{L^d}\right). 
\end{equation}      
This is a topological invariant, namely its value is independent of the choice of the metric $\mathbf{d}$.
Notice that we do not assume the nonnegativity of $\varphi$ in the definition 
(\ref{eq: definition of mean dimension with potential}).

We need to check that the limits in the definition (\ref{eq: definition of mean dimension with potential}) exist.
The limit with respect to $\varepsilon$ exists because $\widim_\varepsilon(\mathcal{X}, \mathbf{d}_L, \varphi_L)$
is monotone in $\varepsilon$. 
We prove the existence of the limit with respect to $L$ in the next lemma.

\begin{lemma} \label{lemma: existence of limit with respect to L}
    The limit $\displaystyle \lim_{L\to \infty} \widim_\varepsilon(\mathcal{X}, \mathbf{d}_L, \varphi_L)/L^d$ exists in the definition 
    (\ref{eq: definition of mean dimension with potential}).
\end{lemma}

\begin{proof}
Let $c$ be the minimum of $\varphi(x)$ over $x\in \mathcal{X}$ and set $\psi(x)  = \varphi(x) -c$.
Then $\psi$ is a nonnegative function with 
\[   \widim_\varepsilon\left(\mathcal{X}, \mathbf{d}_A, \psi_A\right) 
      = \widim_\varepsilon\left(\mathcal{X}, \mathbf{d}_A, \varphi_A\right)  - c\, \mathbf{m}(A). \]
Set $h(A) = \widim_\varepsilon\left(\mathcal{X}, \mathbf{d}_A, \psi_A\right)$.
It is enough to prove that the limit $\displaystyle \lim_{L\to \infty} h\left([0,L)^d\right)/L^d$ exists. 
For $0<L<R$, let $n= \lfloor R/L \rfloor$ be the integer part of $R/L$.
We have
\[  [0, R)^d \subset \bigcup_{u\in \mathbb{Z}^d\cap [0, n]^d} \left(Lu+ [0, L)^d\right). \]
Since $\psi$ is nonnegative, $h(A)$ satisfies the subadditivity, monotonicity and invariance. Hence
\[   h\left([0,R)^d\right) \leq (n+1)^d\cdot  h\left([0,L)^d\right).  \]
Dividing this by $R^d$ and letting $R\to \infty$, we get 
\[  \limsup_{R\to \infty} \frac{h\left([0,R)^d\right)}{R^d}  \leq \frac{h\left([0,L)^d\right)}{L^d}. \]
Then letting $L\to \infty$ we get
\[   \limsup_{R\to \infty} \frac{h\left([0,R)^d\right)}{R^d}  \leq  \liminf_{L\to \infty} \frac{h\left([0,L)^d\right)}{L^d}. \]
Therefore the limit $\displaystyle \lim_{L\to \infty} h\left([0,L)^d\right)/L^d$ exists. 
\end{proof}

\begin{remark}
By the Ornstein--Weiss quasi-tiling argument (\cite{Ornstein--Weiss}, \cite[\S 1.3.1]{Gromov})
we can also prove that for any F{\o}lner sequence $A_1, A_2, A_3, \dots$ of $\mathbb{R}^d$ the limit 
\[   \lim_{n\to \infty} \frac{\widim_\varepsilon(\mathcal{X}, \mathbf{d}_{A_n}, \varphi_{A_n})}{\mathbf{m}(A_n)} \]
exists and that its value is independent of the choice of a F{\o}lner sequence.
In particular, we can define the mean dimension with potential by 
\[ \mdim\left(\mathcal{X}, T, \varphi\right) = \lim_{\varepsilon \to 0} 
      \left(\lim_{R \to \infty} \frac{\widim_\varepsilon(\mathcal{X}, \mathbf{d}_{B_R}, \varphi_{B_R})}{\mathbf{m}(B_R)}\right), \]
where $B_R = \{u\in \mathbb{R}^d\mid |u| \leq R\}$.
\end{remark}

\subsection{Main result}   \label{subsection: main result}

Let $\mathcal{X}$ be a compact metrizable space.
Recall that we have denoted by $\mathscr{D}(\mathcal{X})$ the set of metrics $\mathbf{d}$ on $\mathcal{X}$ compatible with the given topology.
Let $T\colon \mathbb{R}^d\times \mathcal{X}\to \mathcal{X}$ be a continuous action.
A Borel probability measure $\mu$ on $\mathcal{X}$ is said to be \textbf{$T$-invariant} if $\mu(T^{-u} A) = \mu(A)$ for all 
$u\in \mathbb{R}^d$ and all Borel subsets $A\subset  \mathcal{X}$.
We define $\mathscr{M}^T(\mathcal{X})$ as the set of all $T$-invariant Borel probability measures $\mu$ on $\mathcal{X}$.

Take a metric $\mathbf{d}\in \mathscr{D}(\mathcal{X})$ and a measure $\mu\in \mathscr{M}^T(\mathcal{X})$.
We randomly choose a point $x\in \mathcal{X}$ according to the distribution $\mu$ and consider the orbit 
$\{T^u x\}_{u\in \mathbb{R}^d}$.
For a positive number $\varepsilon$ we define the rate distortion function $R(\mathbf{d}, \mu, \varepsilon)$
as the minimum bits per unit volume for describing $\{T^u x\}_{u\in \mathbb{R}^d}$ with average distortion bounded by $\varepsilon$
with respect to $\mathbf{d}$.
The precise definition of $R(\mathbf{d},\mu,\varepsilon)$ is given in \S \ref{subsection: rate distortion theory}.

We define the \textbf{upper and lower rate distortion dimensions} by 
\[  \urdim(\mathcal{X}, T, \mathbf{d}, \mu) = \limsup_{\varepsilon\to 0} \frac{R(\mathbf{d}, \mu, \varepsilon)}{\log(1/\varepsilon)},\quad 
    \lrdim(\mathcal{X}, T, \mathbf{d}, \mu) = \liminf_{\varepsilon\to 0} \frac{R(\mathbf{d}, \mu, \varepsilon)}{\log(1/\varepsilon)}. \]
The following is the main result of this paper.

\begin{theorem}[Main theorem]  \label{main theorem}
Let $T\colon \mathbb{R}^d\times \mathcal{X}\to \mathcal{X}$ be a continuous action of $\mathbb{R}^d$ on 
a compact metrizable space $\mathcal{X}$. Let $\varphi\colon \mathcal{X}\to \mathbb{R}$ be a continuous function.
Then for any metric $\mathbf{d}\in \mathscr{D}(\mathcal{X})$
\[   \mdim\left(\mathcal{X}, T, \varphi\right) 
      \leq \sup_{\mu\in \mathscr{M}^T(\mathcal{X})}\left(\lrdim(\mathcal{X}, T, \mathbf{d}, \mu) + \int_\mathcal{X} \varphi\, d\mu\right). \]
\end{theorem}

We propose a conjecture:

\begin{conjecture}  \label{conjecture: double variational principle}
In the setting of Theorem \ref{main theorem} we have 
     \begin{equation*}
         \begin{split}
           \mdim(\mathcal{X}, T, \varphi) & 
           = \min_{\mathbf{d}\in \mathscr{D}(\mathcal{X})}
           \left\{\sup_{\mu\in \mathscr{M}^T(\mathcal{X})}\left(\lrdim(\mathcal{X}, T, \mathbf{d}, \mu) 
           + \int_\mathcal{X} \varphi\, d\mu\right)\right\}  \\
           & =  \min_{\mathbf{d}\in \mathscr{D}(\mathcal{X})}
           \left\{\sup_{\mu\in \mathscr{M}^T(\mathcal{X})}\left(\urdim(\mathcal{X}, T, \mathbf{d}, \mu) 
           + \int_\mathcal{X} \varphi\, d\mu\right)\right\}.
         \end{split}
     \end{equation*}
\end{conjecture}
We think that probably we can prove this conjecture if $T\colon \mathbb{R}^d\times \mathcal{X}\to \mathcal{X}$ is a free minimal action.
The proof will be rather lengthy and technically heavy.
We postpone it to Part II of this series of papers.

Along the way to prove Theorem \ref{main theorem}, we will introduce \textit{mean Hausdorff dimension with potential}
and \textit{metric mean dimension with potential} for $\mathbb{R}^d$-actions and establish their basic properties.
In particular we prove:
\begin{enumerate}
   \item  Mean Hausdorff dimension with potential bounds $\mdim\left(\mathcal{X}, T, \varphi\right)$ from above
   (Theorem \ref{theorem: mdim is bounded by mdimh}).
   \item We can construct invariant probability measures which capture the complexity of dynamics expressed by 
   mean Hausdorff dimension with potential (\textit{Dynamical Frostman’s lemma}; Theorem \ref{theorem: dynamical Frostman lemma}).
   \item  Metric mean dimension with potential bounds 
            $\displaystyle \urdim(\mathcal{X}, T, \mathbf{d}, \mu) + \int_\mathcal{X} \varphi\, d\mu$ from above
            (Proposition \ref{proposition: rate distortion dimension and metric mean dimension}).
   \item  Metric mean dimension with potential can be calculated by using only “local” information
             (Theorem \ref{theorem: local formula of metric mean dimension with potential}).
\end{enumerate}
The results (1) and (2) will be used in the proof of Theorem \ref{main theorem}.
The results (3) and (4) are not used in the proof of Theorem \ref{main theorem}.
We plan to use (3) in Part II of this series of papers.
The result (4) may be useful when we study geometric examples of 
\cite{Gromov, Matsuo--Tsukamoto Brody curves, Tsukamoto Brody curves} in a future.

\subsection{Organization of the paper}
In \S \ref{section: mutual information and rate distortion theory} 
we prepare basic definitions and results on mutual information and 
rate distortion theory. 
In \S \ref{section: metric mean dimension with potential and mean Hausdorff dimension with potential}
we introduce mean Hausdorff dimension with potential and metric mean dimension with potential for $\mathbb{R}^d$-actions.
We also state their fundamental properties in \S \ref{section: metric mean dimension with potential and mean Hausdorff dimension with potential}.
The proofs will be given in \S \ref{section: proof of mdim is bounded by mdimh} and 
\S \ref{section: proof of dynamical Frostman lemma}.
Theorem \ref{main theorem} (Main Theorem) follows from the properties
of mean Hausdorff dimension with potential stated in
\S \ref{section: metric mean dimension with potential and mean Hausdorff dimension with potential}.
In \S \ref{section: mean dimension of Z^d-actions} we prepare some basic results on mean dimension theory of 
$\mathbb{Z}^d$-actions. They will be used in \S \ref{section: proof of mdim is bounded by mdimh}.
In \S \ref{section: proof of mdim is bounded by mdimh} we prove that $\mdim(\mathcal{X}, T, \varphi)$ is bounded from above by 
mean Hausdorff dimension with potential.
In \S \ref{section: proof of dynamical Frostman lemma} we prove dynamical Frostman’s lemma.
In \S \ref{section: local nature of metric mean dimension with potential} we prove that metric mean dimension with potential 
can be calculated by using certain local information.
\S \ref{section: local nature of metric mean dimension with potential} is independent of 
the proof of Theorem \ref{main theorem}.

\section{Mutual information and rate distortion theory}  \label{section: mutual information and rate distortion theory}

We prepare basics of rate distortion theory in this section.
Throughout this paper $\log x$ denotes the logarithm of base two.
The natural logarithm is denoted by $\ln x$:
\[   \log x = \log_2 x, \quad \ln x = \log_e x. \]
This section is rather long.
This is partly because we have to be careful of measure theoretic details.
Hopefully this section will become a useful reference in a future study of mean dimension of $\mathbb{R}^d$-actions.
At the first reading, readers may skip the whole of Subsection \ref{subsection: measure theoretic preparations} and most of 
Subsection \ref{subsection: mutual information}.
The crucial parts of this section are only the definition of mutual information in \S \ref{subsection: mutual information} and 
the definition of rate distortion function in \S \ref{subsection: rate distortion theory}.
All the rest of this section are technical details.

\subsection{Measure theoretic preparations}   \label{subsection: measure theoretic preparations}

We need to prepare some basic results on measure theory.
A \textbf{measurable space} is a pair $(\mathcal{X}, \mathcal{A})$ of a set $\mathcal{X}$ and its $\sigma$-algebra $\mathcal{A}$. 
Two measurable spaces $(\mathcal{X}, \mathcal{A})$ and $(\mathcal{Y}, \mathcal{B})$ are said to be \textbf{isomorphic} if 
there exists a bijection $f\colon \mathcal{X}\to \mathcal{Y}$ such that both $f$ and $f^{-1}$ are measurable 
(i.e. $f(\mathcal{A}) = \mathcal{B}$).

For a topological space $\mathcal{X}$, its \textbf{Borel $\sigma$-algebra} $\mathcal{B}_\mathcal{X}$ 
is the minimum $\sigma$-algebra containing all open subsets of $\mathcal{X}$. 
A \textbf{Polish space} is a topological space $\mathcal{X}$ admitting a metric $\mathbf{d}$ for which $(\mathcal{X}, \mathbf{d})$
is a complete separable metric space. 

A measurable space $(\mathcal{X}, \mathcal{A})$ is said to be a \textbf{standard Borel space} if there exists a Polish space $\mathcal{Y}$
for which $(\mathcal{X}, \mathcal{A})$ is isomorphic to $(\mathcal{Y}, \mathcal{B}_\mathcal{Y})$ as measurable spaces.
It is known that any two uncountable standard Borel spaces are isomorphic to each other
(the Borel isomorphism theorem \cite[Theorem 3.3.13]{Srivastava}).
Therefore every standard Borel space is isomorphic to one of the following measurable spaces:
  \begin{itemize}
     \item  A finite set $A$ with its discrete $\sigma$-algebra $2^A := \{\text{subset of $A$}\}$.
     \item  The set of natural numbers $\mathbb{N}$ with its discrete $\sigma$-algebra 
               $2^{\mathbb{N}} := \{\text{subset of $\mathbb{N}$}\}$.
     \item  The Cantor set $\mathcal{C} = \{0,1\}^{\mathbb{N}}$ with its Borel $\sigma$-algebra $\mathcal{B}_{\mathcal{C}}$.
     (Here $\{0,1\}$ is endowed with the discrete topology and the topology of $\mathcal{C}$ is the product topology.) 
  \end{itemize}  

An importance of standard Borel spaces is that we can prove the existence of \textit{regular conditional distribution} 
under the assumption of “standard Borel”. 
Let $(\mathcal{X}, \mathcal{A})$ and $(\mathcal{Y}, \mathcal{B})$ be measurable spaces.
A \textbf{transition probability} on $\mathcal{X}\times \mathcal{Y}$ is a map 
$\nu:\mathcal{X}\times  \mathcal{B}\to [0,1]$ such that 
 \begin{itemize}
    \item for every $x\in \mathcal{X}$, the map $\mathcal{B}\ni B\mapsto \nu(x, B)\in [0,1]$ is a probability measure on 
    $(\mathcal{Y}, \mathcal{B})$,   
    \item for every $B \in \mathcal{B}$, the map $\mathcal{X}\ni x \mapsto \nu(x, B)\in [0,1]$ is measurable.
 \end{itemize}
We often denote $\nu(x, B)$ by $\nu(B|x)$.

For two measurable spaces $(\mathcal{X}, \mathcal{A})$ and $(\mathcal{Y}, \mathcal{B})$ we denote their product by 
$(\mathcal{X}\times \mathcal{Y}, \mathcal{A}\otimes \mathcal{B})$ where $\mathcal{A}\otimes \mathcal{B}$ is the minimum 
$\sigma$-algebra containing all the rectangles $A\times B$ $(A\in \mathcal{A}, B\in \mathcal{B})$.
For any $E\in \mathcal{A}\otimes \mathcal{B}$, it is known that the \textbf{section} $E_x:= \{y\in \mathcal{Y}\mid (x, y)\in E\}$
belongs to $\mathcal{B}$ for every $x\in \mathcal{X}$.
(This fact is a part of the Fubini theorem. It can be easily proved by using Dynkin’s $\pi$-$\lambda$ theorem
\cite[p.402 Theorem A.1.4]{Durrett}.)
Moreover, if $(\mathcal{Y}, \mathcal{B})$ is a standard Borel space, 
then for any transition probability $\nu$ on $\mathcal{X}\times \mathcal{Y}$
and any $E\in \mathcal{A}\otimes \mathcal{B}$ the map 
$\mathcal{X}\ni x\mapsto \nu(E_x|x) \in [0,1]$ is measurable \cite[Proposition 3.4.24]{Srivastava}.

A \textbf{probability space} is a triplet $(\Omega, \mathcal{F}, \mathbb{P})$ where $(\Omega, \mathcal{F})$ is a measurable space and 
$\mathbb{P}$ is a probability measure defined on it.
Let $X\colon \Omega \to \mathcal{X}$ be a measurable map from a probability space $(\Omega, \mathcal{F}, \mathbb{P})$ to a
Borel space $(\mathcal{X}, \mathcal{A})$.
We denote the push-forward measure $X_*\mathbb{P}$ by $\mathrm{Law} X$ and call it the \textbf{law of $X$} 
or the \textbf{distribution of $X$}.
(Here $X_*\mathbb{P}(A) = \mathbb{P}\left(X\in A\right) = \mathbb{P}\left(X^{-1}(A)\right)$ for $A\in \mathcal{A}$.)

The next theorem is a fundamental result.
It guarantees the existence of regular conditional probability.
For the proof, see \cite[p.15 Theorem 3.3 and its Corollary]{Ikeda--Watanabe} or 
Gray \cite[p. 182 Corollary 6.2]{Gray_probability}.

\begin{theorem}[Existence of regular conditional distribution]
Let $(\Omega, \mathcal{F}, \mathbb{P})$ be a probability space, and 
$(\mathcal{X}, \mathcal{A})$ and $(\mathcal{Y}, \mathcal{B})$ standard Borel spaces.
Let $X:\Omega\to \mathcal{X}$ and $Y:\Omega \to \mathcal{Y}$ be measurable maps,
and set $\mu := \mathrm{Law} X$. 
Then there exists a transition probability $\nu$ on $\mathcal{X}\times \mathcal{Y}$ such that for any 
$E\in \mathcal{A}\otimes \mathcal{B}$ we have 
\[  \mathbb{P}\left((X,Y)\in E\right) = \int_\mathcal{X} \nu(E_x|x)\, d\mu(x).  \]
If a transition probability $\nu^\prime$ on $\mathcal{X}\times \mathcal{Y}$ satisfies the same property then 
there exists a $\mu$-null set $N\in \mathcal{A}$ such that 
$\nu(B|x) = \nu^\prime(B|x)$ for all $x\in \mathcal{X}\setminus N$ and $B\in \mathcal{B}$.
\end{theorem}

The transition probability $\nu(\cdot|x)$ in this theorem is called the \textbf{regular conditional distribution of $Y$ given $X=x$}.
We sometimes denote $\nu(B|x)$ by
$\mathbb{P}(Y\in B|X=x)$ for $x\in \mathcal{X}$ and $B\in \mathcal{B}$.
If $\mathcal{X}$ and $\mathcal{Y}$ are finite sets, then this coincides with the elementary notion of conditional probability:
\[   \mathbb{P}(Y\in B|X=x)  =  \frac{\mathbb{P}(X=x, Y\in B)}{\mathbb{P}(X=x)}, \quad \left(\text{if }\mathbb{P}(X=x)\neq 0\right). \]
In this case we usually denote $\nu\left(\{y\}|x\right)$ by $\nu(y|x)$ $(x\in \mathcal{X}, y\in \mathcal{Y})$ and call it 
a \textbf{conditional probability mass function}\footnote{For convenience in the sequel, we define this notion more precisely.
Let $\mathcal{X}$ and $\mathcal{Y}$ be finite sets. 
A map $\mathcal{X}\times \mathcal{Y}\ni (x, y) \mapsto  \nu(y|x)\in [0,1]$ is called a conditional probability mass function if 
$\sum_{y\in \mathcal{Y}} \nu(y|x) = 1$ for every $x\in \mathcal{X}$.}.

By using the notion of regular conditional distribution, 
we can introduce the definition of \textit{conditional independence} of random variables.
Let $(\Omega, \mathbb{P})$ be a probability space and $(\mathcal{X},\mathcal{A}), (\mathcal{Y},\mathcal{B}), (\mathcal{Z},\mathcal{C})$
standard Borel spaces.
Let $X\colon \Omega\to \mathcal{X}$, $Y\colon \Omega\to \mathcal{Y}$ and $Z\colon \Omega\to \mathcal{Z}$ be measurable maps.
We say that $X$ and $Y$ are \textbf{conditionally independent given $Z$} if we have 
\begin{equation} \label{eq: conditional independence}
    \mathbb{P}\left((X, Y)\in A\times B| Z=z\right) = \mathbb{P}\left(X\in A|Z=z\right) \cdot \mathbb{P}\left(Y\in B|Z=z\right) 
\end{equation}    
for $Z_*\mathbb{P}$-a.e. $z\in \mathcal{Z}$ and all $A\in \mathcal{A}$ and $B\in \mathcal{B}$.
Here $Z_*\mathbb{P}$ is the push-forward measure of $\mathbb{P}$ by $Z$.
The left-hand side of (\ref{eq: conditional independence})
is the conditional regular distribution of $(X, Y)\colon \Omega \to \mathcal{X}\times \mathcal{Y}$ given $Z=z$.
The right-hand side is the multiple of the conditional distribution of $X$ given $Z=z$ and the conditional distribution of $Y$
given $Z=z$.

At the end of this subsection we explain the log-sum inequality.
This will be used in the next subsection.

\begin{lemma}[Log-sum inequality] \label{lemma: log-sum inequality}
Let $(\mathcal{X}, \mathcal{A})$ be a measurable space.
Let $\mu$ be a measure on it with $0< \mu(\mathcal{X}) < \infty$.
Let $f$ and $g$ be nonnegative measurable functions defined on $\mathcal{X}$.
Suppose that $g$ is $\mu$-integrable and $g(x) >0$ for $\mu$-a.e. $x\in \mathcal{X}$.
Then 
\[  \left(\int_\mathcal{X} f(x) \, d\mu(x) \right) \log \frac{\int_\mathcal{X} f(x) \, d\mu(x)}{\int_\mathcal{X} g(x) \, d\mu(x)}
     \leq \int_\mathcal{X} f(x) \log \frac{f(x)}{g(x)}\, d\mu(x). \]
In particular, if the left-hand side is infinite then the right-hand side is also infinite.     
\end{lemma}

Here we assume $0\log \frac{0}{a} = 0$ for all $a>0$.

\begin{proof}
Set $\phi(t) = t\log t$ for $t\geq 0$. Since $\phi^{\prime \prime}(t) = \log e/t>0$ for $t>0$, this is a convex function.
We define a probability measure $w$ on $\mathcal{X}$ by 
\[  w(A) = \frac{\int_A g \, d\mu}{\int_{\mathcal{X}} g\, d\mu} \quad (A\subset \mathcal{X}). \]
The Radon--Nikodim derivative of $w$ by $\mu$ is given by
\[  \frac{dw}{d\mu} = \frac{g}{\int_{\mathcal{X}} g\, d\mu}. \]
By Jensen’s inequality
\begin{equation}   \label{eq: Jensen inequality}
    \phi\left(\int_{\mathcal{X}} \frac{f}{g}\, dw\right) \leq \int_{\mathcal{X}} \phi\left(\frac{f}{g}\right) \, dw. 
\end{equation}
Here, if the left-hand side is infinite, then the right-hand side is also infinite.
We have 
\[  \phi\left(\int_{\mathcal{X}} \frac{f}{g}\, dw\right)  = \phi\left(\frac{\int_{\mathcal{X}}f\, d\mu}{\int_{\mathcal{X}}g\, d\mu}\right)
     = \frac{\int_{\mathcal{X}}f\, d\mu}{\int_{\mathcal{X}}g\, d\mu} \log \frac{\int_{\mathcal{X}}f\, d\mu}{\int_{\mathcal{X}}g\, d\mu}. \]
The right-hand side of (\ref{eq: Jensen inequality}) is 
\[      \int_{\mathcal{X}} \phi\left(\frac{f}{g}\right) \, dw 
     =  \frac{1}{\int_{\mathcal{X}} g\, d\mu} \int_{\mathcal{X}} g \phi\left(\frac{f}{g}\right)\, d\mu 
     =  \frac{1}{\int_{\mathcal{X}} g\, d\mu} \int_{\mathcal{X}} f \log \frac{f}{g}\, d\mu. \]
Therefore (\ref{eq: Jensen inequality}) provides 
\[   \left(\int_{\mathcal{X}} f\, d\mu\right) \log \frac{\int_{\mathcal{X}} f\, d\mu}{\int_{\mathcal{X}}g\, d\mu}
      \leq  \int_{\mathcal{X}} f \log \frac{f}{g}\, d\mu. \]     
\end{proof}

The following is the finitary version of the log-sum inequality:

\begin{corollary}   \label{cor: log-sum inequality}
Let $a_1, \dots, a_n$ be nonnegative numbers and $b_1, \dots, b_n$ positive numbers.
Then 
\[   \left(\sum_{i=1}^n a_i\right) \log \frac{\sum_{i=1}^n a_i}{\sum_{i=1}^n b_i} \leq 
       \sum_{i=1}^n a_i \log \frac{a_i}{b_i}. \]
\end{corollary}

\begin{proof}
Apply Lemma \ref{lemma: log-sum inequality} to the finite set $\mathcal{X}= \{1,2,\dots, n\}$ with the discrete $\sigma$-algebra and 
the counting measure.
\end{proof}

\subsection{Mutual information}  \label{subsection: mutual information}

Let $(\Omega, \mathcal{F}, \mathbb{P})$ be a probability space.
We assume that all random variables in this subsection are defined on $(\Omega, \mathcal{F}, \mathbb{P})$.
In this paper a finite set is always assumed to be endowed with the discrete topology and the discrete $\sigma$-algebra (i.e. the set of 
all subsets).
The purpose of this subsection is to define and study mutual information.
A basic reference of mutual information is the book of Cover--Thomas \cite{Cover--Thomas}.
A mathematically sophisticated presentation is given in the book of Gray \cite{Gray_entropy}.

First we define the Shannon entropy.
Let $(\mathcal{X}, \mathcal{A})$ be a finite set with the discrete $\sigma$-algebra, and let 
$X\colon \Omega \to \mathcal{X}$ be a measurable map.
We define the \textbf{Shannon entropy of $X$} by 
\[  H(X) = -\sum_{x\in \mathcal{X}} \mathbb{P}(X=x) \log \mathbb{P}(X=x). \]
Here we assume $0\log 0 = 0$ as usual.

Next we define the mutual information.
Let $(\mathcal{X}, \mathcal{A})$ and $(\mathcal{Y}, \mathcal{B})$ be two measurable spaces and 
let $X\colon \Omega \to \mathcal{X}$ and $Y\colon \Omega \to \mathcal{Y}$ be measurable maps.
We want to define the mutual information $I(X;Y)$.
Intuitively $I(X;Y)$ measure the amount of information shared by the random variables $X$ and $Y$.

\begin{itemize}
  \item \textbf{Case I:} 
   Suppose $(\mathcal{X}, \mathcal{A})$ and $(\mathcal{Y}, \mathcal{B})$ are finite sets with 
   the discrete $\sigma$-algebras.
   Then we define 
    \begin{equation}  \label{eq: mutual information}
        I(X;Y) = H(X) + H(Y) - H(X, Y), 
    \end{equation}    
    where $H(X, Y)$ is the Shannon entropy of the measurable map $(X, Y): \Omega \to \mathcal{X}\times \mathcal{Y}$.
    Since $H(X, Y) \leq  H(X) + H(Y)$, the mutual information $I(X;Y)$ is always nonnegative.
    The explicit formula is given by
   \[  I(X;Y) =
   \sum_{x\in \mathcal{X}, y\in \mathcal{Y}}  \mathbb{P}(X=x, Y=y) \log \frac{\mathbb{P}(X=x, Y=y)}{\mathbb{P}(X=x) \mathbb{P}(Y=y)}. \]
    Here we assume $0\log \frac{0}{a}=0$ for any $a\geq 0$.   
    The mutual information $I(X;Y)$ satisfies 
    the following natural monotonicity\footnote{This is a special case of 
    the data-processing inequality \cite[Theorem 2.8.1]{Cover--Thomas}.}:
   Let $\mathcal{X}^\prime$ and $\mathcal{Y}^\prime$ be finite sets (endowed with the discrete $\sigma$-algebras), 
   and let $f\colon \mathcal{X}\to \mathcal{X}^\prime$ and $g\colon \mathcal{Y}\to \mathcal{Y}^\prime$ be any maps. 
   Then it follows from the log-sum inequality (Corollary \ref{cor: log-sum inequality}) that
   \begin{equation} \label{eq: data-processing inequality}
       I\left(f(X); g(Y)\right) \leq I(X;Y). 
   \end{equation}    
   \item \textbf{Case II:} Here we define $I(X;Y)$ for general random variables $X$ and $Y$.
   (Namely $\mathcal{X}$ and $\mathcal{Y}$ may be infinite sets.)
   Let $\mathcal{X}^\prime$ and $\mathcal{Y}^\prime$ be arbitrary finite sets, and let $f\colon \mathcal{X}\to \mathcal{X}^\prime$ and 
   $g\colon \mathcal{Y}\to \mathcal{Y}^\prime$ be any measurable maps.
   Then we can consider $I\left(f(X); g(Y)\right)$ by Case I.
   We define $I(X;Y)$ as the supremum of $I\left(f(X); g(Y)\right)$ over all finite sets $\mathcal{X}^\prime, \mathcal{Y}^\prime$ and 
   measurable maps $f\colon \mathcal{X}\to \mathcal{X}^\prime$ and $g\colon \mathcal{Y}\to \mathcal{Y}^\prime$.
   The mutual information $I(X, Y)$ is always nonnegative and symmetric: $I(X;Y) = I(Y;X)\geq 0$.
 
   If $f\colon \mathcal{X}\to \mathcal{X}^\prime$ and $g\colon \mathcal{Y}\to \mathcal{Y}^\prime$ are
   measurable maps to some other measurable spaces $(\mathcal{X}^\prime, \mathcal{A}^\prime)$ and 
   $(\mathcal{Y}^\prime, \mathcal{B}^\prime)$ (not necessarily finite sets) then we have $I\left(f(X); g(Y)\right)  \leq  I\left(X;Y\right)$.
   
   If $\mathcal{X}$ and $\mathcal{Y}$ are finite sets, then the definition of Case II is compatible 
   with Case I by the monotonicity (\ref{eq: data-processing inequality}).
\end{itemize}

If $(\mathcal{X}, \mathcal{A})$ and $(\mathcal{Y}, \mathcal{B})$ are standard Borel spaces, then we can consider 
the regular conditional distribution of $Y$ given $X=x$. We denote 
\[  \nu(B|x) = \mathbb{P}(Y\in B|X=x) \quad (x\in \mathcal{X}, B\in \mathcal{B}). \]
Let $\mu = X_*\mathbb{P}$ be the push-forward measure of $\mathbb{P}$ by $X$. 
The distribution of $(X, Y)$ is determined by $\mu$ and $\nu$.
Hence the mutual information $I(X;Y)$ is also determined by $\mu$ and $\nu$.
Therefore we sometimes denote $I(X;Y)$ by $I(\mu, \nu)$.
An importance of this description comes from the fact that $I(\mu, \nu)$ is a concave function in $\mu$ and
a convex function in $\nu$ (Proposition \ref{proposition: concavity and convexity of I(mu, nu)} below).

In the rest of this subsection we prepare several basic properties of mutual information.
They are rather heavy.
Readers may skip to the next subsection at the first reading.

If $(\mathcal{X}, \mathcal{A})$ is a standard Borel space and \textit{if $\mathcal{Y}$ is a finite set}, then
we can express $I(X; Y)$ in another convenient way.
For $x\in \mathcal{X}$ we set
\[  H(Y|X=x) = - \sum_{y\in \mathcal{Y}} \mathbb{P}\left(Y=y|X=x\right) \log \mathbb{P}\left(Y=y|X=x\right). \]
We define the \textbf{conditional entropy of $Y$ given $X$} by 
\[  H(Y|X) = \int_{\mathcal{X}} H(Y|X=x) \, d\mu(x),  \quad (\mu := X_*\mathbb{P}). \]
The next theorem is given in the book of Gray \cite[p. 213, Lemma 7.20]{Gray_entropy}.

\begin{theorem} \label{theorem: mutual information and conditional entropy}
Let $X$ and $Y$ be random variables taking values in a standard Borel space $(\mathcal{X}, \mathcal{A})$ and a 
finite set $\mathcal{Y}$ respectively. Then we have 
\[   I(X;Y) = H(Y) - H(Y|X). \]
\end{theorem}
When both $\mathcal{X}$ and $\mathcal{Y}$ are finite sets, this theorem is a very well-known result.
A point of the theorem is that we do not need to assume that $\mathcal{X}$ is a finite set.

The following is also a basic result. This is given in \cite[p. 211, Lemma 7.18]{Gray_entropy}.

\begin{theorem}  \label{theorem: approximation of mutual information}
Let $(\mathcal{X}, \mathcal{A})$ and $\left(\mathcal{Y}, \mathcal{B}\right)$ be standard Borel spaces.
Then there exist sequences of measurable maps $f_n\colon \mathcal{X}\to \mathcal{X}_n$ and 
$g_n\colon \mathcal{Y}\to \mathcal{Y}_n$ to some finite sets $\mathcal{X}_n$ and $\mathcal{Y}_n$ $(n\geq 1)$ for which 
the following statement holds:
If $X$ and $Y$ are random variables taking values in $(\mathcal{X}, \mathcal{A})$ and $\left(\mathcal{Y}, \mathcal{B}\right)$
respectively, then 
\[   I(X;Y) = \lim_{n\to \infty}  I\left(f_n(X); g_n(Y)\right). \]
\end{theorem}

\begin{proof}[Sketch of the proof]
Standard Borel spaces are isomorphic to either countable sets or the Cantor sets.
The case of countable sets is easier. So we assume that both $\mathcal{X}$ and $\mathcal{Y}$ are the 
Cantor set $\{0,1\}^{\mathbb{N}}$.
Let $f_n\colon \mathcal{X}\to \{0,1\}^n$ and $g_n\colon \mathcal{Y}\to \{0,1\}^n$ be the natural projections to the first 
$n$ coordinates. Then we can check that $f_n$ and $g_n$ satisfy the statement.
\end{proof}

\begin{lemma}  \label{lemma: convergence of mutual information}
Let $X_n$ and $Y_n$ $(n\geq 1)$ be sequences of random variables taking values in finite sets 
$\mathcal{X}$ and $\mathcal{Y}$ respectively.
Suppose $(X_n, Y_n)$ converges to $(X, Y)$ in law.
(Namely $\mathbb{P}\left(X_n = x, Y_n =y\right) \to \mathbb{P}(X=x, Y=y)$ as $n\to \infty$ 
for all $(x, y)\in \mathcal{X}\times \mathcal{Y}$.)
Then $I(X_n;Y_n)$ converges to $I(X;Y)$ as $n\to \infty$.
\end{lemma}

\begin{proof}
This immediately follows from the definition (\ref{eq: mutual information}) in Case I above.
\end{proof}

\begin{lemma}[Subadditivity of mutual information]  \label{lemma: subadditivity of mutual information}
Let $X, Y, Z$ be random variables taking values in standard Borel spaces $(\mathcal{X}, \mathcal{A}), 
(\mathcal{Y},\mathcal{B}), (\mathcal{Z}, \mathcal{C})$ respectively. Suppose that $X$ and $Y$ are conditionally independent 
given $Z$. Then 
\[   I\left(X, Y;Z\right)  \leq  I(X; Z) + I(Y; Z), \]
where $I\left(X, Y; Z\right) = I\left((X, Y); Z\right)$ is the mutual information between the random variables $(X, Y)$ and $Z$.
\end{lemma}

\begin{proof}
Let $f\colon \mathcal{X}\to \mathcal{X}^\prime$ and $g\colon \mathcal{Y}\to \mathcal{Y}^\prime$
be measurable maps to some finite sets $\mathcal{X}^\prime$ and $\mathcal{Y}^\prime$.
Then by Theorem \ref{theorem: mutual information and conditional entropy}
\[  I\left(f(X), g(Y); Z\right) = H\left(f(X), g(Y)\right) - H\left(f(X), g(Y)|Z\right). \]
We have \cite[Theorem 2.6.6]{Cover--Thomas}
\[  H\left(f(X), g(Y)\right)  \leq H\left(f(X)\right) + H\left(g(Y)\right). \]
The random variables $f(X)$ and $g(Y)$ are conditionally independent given $Z$. Hence 
\[ H\left(f(X), g(Y)|Z\right)  = H\left(f(X)|Z\right) + H\left(g(Y)|Z\right). \]
Therefore 
\begin{equation*}
   \begin{split}
      I\left(f(X), g(Y); Z\right) & \leq \left\{H\left(f(X)\right) - H\left(f(X)|Z\right)\right\} + 
                                         \left\{H\left(g(Y)\right) - H\left(g(Y)|Z\right)\right\} \\
                                    & = I\left(f(X); Z\right) + I\left(g(Y); Z\right). 
   \end{split}
\end{equation*}    
                                     
We have $I\left(X, Y; Z\right) = \sup_{f, g} I\left(f(X), g(Y); Z\right)$ where $f\colon \mathcal{X}\to \mathcal{X}^\prime$ and 
$g\colon \mathcal{Y}\to \mathcal{Y}^\prime$ run over all measurable maps to some finite sets.
This follows from the fact that $\mathcal{A}\otimes \mathcal{B}$ is generated by rectangles
$A\times B$ $(A\in \mathcal{A}, B\in \mathcal{B})$ \cite[p.175 Lemma 7.3]{Gray_entropy}.
Therefore we get 
\[    I\left(X, Y;Z\right)  \leq  I(X; Z) + I(Y; Z).  \]
\end{proof}

As we briefly mentioned above, the mutual information $I(\mu, \nu)$ is a concave function in a probability measure $\mu$
and a convex function in a transition probability $\nu$.
Next we are going to establish this fact. We need some preparations.

For a finite set $\mathcal{Y}$, a \textbf{probability mass function} $p$ on $\mathcal{Y}$ is a nonnegative function on $\mathcal{Y}$
satisfying $\sum_{y\in \mathcal{Y}} p(y) = 1$. For a probability mass function $p$ on $\mathcal{Y}$ we define 
\[  H(p) = -\sum_{y\in \mathcal{Y}} p(y) \log p(y). \]

\begin{lemma}[Concavity of the Shannon entropy] \label{lemma: entropy is concave}
 Let $\mathcal{Y}$ be a finite set and let $(\mathcal{Z}, \mathcal{C}, m)$ be a probability space.
 Suppose that we are given a probability mass function $p_z$ on $\mathcal{Y}$ for each $z\in \mathcal{Z}$
 and that the map $\mathcal{Z}\ni z\mapsto p_z(y)\in [0,1]$ is measurable for each $y\in \mathcal{Y}$.
 We define a probability mass function $p$ on $\mathcal{Y}$ by 
 \[  p(y) = \int_{\mathcal{Z}} p_z(y) \, dm(z). \]
 Then
 \[   H(p) \geq   \int_{\mathcal{Z}}  H(p_z) \, dm(z). \]
\end{lemma}

\begin{proof}
From the log-sum inequality (Lemma \ref{lemma: log-sum inequality}), 
\[   - p(y) \log p(y) \geq - \int_{\mathcal{Z}} p_z(y) \log p_z(y)\, dm(z). \]
Summing this over $y\in \mathcal{Y}$, we get the statement.
\end{proof}

\begin{lemma}  \label{lemma: mutual information is convex in nu finite case}
Let $\mathcal{X}$ and $\mathcal{Y}$ be finite sets and $(\mathcal{Z}, \mathcal{C}, m)$ a probability space.
Let $\mu$ be a probability mass function on $\mathcal{X}$.
Suppose that, for each $z\in \mathcal{Z}$, we are given a conditional probability mass function
$\nu_z(y|x)$ in $x\in \mathcal{X}$ and $y\in \mathcal{Y}$ such that the map 
$\mathcal{Z}\ni  z\mapsto \nu_z(y|x)\in [0,1]$ is measurable for each $(x, y)\in \mathcal{X}\times \mathcal{Y}$.
We define 
\[   \nu(y|x) = \int_{\mathcal{Z}} \nu_z(y|x)\, dm(z), \quad (x\in \mathcal{X}, y\in \mathcal{Y}).  \]
Then 
\[   I(\mu, \nu)  \leq \int_{\mathcal{Z}} I(\mu, \nu_z)\, dm(z). \]
\end{lemma}

\begin{proof}
For $y\in \mathcal{Y}$ we set 
\[  p_z(y) = \sum_{x\in \mathcal{X}} \mu(x) \nu_z(y|x), \quad p(y) = \sum_{x\in \mathcal{X}} \mu(x) \nu(y|x). \]
We have 
\[   I(\mu, \nu) = \sum_{x, y} \mu(x) \nu(y|x) \log \frac{\mu(x) \nu(y|x)}{\mu(x) p(y)}, \quad 
      I(\mu, \nu_z) = \sum_{x, y} \mu(x) \nu_z(y|x) \log \frac{\mu(x)\nu_z(y|x)}{\mu(x) p_z(y)}.  \]
Here we assume $0\log \frac{a}{0} = 0$ for all $a\geq 0$.

We estimate each summand of $I(\mu, \nu)$ and $I(\mu, \nu_z)$.
We fix $(x, y)\in \mathcal{X}\times \mathcal{Y}$ with $\mu(x) p(y) >0$.
We define a subset $\mathcal{Z}^\prime \subset \mathcal{Z}$ by 
\[  \mathcal{Z}^\prime = \{z\mid  p_z(y) >0\}  \supset \{z\mid   \nu_z(y|x)>0\}. \]
Since $\mu(x) p(y) >0$, we have $m\left(\mathcal{Z}^\prime\right) >0$.
We have 
\[  \mu(x) \nu(y|x) = \int_{\mathcal{Z}^\prime} \mu(x)\nu_z(y|x)\, dm(z), \quad 
     \mu(x)p(y)  = \int_{\mathcal{Z}^\prime} \mu(x) p_z(y) \, dm(z). \]
By the log-sum inequality (Lemma \ref{lemma: log-sum inequality})
\begin{equation*}
   \begin{split}
     \mu(x) \nu(y|x) \log \frac{\mu(x) \nu(y|x)}{\mu(x) p(y)}  & \leq 
     \int_{\mathcal{Z}^\prime} \mu(x) \nu_z(y|x) \log \frac{\mu(x)\nu_z(y|x)}{\mu(x) p_z(y)} dm(z)  \\
     & = \int_{\mathcal{Z}} \mu(x) \nu_z(y|x) \log \frac{\mu(x)\nu_z(y|x)}{\mu(x) p_z(y)} dm(z).
    \end{split}
\end{equation*}     
Taking sums over $(x, y)\in \mathcal{X}\times \mathcal{Y}$, we get the statement.
\end{proof}

\begin{proposition}[$I(\mu, \nu)$ is concaive in $\mu$ and convex in $\nu$] \label{proposition: concavity and convexity of I(mu, nu)}
Let $(\mathcal{X}, \mathcal{A})$ and $(\mathcal{Y}, \mathcal{B})$ be standard Borel spaces, and let 
$(\mathcal{Z}, \mathcal{C}, m)$ be a probability space.

\begin{enumerate}
   \item Let $\nu$ be a transition probability on $\mathcal{X}\times \mathcal{Y}$.
   Suppose that we are given a probability measure $\mu_z$ on $\mathcal{X}$ for each $z\in \mathcal{Z}$ such that 
   the map $\mathcal{Z}\ni z\mapsto  \mu_z(A)\in [0,1]$ is measurable for every $A\in \mathcal{A}$.
   We define a probability measure $\mu$ on $(\mathcal{X}, \mathcal{A})$ by 
\[  \mu(A)  = \int_{\mathcal{Z}} \mu_z(A)\, dm(z), \quad (A\in \mathcal{A}). \]
   Then we have 
  \begin{equation}   \label{eq: concavity of mutual information}
     I(\mu, \nu) \geq \int_{\mathcal{Z}} I(\mu_z, \nu)\, dm(z). 
  \end{equation}  
   \item  Let $\mu$ be a probability measure on $\mathcal{X}$. Suppose that we are given a 
   transition probability $\nu_z$ on $\mathcal{X}\times \mathcal{Y}$ for each $z\in \mathcal{Z}$
   such that 
   the map $\mathcal{X} \times \mathcal{Z} \ni (x, z) \mapsto  \nu_z(B|x) \in [0,1]$ is 
   measurable with respect to $\mathcal{A}\otimes \mathcal{C}$ for each $B \in \mathcal{B}$.
   We define a transition probability $\nu$ on $\mathcal{X}\times \mathcal{Y}$ by
   \[  \nu(B|x) = \int_{\mathcal{Z}} \nu_z(B|x) dm(z), \quad (x\in \mathcal{X}, B\in \mathcal{B}). \]
   Then we have 
\[  I(\mu, \nu)  \leq \int_{\mathcal{Z}} I(\mu, \nu_z) \, dm(z). \]
\end{enumerate}
\end{proposition}

\begin{proof}
(1) By Lemma \ref{theorem: approximation of mutual information}, there exists a sequence of measurable maps 
$g_n\colon \mathcal{Y}\to \mathcal{Y}_n$ to finite sets $\mathcal{Y}_n$ such that 
\[  I(\mu_z, \nu) = \lim_{n\to \infty} I\left(\mu_z, (g_n)_*\nu\right), \quad 
     I(\mu, \nu) = \lim_{n\to \infty} I\left(\mu, (g_n)_*\nu\right).   \]
Here $(g_n)_*\nu$ is a transition probability on $\mathcal{X}\times \mathcal{Y}_n$ defined by 
\[  (g_n)_*\nu(B|x) = \nu\left((g_n)^{-1}B|x\right), \quad (B\subset \mathcal{Y}_n). \]
It is enough to prove that for each $n$ we have 
\[   I\left(\mu, (g_n)_*\nu\right) \geq  \int_{\mathcal{Z}} I\left(\mu_z, (g_n)_*\nu\right) \, dm(z). \]
If this is proved then we get the above (\ref{eq: concavity of mutual information}) by Fatou’s lemma.
Therefore we can assume that $\mathcal{Y}$ itself is a finite set from the beginning.

We define probability mass functions $p(y)$ and $p_z(y)$ $(z\in \mathcal{Z}$) on $\mathcal{Y}$ by 
\[  p(y) = \int_{\mathcal{X}} \nu(y|x)\, d\mu(x), \quad   p_z(y) = \int_{\mathcal{X}} \nu(y|x)\, d\mu_z(x). \]
We have $p(y) = \int_{\mathcal{Z}} p_z(y)\, dm(z)$.
Then by Theorem \ref{theorem: mutual information and conditional entropy}
\[  I(\mu, \nu) = H(p) - \int_{\mathcal{X}} H\left(\nu(\cdot|x)\right)\, d\mu(x), \quad 
     I\left(\mu_z, \nu\right) = H(p_z) - \int_{\mathcal{X}} H\left(\nu(\cdot|x)\right) \, d\mu_z(x). \]
	Here $H\left(\nu(\cdot|x)\right)  = -\sum_{y\in \mathcal{Y}} \nu(y|x)\log \nu(y|x)$. 
	Notice that in particular this shows that 
	$I\left(\mu_z, \nu\right)$ is a measurable function in the variable $z\in \mathcal{Z}$.
By Lemma \ref{lemma: entropy is concave}, we have $H(p) \geq \int_{\mathcal{Z}} H(p_z)\, dm(z)$.
We also have 
\[  \int_{\mathcal{X}} H\left(\nu(\cdot|x)\right) d\mu(x) 
    = \int_{\mathcal{Z}} \left(\int_{\mathcal{X}} H\left(\nu(\cdot|x)\right)\, d\mu_z(x)\right) dm(z). \]
Thus 
\[   I(\mu, \nu)    \geq \int_{\mathcal{Z}}  I\left(\mu_z, \nu\right) dm(z). \]

(2) Let $f\colon \mathcal{X}\to \mathcal{X}^\prime$ and $g\colon \mathcal{Y}\to \mathcal{Y}^\prime$ be measurable maps 
to finite sets $\mathcal{X}^\prime$ and $\mathcal{Y}^\prime$.
We define a probability mass function $\mu^\prime$ on $\mathcal{X}^\prime$ by
\[  \mu^\prime(x^\prime) = \mu\left(f^{-1}(x^\prime)\right)  \quad (x^\prime\in \mathcal{X}^\prime). \]
We also define conditional probability mass functions $\nu^\prime$ and $\nu^\prime_z$ $(z\in \mathcal{Z})$ 
on $\mathcal{X}^\prime\times \mathcal{Y}^\prime$ by
\[  \nu^\prime(y^\prime|x^\prime) = \frac{\int_{f^{-1}(x^\prime)} \nu\left(g^{-1}(y^\prime)|x\right) d\mu(x)}{\mu\left(f^{-1}(x^\prime)\right)}, 
   \quad   \nu^\prime_z(y^\prime|x^\prime) = \frac{\int_{f^{-1}(x^\prime)} \nu_z\left(g^{-1}(y^\prime)|x\right) d\mu(x)}{\mu\left(f^{-1}(x^\prime)\right)}\]
where $x^\prime \in \mathcal{X}^\prime$ and $y^\prime\in \mathcal{Y}^\prime$.
We have 
\[   \nu^\prime(y^\prime|x^\prime) =  \int_{\mathcal{Z}}  \nu^\prime_z(y^\prime|x^\prime) dm(z).  \]
Then by Lemma \ref{lemma: mutual information is convex in nu finite case}
\[   I(\mu^\prime, \nu^\prime) \leq \int_{\mathcal{Z}} I(\mu^\prime, \nu^\prime_z)\, dm(z). \]
It follows from the definition of mutual information that 
we have $I(\mu^\prime, \nu^\prime_z) \leq  I(\mu, \nu_z)$ for all $z\in \mathcal{Z}$.
Hence 
\[  I(\mu^\prime, \nu^\prime) \leq   \int_{\mathcal{Z}} I(\mu, \nu_z) \, dm(z). \]
(It follows from Theorem \ref{theorem: approximation of mutual information} that $I(\mu, \nu_z)$ is measurable in $z\in \mathcal{Z}$.)
Taking the supremum over $f$ and $g$, we get 
\[  I(\mu, \nu)  \leq  \int_{\mathcal{Z}} I(\mu, \nu_z) \, dm(z). \]
\end{proof}

Next we will establish a method to prove a lower bound on mutual information (Proposition \ref{proposition: Kawabata--Dembo} below).
We need to use the following integral representation of $I(X;Y)$. 
This is given in \cite[p. 176 Lemma 7.4, p. 206 Equation (7.31)]{Gray_entropy}.

\begin{theorem} \label{theorem: integral representation of mutual information}
Let $(\Omega, \mathcal{F}, \mathbb{P})$ be a probability space, and let $(\mathcal{X}, \mathcal{A})$ and $(\mathcal{Y}, \mathcal{B})$ be 
measurable spaces. Let $X\colon \Omega\to \mathcal{X}$ and $Y\colon \Omega\to \mathcal{Y}$ be measurable maps
with distributions $\mu = \mathrm{Law}(X) = X_*\mathbb{P}$ and $\nu = \mathrm{Law}(Y) = Y_*\mathbb{P}$ respectively.
Let $p = \mathrm{Law}(X, Y) = (X, Y)_*\mathbb{P}$ be the distribution of $(X, Y)\colon \Omega \to \mathcal{X}\times \mathcal{Y}$.
Suppose that the mutual information $I(X;Y)$ is finite.
Then $p$ is absolutely continuous with respect to the product measure $\mu\otimes \nu$.
Moreover, letting $f = dp/d(\mu\otimes \nu)$ be the Radon--Nikodim derivative, we have
\[  I(X;Y) = \int_{\mathcal{X}\times \mathcal{Y}} \log f\, dp = \int_{\mathcal{X}\times \mathcal{Y}}f\log f \, d(\mu\otimes \nu). \]
\end{theorem}

We learnt the next result from \cite[Lemma A.1]{Kawabata--Dembo}.
This is a kind of duality of convex programming.

\begin{proposition}  \label{proposition: duality of convex programming}
Let $\varepsilon>0$ and $a\geq 0$ be real numbers.
Let $(\mathcal{X}, \mathcal{A})$ and $(\mathcal{Y}, \mathcal{B})$ be measurable spaces and 
$\rho\colon \mathcal{X}\times \mathcal{Y}\to [0, +\infty)$ a measurable map.
Let $\mu$ be a probability measure on $\mathcal{X}$.
Suppose a measurable map $\lambda\colon \mathcal{X}\to  [0, +\infty)$ satisfies 
\[  \forall y\in \mathcal{Y}: \quad  \int_{\mathcal{X}} \lambda(x) 2^{-a \rho(x, y)}\, d\mu(x) \leq  1. \]
If $X$ and $Y$ are random variables taking values in $\mathcal{X}$ and $\mathcal{Y}$ respectively and 
satisfying $\mathrm{Law}(X) = \mu$ and $\mathbb{E}\rho(X, Y) < \varepsilon$ then we have
\begin{equation}  \label{eq: duality of convex programming}
     I(X;Y) \geq  -a\varepsilon + \int_{\mathcal{X}} \log \lambda(x) \, d\mu(x). 
\end{equation}    
\end{proposition}

\begin{proof}
Let $\nu = \mathrm{Law}(Y)$ and $p = \mathrm{Law}(X, Y)$ be the distributions of $Y$ and $(X, Y)$ respectively.
If $I(X;Y)$ is infinite then the statement is trivial. So we assume $I(X;Y)<\infty$.
Then by Theorem \ref{theorem: integral representation of mutual information}
the measure $p$ is absolutely continuous with respect to $\mu\otimes \nu$.
Let $f = dp/d(\mu\otimes \nu)$ be the Radon--Nikodim derivative.
We have 
\[  I(X;Y) = \int_{\mathcal{X}\times \mathcal{Y}} \log f\, dp. \]
Set $g(x, y) = \lambda(x) 2^{-a\rho(x, y)}$. 
Since $-\varepsilon < -\mathbb{E}\rho(X, Y) = -\int_{\mathcal{X}\times \mathcal{Y}} \rho(x, y) dp(x, y)$, we have 
\begin{align*}
  \int_{\mathcal{X}\times \mathcal{Y}} \log g(x, y) \, dp(x, y)  & 
    \geq  -a\varepsilon + \int_{\mathcal{X}\times \mathcal{Y}} \log \lambda(x) \, dp(x, y) \\
    & = -a\varepsilon + \int_{\mathcal{X}}\log \lambda(x) \, d\mu(x). 
\end{align*}
Therefore 
\[  I(X;Y) + a\varepsilon - \int_{\mathcal{X}}\log \lambda(x) \, d\mu(x) 
    \geq  \int_{\mathcal{X}\times \mathcal{Y}} (\log f - \log g) \, dp 
    = \int_{\mathcal{X}\times \mathcal{Y}}f \log (f/g) \, d\mu(x) d\nu(y).\]
Since $\ln t \leq t-1$, we have $\ln (1/t) \geq 1-t$ and hence 
$f\ln(f/g) \geq  f-g$. Then 
\begin{align*}
    (\ln 2) \int_{\mathcal{X}\times \mathcal{Y}} f \log (f/g) \, d\mu(x) d\nu(y)
     & = \int_{\mathcal{X}\times \mathcal{Y}}  f \ln (f/g) \, d\mu(x) d\nu(y)   \\
     & \geq   \int_{\mathcal{X}\times \mathcal{Y}} \left(f(x,y)-g(x,y)\right) d\mu(x) d\nu(y) \\
     & = 1 - \int_{\mathcal{Y}} \left(\int_{\mathcal{X}} g(x,y)\, d\mu(x)\right) d\nu(y) \geq  0.
\end{align*}
In the last inequality we have used the assumption $\int_{\mathcal{X}} g(x,y)\, d\mu(x) \leq 1$.
\end{proof}

The next proposition is a key result. We will use it for connecting geometric measure theory to 
rate distortion theory. This result is essentially due to Kawabata--Dembo \cite[Proposition 3.2]{Kawabata--Dembo}.
Recall that, for a metric space $(\mathcal{X}, \mathbf{d})$, we use the notation 
\[  \diam E = \sup\{\mathbf{d}(x, y)\mid  x, y\in E\} \quad (E\subset \mathcal{X}). \]

\begin{proposition}[Kawabata--Dembo estimate]  \label{proposition: Kawabata--Dembo}
Let $\varepsilon$ and $\delta$ be positive numbers with $2\varepsilon \log (1/\varepsilon) \leq \delta$.
Let $s$ be a nonnegative real number.
Let $(\mathcal{X}, \mathbf{d})$ be a separable metric space with a Borel probability measure $\mu$ satisfying 
\begin{equation}  \label{eq: power law}
    \mu(E) \leq  \left(\diam E\right)^s \quad \text{for all Borel sets $E\subset \mathcal{X}$ with $\diam E < \delta$}. 
\end{equation}    
Let $X$ and $Y$ be random variables taking values in $\mathcal{X}$ and satisfying 
$\mathrm{Law} X = \mu$ and $\mathbb{E}\mathbf{d}(X, Y) < \varepsilon$.
Then 
\[  I(X;Y) \geq  s\log (1/\varepsilon)  -K(s+1). \]
Here $K$ is a universal positive constant independent of the given data (i.e. $\varepsilon, \delta, s, (\mathcal{X}, \mathbf{d}), \mu$).
\end{proposition}

\begin{proof}
The proof is almost identical with \cite[Lemma 2.10]{Lindenstrauss--Tsukamoto double}.
But we repeat it for completeness.
If $s=0$ then the statement is trivial.
So we can assume $s>0$. We use Proposition \ref{proposition: duality of convex programming}.
Set $a= s/\varepsilon$ and we estimate $\int_{\mathcal{X}}2^{-a d(x, y)}d\mu(x)$ for each $y\in \mathcal{X}$.
By the Fubini theorem (see \cite[1.15 Theorem]{Mattila})
\[  \int_{\mathcal{X}}2^{-a d(x, y)}d\mu(x) = \int_0^1 \mu\{x\mid 2^{-a d(x,y)}\geq u\} \, du. \]
Changing the variable $u = 2^{-a v}$, we have $du = -a (\ln 2) 2^{-av} dv$ and hence 
\begin{align*}
    \int_0^1 \mu\{x\mid 2^{-a d(x,y)}\geq u\} \, du & = \int_0^\infty \mu\{x\mid d(x, y) \leq v\} a(\ln 2) 2^{-av}\, dv \\
    & = a\ln 2 \left(\int_0^{\delta/2} + \int_{\delta/2}^\infty\right)  \mu\{x\mid d(x, y) \leq v\} 2^{-av} \, dv.
\end{align*}    
By using (\ref{eq: power law})
\begin{align*}
    a\ln 2 \int_0^{\delta/2} \mu\{x\mid d(x, y) \leq v\} 2^{-av} \, dv  &\leq a\ln 2 \int_0^{\delta/2} (2v)^s 2^{-av}\, dv \\
   & = \int_0^{\frac{a\delta\ln 2}{2}}\left(\frac{2t}{a\ln 2}\right)^s e^{-t}\, dt,
   \quad  (t = a(\ln 2) v) \\
   & \leq \left(\frac{2}{a\ln 2}\right)^s \int_0^\infty t^s e^{-t}\, dt  \\
   & = \left(\frac{2\varepsilon}{\ln 2}\right)^s s^{-s} \Gamma(s+1), \quad    \left(a = \frac{s}{\varepsilon}\right). 
\end{align*}
On the other hand
\begin{align*}
  a\ln 2 \int_{\delta/2}^\infty  \mu\{x\mid d(x, y) \leq v\} 2^{-av} \, dv  & \leq  a \ln 2 \int_{\delta/2}^\infty 2^{-av} dv  \\
  & = 2^{-a\delta/2}  \\
  & = \left(2^{-\delta/(2\varepsilon)}\right)^s , \quad   \left(a = \frac{s}{\varepsilon}\right).
\end{align*}
Since $\delta\geq  2\varepsilon\log(1/\varepsilon)$, we have $-\frac{\delta}{2\varepsilon}\leq \log \varepsilon$. Hence
$2^{-\delta/(2\varepsilon)} \leq \varepsilon$.

Summing the above estimates, we get 
\[  \int_{\mathcal{X}}2^{-a d(x, y)}d\mu(x)  \leq \varepsilon^s\left\{1 + \left(\frac{2}{\ln 2}\right)^s s^{-s} \Gamma(s+1)\right\}. \]
Using the Stirling formula $s^{-s}\Gamma(s+1) \sim e^{-s}\sqrt{2\pi s}$, we can find a constant $c>1$ such that 
the term $\{\cdots\}$ is bounded by $c^{s+1}$ from above and hence
\[   \int_{\mathcal{X}}2^{-a d(x, y)}d\mu(x)  \leq  c^{s+1} \varepsilon^s. \]
We set $\lambda(x) = c^{-1-s} \varepsilon^{-s}$ for $x\in \mathcal{X}$. (This is a constant function.)
Then for all $y\in \mathcal{X}$
\[   \int_{\mathcal{X}}\lambda(x) 2^{-a d(x,y)}d\mu(x) \leq  1. \]
We apply Proposition \ref{proposition: duality of convex programming} and get 
\begin{align*}
  I(X;Y) & \geq -a\varepsilon + \int_{\mathcal{X}}\log \lambda \, d\mu \\
          & =  -s + \log \lambda ,  \quad  \left(a = \frac{s}{\varepsilon}\right) \\
          & = s\log(1/\varepsilon) - (1+\log c) s - \log c.
\end{align*}
Then the constant $K := 1+ \log c$ satisfies the statement.
\end{proof}

\subsection{Rate distortion theory}  \label{subsection: rate distortion theory}

In this subsection we introduce a rate distortion function.
A basic of rate distortion theory can be found in the book of Cover--Thomas \cite[Chapter 10]{Cover--Thomas}.
The rate distortion theory for continuous-time stochastic processes are 
presented in the paper of Pursley--Gray \cite{Pursley--Gray}. 

Recall that we have denoted the Lebesgue measure on $\mathbb{R}^d$ by $\mathbf{m}$.
For a measurable function $f(u)$ on $\mathbb{R}^d$ 
we usually denote its integral with respect to $\mathbf{m}$ by
\[  \int_{\mathbb{R}^d} f(u) du. \]

Let $(\mathcal{X}, \mathbf{d})$ be a compact metric space.
Let $A$ be a Borel subset of $\mathbb{R}^d$ of finite measure $\mathbf{m}(A) < \infty$.
We define $L^1(A, \mathcal{X})$ as the space of all measurable maps $f\colon A\to \mathcal{X}$.
We identify two maps if they coincide $\mathbf{m}$-almost everywhere.
We define a metric on $L^1(A, \mathcal{X})$ by 
\[  D(f, g) = \int_{A} \mathbf{d}\left(f(u), g(u)\right) du \quad  \left(f, g \in L^1(A, \mathcal{X})\right). \]

We need to check the following technical fact.
\begin{lemma}
 $\left(L^1(A, \mathcal{X}), D\right)$ is a complete separable metric space.
 Hence it is a standard Borel space with respect to the Borel $\sigma$-algebra.
\end{lemma}

\begin{proof}
First we need to understand what happens if we change the metric $\mathbf{d}$ on $\mathcal{X}$.
Let $\mathbf{d}^\prime$ be another metric on $\mathcal{X}$ compatible with the given topology.
We define a metric $D^\prime$ on $L^1(A, \mathcal{X})$ by 
\[  D^\prime(f, g) = \int_{A} \mathbf{d}^\prime \left(f(u), g(u)\right) du. \]
Let $\varepsilon$ be a positive number. There exists $\delta>0$ such that 
$\mathbf{d}(x, y) < \delta \Longrightarrow \mathbf{d}^\prime(x, y) < \varepsilon$.

Suppose $f, g\in L^1(A, \mathcal{X})$ satisfy $D(f, g) < \varepsilon \delta$. Then 
\[  \mathbf{m}\{u\in A\mid  \mathbf{d}\left(f(u), g(u)\right) \geq \delta\} \leq 
     \frac{1}{\delta} \int_{A} \mathbf{d}\left(f(u), g(u)\right) du  < \varepsilon. \]
We have $\mathbf{d}^\prime\left(f(u), g(u)\right) < \varepsilon$ on $\{u\in A\mid  \mathbf{d}\left(f(u), g(u)\right) < \delta\}$.
Hence 
\[  D^\prime(f, g) < \varepsilon \left(\diam(\mathcal{X}, \mathbf{d}^\prime) + \mathbf{m}(A)\right). \] 
So the identity map $\mathrm{id}\colon \left(L^1(A, \mathcal{X}), D\right) \to \left(L^1(A, \mathcal{X}), D^\prime\right)$ 
is uniformly continuous.
The same is true if we exchange $D$ and $D^\prime$.
Therefore if $\left(L^1(A, \mathcal{X}), D^\prime\right)$ is complete and separable then
so is $\left(L^1(A, \mathcal{X}), D\right)$.

Every compact metric space topologically embeds into the Hilbert cube $[0,1]^{\mathbb{N}}$.
We define a metric $\mathbf{d}^\prime$ on $[0,1]^{\mathbb{N}}$ by 
\[  \mathbf{d}^\prime(x, y) = \sum_{n=1}^\infty 2^{-n} |x_n-y_n|. \]
Let $L^1\left(A, [0,1]^{\mathbb{N}}\right)$ be the space of measurable maps from $A$ to $[0,1]^{\mathbb{N}}$.
We define a metric $D^\prime$ on $L^1\left(A, [0,1]^{\mathbb{N}}\right)$ as above.
The space $L^1(A, \mathcal{X})$ is identified with a closed subspace of $L^1\left(A, [0,1]^{\mathbb{N}}\right)$.
So it is enough to show that $\left(L^1\left(A, [0,1]^{\mathbb{N}}\right), D^\prime\right)$ is a complete separable metric space.
This follows from the standard fact that $L^1(A, [0,1])$ is complete and separable with respect to the $L^1$-norm.
\end{proof}
In the following we always assume that $L^1(A, \mathcal{X})$ is endowed with the Borel $\sigma$-algebra 
(and hence it is a standard Borel space).

Let $(\mathcal{X}, \mathbf{d})$ be a compact metric space, and let $T\colon \mathbb{R}^d\times \mathcal{X}\to \mathcal{X}$
be a continuous action of $\mathbb{R}^d$.
Let $\mu$ be a $T$-invariant Borel probability measure on $\mathcal{X}$.

Let $\varepsilon>0$ and let $A$ be a bounded Borel subset of $\mathbb{R}^d$ with $\mathbf{m}(A)>0$.
We define $R(\varepsilon, A)$ as the infimum of the mutual information $I(X; Y)$ where $X$ and $Y$ are random variables defined on 
some probability space $(\Omega, \mathcal{F}, \mathbb{P})$ such that 
\begin{itemize}
   \item $X$ takes values in $\mathcal{X}$ and its distribution is given by $\mu$,
   \item $Y$ takes values in $L^1(A, \mathcal{X})$ and satisfies
   \[   \mathbb{E}\left(\frac{1}{\mathbf{m}(A)}\int_A\mathbf{d}(T^u X, Y_u)\, du\right) < \varepsilon. \]
\end{itemize}
Here $Y_u = Y_u(\omega)$ ($\omega\in \Omega$) is the value of $Y(\omega)\in L^1(A, \mathcal{X})$ at $u\in A$.
We set $R(\varepsilon,A) = 0$ if $\mathbf{m}(A) = 0$.

\begin{remark} \label{remark: we can assume that Y takes only finitely many values}
In the above definition of $R(\varepsilon, A)$, we can assume that $Y$ takes only finitely many values.
Indeed let $X$ and $Y$ be random variables satisfying the conditions in the definition of $R(\varepsilon, A)$.
We take a positive number $\tau$ satisfying 
\[   \mathbb{E}\left(\frac{1}{\mathbf{m}(A)}\int_A\mathbf{d}(T^u X, Y_u)\, du\right) < \varepsilon - 2\tau. \]
Since $L^1(A, \mathcal{X})$ is separable, it contains a dense countable subsets $\{f_1, f_2, f_3, \dots\}$.
We define a map $F\colon L^1(A, \mathcal{X})\to \{f_1, f_2, f_3, \dots\}$ by $F(f)  = f_n$ where $n$ is the smallest 
natural number satisfying $D(f, f_n) < \tau\cdot \mathbf{m}(A)$. 
Set $Y^\prime = F(Y)$. Then we have 
\[   \mathbb{E}\left(\frac{1}{\mathbf{m}(A)}\int_A\mathbf{d}(T^u X, Y^\prime_u)\, du\right) < \varepsilon - \tau. \]
Define $p_n = \mathbb{P}(Y^\prime = f_n)$. We choose $n_0$ such that 
\[  \sum_{n>n_0} p_n \diam(\mathcal{X}, \mathbf{d}) < \tau. \]
We define $G\colon \{f_1, f_2, f_3, \dots\}\to \{f_1, f_2, \dots, f_{n_0}\}$ by 
\begin{equation*}
   G(f) = \begin{cases}  f & \text{if } f\in \{f_1, f_2, \dots, f_{n_0}\} \\
                                 f_{n_0} & \text{otherwise}   \end{cases}. 
\end{equation*}
Set $Y^{\prime\prime} = G(Y^\prime)$. Then $Y^{\prime\prime}$ takes only finitely many values (i.e. $f_1, \dots, f_{n_0}$)
and we have
\[   \mathbb{E}\left(\frac{1}{\mathbf{m}(A)}\int_A\mathbf{d}(T^u X, Y^{\prime\prime}_u)\, du\right) < \varepsilon, \]
\[  I(X; Y^{\prime\prime}) \leq I(X; Y^\prime) \leq  I(X; Y). \]
Therefore, when we consider the infimum in the definition of $R(\varepsilon, A)$, we only need to take into account 
such random variables $Y^{\prime\prime}$.
\end{remark}

For a bounded Borel subset $A\subset \mathbb{R}^d$ and $r>0$ 
we define $N_r(A)$ as the $r$-neighborhood of $A$ with respect to the $\ell^\infty$-norm, i.e.
$N_r(A) = \{u +v\mid u\in A, v\in (-r, r)^d\}$.
\begin{lemma}   \label{lemma: subadditivity of R(A, varepsilon)}
We have:
  \begin{enumerate}
   \item  $R(\varepsilon, A) \leq \log \#\left(\mathcal{X}, \mathbf{d}_A,\varepsilon\right)
    \leq \mathbf{m}\left(N_{1/2}(A)\right) \log \#(\mathcal{X}, \mathbf{d}_{(-1, 1)^d},\varepsilon)$.
   \item  $R(\varepsilon, a+A) = R(\varepsilon, A)$ for any $a\in \mathbb{R}^d$.
   \item  If $A\cap B = \emptyset$ then $R(\varepsilon, A\cup B) \leq  R(\varepsilon, A) + R(\varepsilon,B)$.
  \end{enumerate}
\end{lemma}

\begin{proof}
(1) Let $\mathcal{X}= U_1\cup U_2\cup \dots \cup U_n$ be an open cover with $n = \#\left(\mathcal{X}, \mathbf{d}_A,\varepsilon\right)$
and $\diam (U_k, \mathbf{d}_A) < \varepsilon$ for all $1\leq k \leq n$. 
Take a point $x_k\in U_k$ for each $k$ and define a map $f\colon \mathcal{X}\to \{x_1,\dots, x_n\}$
by $f(x)  = \{x_k\}$ for $x\in U_k\setminus (U_1\cup \dots \cup U_{k-1})$.
Let $X$ be a random variable taking values in $\mathcal{X}$ according to $\mu$.
We set $Y_u = T^u f(X)$ for $u\in A$. Then $X$ and $Y$ satisfy the conditions of the definition of $R(\varepsilon, A)$.
We have 
\[  R(\varepsilon, A)  \leq I(X;Y) \leq H(Y)  \leq \log n =  \log \#\left(\mathcal{X}, \mathbf{d}_A,\varepsilon\right). \]
We estimate $\log \#\left(\mathcal{X}, \mathbf{d}_A,\varepsilon\right)$.
Let $\{u_1, \dots, u_a\}$ be a maximal $1$-separated subset of $A$ where “$1$-separated” means $\norm{u_i-u_j}_\infty \geq 1$ for $i\neq  j$.
Then $A\subset \bigcup_{i=1}^a \left(u_i+(-1, 1)^d\right)$ and hence 
\[  \log \#\left(\mathcal{X}, \mathbf{d}_A,\varepsilon\right) \leq  a \log \#\left(\mathcal{X}, \mathbf{d}_{(-1, 1)^d},\varepsilon\right). \]
The sets $u_i+(-1/2,1/2)^d$ $(1\leq i\leq a)$ are mutually disjoint and contained in $N_{1/2}(A)$. Therefore
$a\leq \mathbf{m}\left(N_{1/2}(A)\right)$.

(2) Let $X$ and $Y$ be random variables satisfying the conditions of the definition of $R(\varepsilon, A)$
(i.e., $X$ is distributed according to $\mu$ and the average distance between $\{T^u X\}_{u\in A}$ and $Y$ is bounded by $\varepsilon$).
We define new random variables $X^\prime$ and $Y^\prime$ by 
\[  X^\prime = T^{-a} X, \quad  Y^\prime_v = Y_{v-a} \quad  (v\in a + A). \]
Since $\mu$ is $T$-invariant, we have $\mathrm{Law} X^\prime  = \mathrm{Law} X = \mu$. 
The random variable $Y^\prime$ takes values in $L^1(a+A, \mathcal{X})$ and 
\begin{align*}
  \int_{a+A}\mathbf{d}(T^v X^\prime, Y^\prime_v) \, dv  & = \int_{a+A}\mathbf{d}(T^{v-a}X, Y_{v-a}) \, dv \\
  & = \int_A  \mathbf{d}(T^u X, Y_u)\, du, \quad  (u=v-a).
\end{align*}
We have $I(X^\prime; Y^\prime) = I(X;Y)$. Therefore $R(\varepsilon, a+A) = R(\varepsilon, A)$.

(3) Let $X$ and $Y$ be random variables satisfying the conditions of the definition of $R(\varepsilon, A)$ as above, and let 
$X^\prime$ and $Y^\prime$ be random variables satisfying the conditions of the definition of $R(\varepsilon, B)$.
We denote by $\mathbb{P}(Y\in E\mid X=x)$ 
$(E\subset L^1(A, \mathcal{X}))$ the regular conditional distribution of $Y$ given $X=x$.
Similarly for $\mathbb{P}(Y^\prime\in F\mid X^\prime=x)$.

We naturally identify $L^1(A\cup B, \mathcal{X})$ with $L^1(A, \mathcal{X})\times L^1(B, \mathcal{X})$.
We define a transition probability $\nu$ on $\mathcal{X}\times L^1(A\cup B, \mathcal{X})$ by 
\[  \nu(E\times F|x) = \mathbb{P}(Y\in E\mid X=x) \mathbb{P}(Y^\prime \in F\mid X^\prime=x),  \]
for $E\times F\subset L^1(A, \mathcal{X})\times L^1(B, \mathcal{X}) = L^1(A\cup B, \mathcal{X})$ and $x\in \mathcal{X}$. 
We define a probability measure $Q$ on $\mathcal{X}\times L^1(A\cup B, \mathcal{X})$ by
\[  Q(G) = \int_{\mathcal{X}} \nu(G_x| x)\, d\mu(x), \quad 
     (G\subset \mathcal{X}\times  L^1(A\cup B, \mathcal{X})),  \]
where $G_x = \{f\in L^1(A\cup B, \mathcal{X})\mid  (x, f)\in  G\}$.
Let $(X^{\prime\prime}, Y^{\prime\prime})$ be the random variable taking values in $\mathcal{X}\times L^1(A\cup B, \mathcal{X})$
according to $Q$.
Then $\mathrm{Law} X^{\prime\prime} = \mu$ and 
\begin{align*}
   \mathbb{E}\left(\int_{A\cup B} \mathbf{d}(T^u X^{\prime\prime}, Y^{\prime\prime}_u)\, du\right)
   & = \mathbb{E}\left(\int_A \mathbf{d}(T^u X, Y_u)\, du\right) + \mathbb{E}\left(\int_B \mathbf{d}(T^u X^\prime, Y^\prime_u)\, du\right) \\
   & < \varepsilon \, \mathbf{m}(A) + \varepsilon \, \mathbf{m}(B) = \varepsilon\, \mathbf{m}(A\cup B). 
\end{align*}   
The random variables $Y^{\prime\prime}|_A$ and $Y^{\prime\prime}|_B$ is conditionally independent given $X^{\prime\prime}$.
Therefore by Lemma \ref{lemma: subadditivity of mutual information}
\[  I(X^{\prime\prime}; Y^{\prime\prime}) = I(X^{\prime\prime}; Y^{\prime\prime}|_A, Y^{\prime\prime}|_B)
     \leq  I(X^{\prime\prime}; Y^{\prime\prime}|_A) + I(X^{\prime\prime}; Y^{\prime\prime}|_B)
     = I(X; Y) + I(X^\prime, Y^\prime). \]
The statement (3) follows from this.
\end{proof}

\begin{lemma}   \label{lemma: existence of the limit in defining rate distortion function}
The limit of $\displaystyle \frac{R\left(\varepsilon, [0, L)^d\right)}{L^d}$ as $L\to \infty$ exists and is equal to the infimum of 
$\displaystyle \frac{R\left(\varepsilon, [0, L)^d\right)}{L^d}$ over $L>0$.
\end{lemma}

\begin{proof}
Let $0<\ell <L$. We divide $L$ by $\ell$ and let $L = q\ell + r$ where $q$ is a natural number and $0\leq r < \ell$.
Set 
\[  \Gamma = \{ (\ell n_1, \dots, \ell n_d)\mid  n_i \in \mathbb{Z}, \, 0 \leq n_i < q \> (1\leq i \leq d)\}. \]
The cubes $u +[0, \ell)^d$ $(u\in \Gamma)$ are disjoint and contained in $[0, L)^d$.
Let $A$ be the complement:
\[  A = [0, L)^d\setminus \bigcup_{u\in \Gamma} \left(u+ [0, \ell)^d\right). \]
The volume of the $1/2$-neighborhood of $A$ is $O(L^{d-1})$:
\[  \mathbf{m}\left(N_{1/2}(A)\right) \leq  d(r+1)(L+1)^{d-1}. \]
By Lemma \ref{lemma: subadditivity of R(A, varepsilon)}
\begin{align*}
   R\left(\varepsilon, [0,L)^d\right) & \leq \sum_{u\in \Gamma} R\left(\varepsilon, u+[0,\ell)^d\right) + R(\varepsilon, A) \\
   & \leq  q^d R\left(\varepsilon, [0,\ell)^d\right) + C (L+1)^{d-1}.
\end{align*}
By dividing this by $L^d$ and letting $L\to \infty$, we get 
\[    \limsup_{L\to \infty} \frac{R\left(\varepsilon, [0, L)^d\right)}{L^d} \leq  \frac{R\left(\varepsilon, [0,\ell)^d\right)}{\ell^d}. \]
Then
\[  \limsup_{L\to \infty} \frac{R\left(\varepsilon, [0, L)^d\right)}{L^d} \leq  
     \inf_{\ell>0} \frac{R\left(\varepsilon, [0,\ell)^d\right)}{\ell^d}  \leq
     \liminf_{\ell\to \infty} \frac{R\left(\varepsilon, [0,\ell)^d\right)}{\ell^d} \]
\end{proof}

Recall that $T\colon \mathbb{R}^d\times \mathcal{X}\to \mathcal{X}$ is a continuous action of $\mathbb{R}^d$ on a 
compact metric space $(\mathcal{X}, \mathbf{d})$ with an invariant probability measure $\mu$.
For $\varepsilon>0$ we define the \textbf{rate distortion function} $R(\mathbf{d}, \mu, \varepsilon)$ $(\varepsilon>0)$ by 
\[  R(\mathbf{d}, \mu, \varepsilon) = \lim_{L\to \infty} \frac{R\left(\varepsilon, [0,L)^d\right)}{L^d}
    = \inf_{L>0}  \frac{R\left(\varepsilon, [0,L)^d\right)}{L^d}. \]
We define the \textbf{upper/lower rate distortion dimensions} of $(\mathcal{X}, T, \mathbf{d}, \mu)$ by 
\[  \urdim(\mathcal{X}, T, \mathbf{d}, \mu) = \limsup_{\varepsilon\to 0} \frac{R(\mathbf{d}, \mu, \varepsilon)}{\log(1/\varepsilon)}, \quad
     \lrdim(\mathcal{X}, T, \mathbf{d}, \mu) =  \liminf_{\varepsilon\to 0} \frac{R(\mathbf{d}, \mu, \varepsilon)}{\log(1/\varepsilon)}. \]

\begin{remark}
A tiling argument similar to the proof of Lemma \ref{lemma: existence of the limit in defining rate distortion function} shows that 
if $\Lambda_1, \Lambda_2, \Lambda_3, \dots$ are a sequence of rectangles of $\mathbb{R}^d$ such that the minimum side length of 
$\Lambda_n$ diverges to infinity then we have 
\[   R(\mathbf{d}, \mu, \varepsilon) = \lim_{n\to \infty} \frac{R(\varepsilon, \Lambda_n)}{\mathbf{m}(\Lambda_n)}. \]
With a bit more effort we can also prove that 
\[  R(\mathbf{d}, \mu, \varepsilon) = \lim_{r\to \infty} \frac{R(\varepsilon, B_r)}{\mathbf{m}(B_r)} \]
where $B_r$ is the Euclidean $r$-ball of $\mathbb{R}^d$ centered at the origin.
But we are not sure whether, for any F{\o}lner sequence $A_1, A_2, A_3,\dots$ of $\mathbb{R}^d$, the limit of 
$\frac{R(\varepsilon, A_n)}{\mathbf{m}(A_n)}$ exists or not. (Maybe not.)
Probably we need to modify the definition of rate distortion function when we study the rate distortion theory 
for actions of general amenable groups. 
\end{remark}

\section{Metric mean dimension with potential and mean Hausdorff dimension with potential}
\label{section: metric mean dimension with potential and mean Hausdorff dimension with potential}

The purpose of this section is to introduce \textit{metric mean dimension with potential} and
\textit{mean Hausdorff dimension with potential} for $\mathbb{R}^d$-actions.
These are dynamical versions of Minkowski dimension and Hausdorff dimension.
Mean Hausdorff dimension with potential is a main ingredient of the proof of Theorem \ref{main theorem}.
Metric mean dimension with potential is not used in the proof of Theorem \ref{main theorem}.
But it is also an indispensable tool of mean dimension theory.
Therefore we develop its basic theory. 
We plan to use it in Part II of this series of papers.

Let $(\mathcal{X}, \mathbf{d})$ be a compact metric space and $\varphi\colon \mathcal{X}\to \mathbb{R}$ a continuous function.
Let $\varepsilon$ be a positive number. 
We define \textbf{covering number with potential} by
\[ \#\left(\mathcal{X}, \mathbf{d}, \varphi, \varepsilon\right) =  
    \inf\left\{\sum_{i=1}^n (1/\varepsilon)^{\sup_{U_i}\varphi} \middle|\, 
    \parbox{3in}{\centering  $\mathcal{X} = U_1\cup \dots \cup U_n$ is an open cover with $\diam\, U_i < \varepsilon$
     for all $1\leq i\leq n$} \right\}. \]
Let $s$ be a real number larger than the maximum value of $\varphi$.
We define 
\[   \mathcal{H}^s_\varepsilon(\mathcal{X}, \mathbf{d},\varphi) = 
      \inf\left\{\sum_{n=1}^\infty \left(\diam E_n\right)^{s-\sup_{E_n} \varphi} \middle|\, 
       \mathcal{X} = E_1\cup E_2\cup E_3\cup \dots \text{ with } \diam E_n < \varepsilon\right\}. \]
Here we assume that the empty set has diameter zero.
We define \textbf{Hausdorff dimension with potential at the scale $\varepsilon$} by 
\[  \dimh(\mathcal{X},\mathbf{d},\varphi,\varepsilon) 
    = \inf\{s  \mid  \mathcal{H}^s_{\varepsilon}(\mathcal{X},\mathbf{d},\varphi) < 1, \, s >\max\varphi\}. \]
    
When the function $\varphi$ is identically zero ($\varphi \equiv  0$), we denote $\#\left(\mathcal{X}, \mathbf{d}, \varphi, \varepsilon\right)$,  
 $\mathcal{H}^s_\varepsilon(\mathcal{X}, \mathbf{d},\varphi)$ and 
$\dimh(\mathcal{X},\mathbf{d},\varphi,\varepsilon)$ by $\#\left(\mathcal{X}, \mathbf{d},\varepsilon\right)$,
 $\mathcal{H}^s_\varepsilon(\mathcal{X}, \mathbf{d})$ and 
$\dimh(\mathcal{X},\mathbf{d},\varepsilon)$ respectively:
\begin{align*}
     \#\left(\mathcal{X}, \mathbf{d},\varepsilon\right) 
      & =  \min \left\{n \middle|\, 
    \parbox{3in}{\centering  $\mathcal{X} = U_1\cup \dots \cup U_n$ is an open cover with $\diam\, U_i < \varepsilon$
                   for all $1\leq i\leq n$} \right\},  \\
     \mathcal{H}^s_\varepsilon(\mathcal{X}, \mathbf{d})  & = 
     \inf\left\{\sum_{n=1}^\infty \left(\diam E_n\right)^{s} \middle|\, 
       \mathcal{X} = E_1\cup E_2\cup E_3\cup \dots \text{ with } \diam E_n < \varepsilon\right\},  \\
     \dimh(\mathcal{X},\mathbf{d},\varepsilon)  
     & =  \inf\{s  \mid  \mathcal{H}^s_{\varepsilon}(\mathcal{X},\mathbf{d}) < 1, \, s > 0 \}. 
\end{align*}     
We assume that $\dimh(\mathcal{X},\mathbf{d},\varepsilon) = -\infty$ if $\mathcal{X}$ is the empty set.
    
\begin{lemma} \label{lemma: Hausdorff dimension and covering number}
  For $0<\varepsilon < 1$, we have 
  $\displaystyle \dimh(\mathcal{X},\mathbf{d},\varphi,\varepsilon) 
   \leq \frac{\log \#(\mathcal{X},\mathbf{d},\varphi,\varepsilon)}{\log(1/\varepsilon)}$.
\end{lemma}
\begin{proof}
Suppose $\displaystyle \frac{\log \#(\mathcal{X},\mathbf{d},\varphi,\varepsilon)}{\log(1/\varepsilon)} <s$.
We have $s>\max \varphi$.  
There is an open cover $\mathcal{X} = U_1\cup\dots\cup U_n$ with $\diam U_i < \varepsilon$ and 
\[   \sum_{i=1}^n (1/\varepsilon)^{\sup_{U_i} \varphi} < (1/\varepsilon)^s. \]
Then 
\[   \sum_{i=1}^n (\diam U_i)^{s- \sup_{U_i} \varphi}  \leq  \sum_{i=1}^n \varepsilon^{s- \sup_{U_i} \varphi} < 1. \]
Hence $\mathcal{H}^s_{\varepsilon}(\mathcal{X},\mathbf{d},\varepsilon) < 1$ and we have 
$\dimh(\mathcal{X},\mathbf{d},\varphi,\varepsilon) \leq  s$.
\end{proof}

Let $T\colon \mathbb{R}^d\times \mathcal{X}\to \mathcal{X}$ be a continuous action of the group $\mathbb{R}^d$.
For a bounded Borel subset $A$ of $\mathbb{R}^d$, as in \S \ref{subsection: mean dimension with potential of R^d-actions}, we define 
a metric $\mathbf{d}_A$ and a function $\varphi_A$ on $\mathcal{X}$ by 
\[  \mathbf{d}_A(x, y) = \sup_{u\in A}\mathbf{d}(T^u x, T^u y), \quad 
     \varphi_A(x)  = \int_{A}\varphi(T^u x)\, du. \]
In particular, for a positive number $L$, we set $\mathbf{d}_L = \mathbf{d}_{[0,L)^d}$ and $\varphi_L = \varphi_{[0,L)^d}$.
We define \textbf{upper/lower metric mean dimension with potential} by 
\begin{equation} \label{eq: metric mean dimension with potential}
   \begin{split}
     \umdimm(\mathcal{X}, T, \mathbf{d}, \varphi) &= \limsup_{\varepsilon\to 0} 
     \left(\lim_{L\to \infty} \frac{\log \#\left(\mathcal{X}, \mathbf{d}_L, \varphi_L,\varepsilon\right)}{L^d \log(1/\varepsilon)}\right),  \\
     \lmdimm(\mathcal{X}, T, \mathbf{d}, \varphi) &= \liminf_{\varepsilon\to 0} 
     \left(\lim_{L\to \infty} \frac{\log \#\left(\mathcal{X}, \mathbf{d}_L, \varphi_L,\varepsilon\right)}{L^d \log(1/\varepsilon)}\right). 
   \end{split}  
\end{equation}     
Here the limit with respect to $L$ exists. The proof is similar to the proof of Lemma \ref{lemma: existence of limit with respect to L}.
(The quantity $\log \#\left(\mathcal{X}, \mathbf{d}_A, \varphi_A, \varepsilon\right)$ satisfies the natural subadditivity, monotonicity 
and invariance if $\varphi(x)$ is a nonnegative function.)

We define \textbf{upper/lower mean Hausdorff dimension with potential} by 
\begin{align*}
   \umdimh(\mathcal{X}, T, \mathbf{d}, \varphi) & = \lim_{\varepsilon\to 0} 
   \left(\limsup_{L\to \infty} \frac{\dimh\left(\mathcal{X},\mathbf{d}_L, \varphi_L,\varepsilon\right)}{L^d}\right), \\
   \lmdimh(\mathcal{X}, T, \mathbf{d}, \varphi) & = \lim_{\varepsilon\to 0} 
   \left(\liminf_{L\to \infty} \frac{\dimh\left(\mathcal{X},\mathbf{d}_L, \varphi_L,\varepsilon\right)}{L^d}\right).
\end{align*}

\begin{remark}
We are not sure whether or not these definitions of the upper and lower mean Hausdorff dimensions with potential 
coincide with the following:
\[
 \lim_{\varepsilon\to 0} 
   \left(\limsup_{r\to \infty} \frac{\dimh\left(\mathcal{X},\mathbf{d}_{B_r}, \varphi_{B_r},\varepsilon\right)}{\mathbf{m}(B_r)}\right), 
  \quad  \lim_{\varepsilon\to 0} 
   \left(\liminf_{r\to \infty} \frac{\dimh\left(\mathcal{X},\mathbf{d}_{B_r}, \varphi_{B_r},\varepsilon\right)}{\mathbf{m}(B_r)}\right).
\]
(Maybe not in general.) Here $B_r$ is the Euclidean $r$-ball of $\mathbb{R}^d$ centered at the origin.
\end{remark}

The next lemma is a dynamical version of the fact that Hausdorff dimension is smaller than or equal to Minkowski dimension.

\begin{lemma} \label{lemma: mean Hausdorff dimension is bounded by metric mean dimension}
  $\umdimh(\mathcal{X},T,\mathbf{d},\varphi) \leq \lmdimm(\mathcal{X},T,\mathbf{d},\varphi)$.
\end{lemma}
\begin{proof}
This follows from Lemma \ref{lemma: Hausdorff dimension and covering number}.
\end{proof}

For the proof of Theorem \ref{main theorem} we need another version of 
mean Hausdorff dimension.
For a positive number $L$
we define a metric $\overline{\mathbf{d}}_L$ on $\mathcal{X}$ by
\[  \overline{\mathbf{d}}_L(x,y) = \frac{1}{L^d}\int_{[0,L)^d} d(T^u x, T^u y)\, du. \]
This is also compatible with the given topology of $\mathcal{X}$.
We define \textbf{upper/lower $L^1$-mean Hausdorff dimension with potential} by
\begin{equation}  \label{eq: L^1 mean Hausdorff dimension}
  \begin{split}
   \overline{\mdim}_{\mathrm{H},L^1}(\mathcal{X}, T, \mathbf{d}, \varphi) 
& = \lim_{\varepsilon\to 0} 
   \left(\limsup_{L\to \infty} \frac{\dimh\left(\mathcal{X},\overline{\mathbf{d}}_L, \varphi_L,\varepsilon\right)}{L^d}\right), \\
   \underline{\mdim}_{\mathrm{H},L^1}(\mathcal{X}, T, \mathbf{d}, \varphi) 
& = \lim_{\varepsilon\to 0} 
   \left(\liminf_{L\to \infty} \frac{\dimh\left(\mathcal{X},\overline{\mathbf{d}}_L, \varphi_L,\varepsilon\right)}{L^d}\right).
   \end{split}
\end{equation}
We have $\overline{\mathbf{d}}_L(x,y) \leq \mathbf{d}_L(x,y)$ and hence 
\[ \overline{\mdim}_{\mathrm{H},L^1}(\mathcal{X}, T, \mathbf{d}, \varphi)  \leq 
    \umdimh(\mathcal{X}, T, \mathbf{d}, \varphi), \quad 
      \underline{\mdim}_{\mathrm{H},L^1}(\mathcal{X}, T, \mathbf{d}, \varphi)  \leq 
      \lmdimh(\mathcal{X}, T, \mathbf{d}, \varphi). \]

It is well-known that topological dimension is smaller than or equal to Hausdorff dimension.
The next result is its dynamical version. The proof will be given in 
\S \ref{section: proof of mdim is bounded by mdimh}.

\begin{theorem}\label{theorem: mdim is bounded by mdimh}
 $\mdim(\mathcal{X}, T, \varphi) \leq  
\underline{\mdim}_{\mathrm{H},L^1}(\mathcal{X}, T, \mathbf{d}, \varphi)$.
\end{theorem}

Notice that this also implies $\mdim(\mathcal{X}, T, \varphi) \leq  
\underline{\mdim}_{\mathrm{H}}(\mathcal{X}, T, \mathbf{d}, \varphi)$.

Metric mean dimension with potential is related to rate distortion dimension by the following result.

\begin{proposition} \label{proposition: rate distortion dimension and metric mean dimension}
For any $T$-invariant Borel probability measure $\mu$ on $\mathcal{X}$ we have 
  \begin{align*}
    \urdim(\mathcal{X},T,\mathbf{d},\mu) +\int_{\mathcal{X}}\varphi\, d\mu & \leq \umdimm(\mathcal{X},T,\mathbf{d},\varphi), \\
    \lrdim(\mathcal{X},T,\mathbf{d},\mu) + \int_{\mathcal{X}}\varphi\, d\mu & \leq \lmdimm(\mathcal{X},T,\mathbf{d},\varphi).
  \end{align*}
\end{proposition}

\begin{proof}
We will use the following well-known inequality. For the proof see \cite[\S 9.3, Lemma 9.9]{Walters}.
  \begin{lemma}  \label{lemma: free energy}
  Let $a_1, \dots, a_n$ be real numbers and $(p_1, \dots, p_n)$ a probability vector.
  For a positive number $\varepsilon$ we have 
  \[  \sum_{i=1}^n \left(-p_i\log p_i + p_i a_i \log(1/\varepsilon)\right) \leq \log \sum_{i=1}^n (1/\varepsilon)^{a_i}. \]
  \end{lemma}
Let $L$ and $\varepsilon$ be positive numbers.
Let $\mathcal{X}= U_1\cup\dots \cup U_n$ be an open cover with 
$\diam(U_i, \mathbf{d}_L)<\varepsilon$.
Take a point $x_i\in U_i$ for each $1\leq i\leq n$.
Set $E_i = U_i\setminus (U_1\cup\dots \cup U_{i-1})$ and 
$p_i = \mu(E_i)$.
We define $f\colon \mathcal{X}\to \{x_1, \dots, x_n\}$ by $f(E_i) =\{x_i\}$.
Let $X$ be a random variable taking values in $\mathcal{X}$ according to $\mu$.
Set $Y_u = T^u f(X)$ for $u\in [0,L)^d$.
Then almost surely we have $\mathbf{d}(T^u X, Y_u) < \varepsilon$ for all $u\in [0,L)^d$.
Therefore 
\[  R\left(\varepsilon, [0,L)^d\right) \leq H(Y) = H\left(f(X)\right)
      = -\sum_{i=1}^n p_i \log p_i. \]
On the other hand 
\[  \int_{\mathcal{X}} \varphi\, d\mu = \int_{\mathcal{X}} \frac{\varphi_L}{L^d}\, d\mu 
     \leq \frac{1}{L^d} \sum_{i=1}^n p_i \sup_{U_i} \varphi_L. \]
By Lemma \ref{lemma: free energy}
\begin{align*}
   R\left(\varepsilon, [0,L)^d\right) 
    + L^d \log(1/\varepsilon) \int_{\mathcal{X}}\varphi \, d\mu  & \leq 
   \sum_{i=1}^n \left(-p_i\log p_i + p_i \sup_{U_i} \varphi_L \log (1/\varepsilon)\right)  \\
   & \leq  \log \sum_{i=1}^n (1/\varepsilon)^{\sup_{U_i} \varphi_L}. 
\end{align*}
Therefore 
\[  \frac{R\left(\varepsilon, [0,L)^d\right)}{L^d \log(1/\varepsilon)} 
     + \int_{\mathcal{X}}\varphi \, d\mu \leq 
\frac{\log\#\left(\mathcal{X}, \mathbf{d}_L,\varphi_L,\varepsilon\right)}{L^d\log(1/\varepsilon)}.
  \]
Letting $L\to \infty$
\[ \frac{R(\mathbf{d},\mu,\varepsilon)}{\log(1/\varepsilon)} 
     + \int_{\mathcal{X}}\varphi \, d\mu \leq  \lim_{L\to \infty}
\frac{\log\#\left(\mathcal{X}, \mathbf{d}_L,\varphi_L,\varepsilon\right)}{L^d\log(1/\varepsilon)}.
  \]
Letting $\varepsilon\to 0$ we get the statement.
\end{proof}

$L^1$-mean Hausdorff dimension with potential is related to rate distortion dimension by 
the next theorem.
We call this result ``dynamical Frostman's lemma" because the classical Frostman's lemma 
\cite[Sections 8.14-8.17]{Mattila} plays an essential role in its proof.
Recall that we have denoted by $\mathscr{M}^T(\mathcal{X})$
the set of $T$-invariant Borel probability measures on $\mathcal{X}$.

\begin{theorem}[Dynamical Frostman's lemma] \label{theorem: dynamical Frostman lemma}
\[  \overline{\mdim}_{\mathrm{H},L^1}(\mathcal{X}, T, \mathbf{d}, \varphi) \leq 
     \sup_{\mu\in \mathscr{M}^T(\mathcal{X})}\left(\lrdim(\mathcal{X},T,\mathbf{d},\mu)
     + \int_{\mathcal{X}}\varphi\, d\mu\right). \]
\end{theorem}

The proof of this theorem is the most important step of the proof of Theorem 
\ref{main theorem}. It will be given in \S \ref{section: proof of dynamical Frostman lemma}.

By combining Theorems \ref{theorem: mdim is bounded by mdimh} and 
\ref{theorem: dynamical Frostman lemma}, we get 
\[ \mdim(\mathcal{X},T,\varphi) \leq 
    \sup_{\mu\in \mathscr{M}^T(\mathcal{X})}\left(\lrdim(\mathcal{X},T,\mathbf{d},\mu)
     + \int_{\mathcal{X}}\varphi\, d\mu\right). \]
This is the statement of Theorem \ref{main theorem}.
Therefore the proof of Theorem \ref{main theorem} is reduced to the proofs of 
Theorems \ref{theorem: mdim is bounded by mdimh} and 
\ref{theorem: dynamical Frostman lemma}.

\begin{conjecture}  \label{conjecture: dynamical PS}
Let $T\colon \mathbb{R}^d\times \mathcal{X}\to \mathcal{X}$ be a continuous action of 
$\mathbb{R}^d$ on a compact metrizable space $\mathcal{X}$.
For any continuous function $\varphi\colon \mathcal{X}\to \mathbb{R}$
there exists a metric $\mathbf{d}$ on $\mathcal{X}$ compatible with the given topology satisfying
\begin{equation} \label{eq: dynamical PS}
   \mdim(\mathcal{X},T,\varphi) = \umdimm(\mathcal{X},T,\mathbf{d},\varphi). 
\end{equation}
\end{conjecture}

Suppose that this conjecture is true and let 
$\mathbf{d}$ be a metric satisfying (\ref{eq: dynamical PS}).
Then by Proposition 
\ref{proposition: rate distortion dimension and metric mean dimension}
\begin{equation*}
         \begin{split}
           \mdim(\mathcal{X}, T, \varphi) & 
           = \sup_{\mu\in \mathscr{M}^T(\mathcal{X})}
           \left(\lrdim(\mathcal{X}, T, \mathbf{d}, \mu) 
           + \int_\mathcal{X} \varphi\, d\mu\right)  \\
           & =  
           \sup_{\mu\in \mathscr{M}^T(\mathcal{X})}
           \left(\urdim(\mathcal{X}, T, \mathbf{d}, \mu) 
           + \int_\mathcal{X} \varphi\, d\mu\right).
         \end{split}
\end{equation*}
Namely Conjecture \ref{conjecture: dynamical PS} implies Conjecture 
\ref{conjecture: double variational principle} in \S \ref{subsection: main result}.
Conjecture \ref{conjecture: dynamical PS} is widely open in general.
We plan to prove it for free minimal $\mathbb{R}^d$-actions
in Part II of this series of papers.

At the end of this section 
we present a small technical result on $L^1$-mean Hausdorff dimension with potential.
This will be used in \S \ref{section: mean dimension of Z^d-actions}.

\begin{lemma} \label{lemma: mean Hausdorff dimension integer parameter}
In the definition (\ref{eq: L^1 mean Hausdorff dimension}) we can restrict the parameter $L$
to natural numbers. Namely we have
  \begin{align*}
   \umdimhl(\mathcal{X}, T, \mathbf{d},\varphi) & = \lim_{\varepsilon\to 0}\left(
   \limsup_{\substack{N\in \mathbb{N} \\ N\to \infty}} 
   \frac{\dimh\left(\mathcal{X}, \overline{\mathrm{d}}_N, \varphi_N, \varepsilon \right)}{N^d}\right) \\
   \lmdimhl(\mathcal{X}, T, \mathbf{d},\varphi) & = \lim_{\varepsilon\to 0}\left(
   \liminf_{\substack{N\in \mathbb{N} \\ N\to \infty}} 
   \frac{\dimh\left(\mathcal{X}, \overline{\mathrm{d}}_N, \varphi_N, \varepsilon \right)}{N^d}\right).
  \end{align*}
Here the parameter $N$ runs over natural numbers.
A similar result also holds for upper and lower mean Hausdorff dimensions with potential.
\end{lemma}

\begin{proof}
We prove the lower case. The upper case is similar.
By adding a positive constant to $\varphi$, we can assume that $\varphi$ is a nonnegative 
function.
Set 
\[ \lmdimhl(\mathcal{X}, T, \mathbf{d},\varphi)^\prime = \lim_{\varepsilon\to 0}\left(
   \liminf_{\substack{N\in \mathbb{N} \\ N\to \infty}} 
   \frac{\dimh\left(\mathcal{X}, \overline{\mathrm{d}}_N, \varphi_N, \varepsilon \right)}{N^d}\right). \]
It is obvious that 
$\lmdimhl(\mathcal{X}, T, \mathbf{d},\varphi) \leq 
\lmdimhl(\mathcal{X}, T, \mathbf{d},\varphi)^\prime$.
We prove the reverse inequality.
Let $L$ be a positive real number. We assume that it is sufficiently large so that 
\[  \left(\frac{L}{L-1}\right)^d  < 2. \]
Let $N:=\lfloor L\rfloor$ be the natural number not greater than $L$.
Then for any $x,y\in \mathcal{X}$ 
we have $\overline{\mathbf{d}}_N(x, y) \leq 2\, \overline{\mathbf{d}}_L(x,y)$.
We also have $\varphi_N(x) \leq \varphi_L(x)$.

Let $0<\varepsilon<1/2$. We prove that for any $s>\max \varphi_L$
\begin{equation} \label{eq: comparison of approximate Hausdorff measures}
     \mathcal{H}^{s-\frac{s}{\log \varepsilon}}_\varepsilon
      \left(\mathcal{X}, \overline{\mathbf{d}}_N, \varphi_N\right)  \leq 
      \mathcal{H}^s_{\varepsilon/2}\left(\mathcal{X},\overline{\mathbf{d}}_L, \varphi_L\right). 
\end{equation}
Indeed let $E$ be a subset of $\mathcal{X}$ with 
$\diam(E, \overline{\mathbf{d}}_L) < \varepsilon/2$.
Then 
\[  \diam(E, \overline{\mathbf{d}}_N) \leq 2\, \diam(E, \overline{\mathbf{d}}_L) 
    < \varepsilon.  \]
Moreover 
\begin{align*}
   \diam(E, \overline{\mathbf{d}}_N)^{s-\frac{s}{\log\varepsilon}-\sup_{E}\varphi_N} & \leq 
   \diam(E, \overline{\mathbf{d}}_N)^{s-\frac{s}{\log\varepsilon}-\sup_{E}\varphi_L}  \quad 
    \text{by $\varphi_N\leq \varphi_L$}\\
   & \leq 
      \left(2\, \diam(E, \overline{\mathbf{d}}_L)
      \right)^{s-\frac{s}{\log\varepsilon}-\sup_{E}\varphi_L}  \\
   & \leq \varepsilon^{-\frac{s}{\log\varepsilon}} \cdot 
    \left(2\, \diam(E, \overline{\mathbf{d}}_L) \right)^{s-\sup_{E}\varphi_L}  \\
   & = 2^{-s} \left(2\, \diam(E, \overline{\mathbf{d}}_L) \right)^{s-\sup_{E}\varphi_L}  \\
   & \leq \diam(E, \overline{\mathbf{d}}_L)^{s-\sup_{E}\varphi_L}   \quad 
    \text{by $\varphi_L\geq 0$}.
\end{align*}
Therefore we have (\ref{eq: comparison of approximate Hausdorff measures}) and hence 
\[  \dimh(\mathcal{X},\overline{\mathbf{d}}_N,\varphi_N,\varepsilon)
    \leq \left(1-\frac{1}{\log\varepsilon}\right)\dimh\left(\mathcal{X},\overline{\mathbf{d}}_L,
       \varphi_L, \frac{\varepsilon}{2}\right). \]
We divide this by $L^d$ and let $L\to \infty$ and $\varepsilon\to 0$.
Then we get $\lmdimhl(\mathcal{X}, T, \mathbf{d},\varphi)^\prime  \leq 
\lmdimhl(\mathcal{X}, T, \mathbf{d},\varphi)$.
\end{proof}

\begin{remark}  \label{remark: mean dimension integer parameter}
In the definitions (\ref{eq: definition of mean dimension with potential}) and (\ref{eq: metric mean dimension with potential}) 
of mean dimension with potential and (upper/lower) metric mean dimension with potential,
the limits with respect to $L$ exist.
Therefore we can also restrict the parameter $L$ to natural numbers when we take the limits.
\end{remark}

\section{Mean dimension of $\mathbb{Z}^d$-actions}
\label{section: mean dimension of Z^d-actions}

In this section we prepare some basic results on mean dimension theory of $\mathbb{Z}^d$-actions.
We need it in the proof of Theorem \ref{theorem: mdim is bounded by mdimh}.
This is a rather technical and indirect approach.
It is desirable to find a more direct proof of Theorem \ref{theorem: mdim is bounded by mdimh}.
However we have not found it so far\footnote{The difficulty lies in Proposition 
\ref{proposition: L^1 mean Hausdorff dimension and mean Hausdorff dimension under tame growth} below.
It is unclear for the author how to formulate and prove an analogous result for $\mathbb{R}^d$-actions.}.

The paper of Huo--Yuan \cite{Huo--Yuan} studies the variational principle for mean dimension of $\mathbb{Z}^d$-actions.
Proposition \ref{proposition: L^1 mean Hausdorff dimension and mean Hausdorff dimension under tame growth} and 
Theorem \ref{theorem: mdim is bounded by mdimh for Z^d actions} below were already mentioned in their paper
\cite[Lemma 2.12 and Lemma 2.15]{Huo--Yuan} in the case that the potential function is zero.

\subsection{Definitions of various mean dimensions for $\mathbb{Z}^d$-actions} \label{subsection: definitions of various mean dimensions
for Z^d-actions}

For a natural number $N$ we set 
\[  [N]^d = \{0,1,2,\dots, N-1\}^d. \]

Let $T\colon \mathbb{Z}^d\times \mathcal{X}\to \mathcal{X}$ be a continuous action of the group $\mathbb{Z}^d$ on a compact 
metrizable space $\mathcal{X}$.
Let $\mathbf{d}$ be a metric on $\mathcal{X}$ and $\varphi\colon \mathcal{X}\to \mathbb{R}$ a continuous function.
For a natural number $N$ we define metrics $\mathbf{d}_N$ and $\overline{\mathbf{d}}_N$ and a function $\varphi_N$ on $\mathcal{X}$ by 
\[  \mathbf{d}_N(x, y) = \max_{u\in [N]^d} \mathbf{d}(T^u x, T^u y), \quad 
     \overline{\mathbf{d}}_N(x,y) = \frac{1}{N^d}\sum_{u\in [N]^d}\mathbf{d}(T^u x, T^u y), \]
\[     \varphi_N(x) = \sum_{u\in [N]^d}\varphi(T^u x). \]
In the sequel we will sometimes consider $\mathbb{Z}^d$-actions and $\mathbb{R}^d$-actions simultaneously.
In that case we use the notations $\mathbf{d}_N^{\mathbb{Z}}, \overline{\mathbf{d}}_N^{\mathbb{Z}}, \varphi^{\mathbb{Z}}_N$
for clarifying that these quantities are defined with respect to $\mathbb{Z}^d$-actions.
(On the other hand, we will use the notations $\mathbf{d}_N^{\mathbb{R}}, \overline{\mathbf{d}}_N^{\mathbb{R}}, \varphi^{\mathbb{R}}_N$
when they are defined with respect to $\mathbb{R}^d$-actions.)

We define mean dimension with potential by 
\[  \mdim(\mathcal{X}, T, \varphi) = \lim_{\varepsilon\to 0} 
     \left(\lim_{N\to \infty} \frac{\widim_\varepsilon(\mathcal{X}, \mathbf{d}_N, \varphi_N)}{N^d}\right). \]
This is a topological invariant (i.e. independent of the choice of $\mathbf{d}$).     
We define upper/lower mean Hausdorff dimension with potential by 
\begin{align*}
   \umdimh(\mathcal{X}, T, \mathbf{d},\varphi) &= \lim_{\varepsilon\to 0} 
   \left(\limsup_{N\to \infty} \frac{\dimh(\mathcal{X},\mathbf{d}_N, \varphi_N,\varepsilon)}{N^d}\right), \\
   \lmdimh(\mathcal{X}, T, \mathbf{d},\varphi) &= \lim_{\varepsilon\to 0} 
   \left(\liminf_{N\to \infty} \frac{\dimh(\mathcal{X},\mathbf{d}_N, \varphi_N,\varepsilon)}{N^d}\right).
\end{align*}
We define upper/lower $L^1$-mean Hausdorff dimension with potential by 
\begin{align*}
   \umdimhl(\mathcal{X}, T, \mathbf{d},\varphi) &= \lim_{\varepsilon\to 0} 
   \left(\limsup_{N\to \infty} \frac{\dimh(\mathcal{X},\overline{\mathbf{d}}_N, \varphi_N,\varepsilon)}{N^d}\right), \\
   \lmdimhl(\mathcal{X}, T, \mathbf{d},\varphi) &= \lim_{\varepsilon\to 0} 
   \left(\liminf_{N\to \infty} \frac{\dimh(\mathcal{X},\overline{\mathbf{d}}_N, \varphi_N,\varepsilon)}{N^d}\right).
\end{align*}
Since $\overline{\mathbf{d}}_N(x, y) \leq \mathbf{d}_N(x,y)$, we have 
\[   \umdimhl(\mathcal{X}, T, \mathbf{d},\varphi) \leq \umdimh(\mathcal{X}, T, \mathbf{d},\varphi), \quad 
       \lmdimhl(\mathcal{X}, T, \mathbf{d},\varphi)  \leq \lmdimh(\mathcal{X}, T, \mathbf{d},\varphi). \]
We can also consider upper/lower metric mean dimension with potential for $\mathbb{Z}^d$-actions.
But we do not need them in this paper.

\subsection{Tame growth of covering numbers} \label{subsection: tame growth of covering numbers}

The purpose of this subsection is to establish a convenient sufficient condition under which 
mean Hausdorff dimension with potential and 
$L^1$-mean Hausdorff dimension with potential coincide.

The following is a key definition \cite[Condition 3]{Lindenstrauss--Tsukamoto IEEE}.

\begin{definition} \label{definition: tame growth of covering numbers}
A compact metric space $(\mathcal{X}, \mathbf{d})$ is said to have \textbf{tame growth of covering numbers}
if for any positive number $\delta$ we have 
\[   \lim_{\varepsilon\to 0} \varepsilon^\delta  \log \#\left(\mathcal{X}, \mathbf{d}, \varepsilon\right)  = 0.  \]
Recall that $\#\left(\mathcal{X}, \mathbf{d}, \varepsilon\right)$ is the minimum number $n$ such that there is an open cover 
$\mathcal{X} = U_1\cup U_2\cup \dots \cup U_n$ with $\diam U_i < \varepsilon$ for all $1\leq i \leq n$.
Notice that this is purely a condition on metric geometry. It does not involve dynamics.
\end{definition}

For example, every compact subset of the Euclidean space $\mathbb{R}^n$ has the tame growth of covering numbers
with respect to the Euclidean metric. 
The Hilbert cube $[0,1]^{\mathbb{N}}$ has the tame growth of covering numbers 
with respect to the metric 
\[  \mathbf{d}\left((x_n)_{n\in \mathbb{N}}, (y_n)_{n\in \mathbb{N}}\right)  =  \sum_{n=1}^\infty 2^{-n}|x_n-y_n|. \]

The next lemma shows that every compact metrizable space admits a metric having the tame growth of covering numbers
\cite[Lemma 3.10]{Lindenstrauss--Tsukamoto double}.

\begin{lemma}  \label{lemma: tame growth of covering numbers}
For any compact metric space $(\mathcal{X},\mathbf{d})$ there exists a metric $\mathbf{d}^\prime$ on $\mathcal{X}$
compatible with the given topology satisfying the following two conditions.
   \begin{itemize}
     \item  For all $x, y\in \mathcal{X}$ we have $\mathbf{d}^\prime(x,y) \leq \mathbf{d}(x,y)$.
     \item  The space $(\mathcal{X}, \mathbf{d}^\prime)$ has the tame growth of covering numbers. 
   \end{itemize}
\end{lemma}

\begin{proof}
Take a countable dense subset $\{x_1, x_2, x_3, \dots\}$ of $\mathcal{X}$.
We define a metric $\mathbf{d}^\prime$ by
\[  \mathbf{d}^\prime(x,y) = \sum_{n=1}^\infty 2^{-n} \left|d(x, x_n) - d(y, x_n)\right|. \]
It is easy to check that this satisfies the statement.
\end{proof}

\begin{proposition}  \label{proposition: L^1 mean Hausdorff dimension and mean Hausdorff dimension under tame growth}
Let $(\mathcal{X},\mathbf{d})$ be a compact metric space having the tame growth of covering numbers.
Let $T\colon  \mathbb{Z}^d\times \mathcal{X}\to  \mathcal{X}$ be a continuous action of the group $\mathbb{Z}^d$.
For any continuous function $\varphi\colon \mathcal{X}\to \mathbb{R}$ we have 
   \begin{align*}
      \umdimhl(\mathcal{X},T,\mathbf{d},\varphi) = \umdimh(\mathcal{X},T,\mathbf{d},\varphi), \\
      \lmdimhl(\mathcal{X},T,\mathbf{d},\varphi) = \lmdimh(\mathcal{X},T,\mathbf{d},\varphi).
   \end{align*}
\end{proposition}

\begin{proof}
The case of $d=1$ was proved in \cite[Lemma 4.3]{Tsukamoto potential}.
The following argument is its simple generalization.
We prove the lower case. The upper case is similar.

It is obvious that $\lmdimhl(\mathcal{X},T,\mathbf{d},\varphi) \leq  \lmdimh(\mathcal{X},T,\mathbf{d},\varphi)$.
We prove the reverse inequality.
By adding a positive constant to $\varphi$, we can assume that $\varphi$ is a nonnegative function.
For a finite subset $A$ of $\mathbb{Z}^d$ we define a metric $\mathbf{d}_A$ on $\mathcal{X}$ by 
\[   \mathbf{d}_A(x,y) = \max_{u\in A} \mathbf{d}(T^u x, T^u y). \]

Let $s$ be an arbitrary positive number with $\lmdimhl(\mathcal{X},T,\mathbf{d},\varphi) <s$. 
Let $0<\delta<1/2$ be arbitrary.
For each positive number $\tau$ we take an open cover $\mathcal{X} = W^\tau_1\cup \dots \cup W^\tau_{M(\tau)}$ with 
$M(\tau) = \#\left(\mathcal{X},\mathbf{d},\tau\right)$ and 
$\diam(W^\tau_i, \mathbf{d}) < \tau$ for all $1\leq i \leq M(\tau)$.
From the condition of tame growth of covering numbers, we can find $0<\varepsilon_0<1$ satisfying 
\begin{align}
  & M(\tau)^{\tau^\delta} < 2 \quad \text{for all $0<\tau<\varepsilon_0$},  \label{eq: choice of varepsilon_0 covering number} \\
  &  2^{2+(1+2\delta)s}\, \varepsilon_0^{s \delta(1-2\delta)} < 1.  \label{eq: varepsilon_0 is small}
\end{align}

Let $0<\varepsilon<\varepsilon_0$.
We have $\dimh(\mathcal{X},\overline{\mathbf{d}}_N, \varphi_N, \varepsilon) < sN^d$ for infinitely many $N$.
Pick up such an $N$.
Then there is a covering $\mathcal{X} = \bigcup_{n=1}^\infty E_n$ with $\tau_n:= \diam(E_n,\overline{\mathbf{d}}_N) < \varepsilon$
for all $n\geq 1$ and 
\begin{equation} \label{eq: choice of tau_n}
   \sum_{n=1}^\infty \tau_n^{sN^d -\sup_{E_n} \varphi_N}  < 1.
\end{equation}
Set $L_n = \tau_n^{-\delta}$. Pick $x_n\in E_n$ for each $n$.
Every $x\in E_n$ satisfies $\overline{\mathbf{d}}_N(x,x_n) \leq \tau_n$ and hence 
\[  \left|\{u\in [N]^d\mid \mathbf{d}(T^u x, T^u x_n) \geq L_n \tau_n\}\right| \leq \frac{N^d}{L_n}.  \]
Namely there is $A\subset  [N]^d$ (depending on $x$) such that $|A| \leq N^d/L_n$ and $\mathbf{d}_{[N]^d\setminus A}(x, x_n) < L_n \tau_n$.
Therefore 
\[  E_n\subset \bigcup_{\substack{A\subset [N]^d \\ |A|\leq N^d/L_n}}  B_{L_n\tau_n}(x_n, \mathbf{d}_{[N]^d\setminus A}). \]
Here $B_{L_n\tau_n}(x_n, \mathbf{d}_{[N]^d\setminus A})$ is the ball of radius $L_n \tau_n$ with respect to $\mathbf{d}_{[N]^d\setminus A}$
centered at $x_n$.
For $A = \{a_1, \dots, a_r\}\subset [N]^d$ and $1\leq i_1, \dots, i_r \leq M(\tau_n)$ we set 
\[ W(A, \tau_n, i_1, \dots,i_r) = T^{-a_1}W_{i_1}^{\tau_n}\cap \dots \cap T^{-a_r}W_{i_r}^{\tau_n}. \]  
We have 
\[  \mathcal{X} = \bigcup_{1\leq i_1, \dots, i_r\leq M(\tau_n)} W(A, \tau_n, i_1, \dots, i_r), \quad 
 (\text{here $A$ and $\tau_n$ are fixed}), \]
and hence
\[  B_{L_n\tau_n}(x_n, \mathbf{d}_{[N]^d\setminus A})  
 = \bigcup_{1\leq i_1, \dots, i_r\leq M(\tau_n)}   B_{L_n\tau_n}(x_n, \mathbf{d}_{[N]^d\setminus A}) \cap W(A, \tau_n, i_1, \dots, i_r). \]
Then 
\[  \mathcal{X} = \bigcup_{n=1}^\infty \bigcup_{\substack{A\subset [N]^d \\ r:= |A|\leq N^d/L_n}} 
     \bigcup_{1\leq i_1, \dots, i_r\leq M(\tau_n)}  E_n\cap  B_{L_n\tau_n}(x_n, \mathbf{d}_{[N]^d\setminus A}) \cap W(A, \tau_n, i_1, \dots, i_r). \]
The diameter of $E_n\cap B_{L_n\tau_n}(x_n, \mathbf{d}_{[N]^d\setminus A}) \cap W(A, \tau_n, i_1, \dots, i_r)$ 
with respect to $\mathbf{d}_N$ is less than or equal to $2L_n \tau_n = 2\tau_n^{1-\delta} < 2\varepsilon^{1-\delta}$.
Hence 
\[  \mathcal{H}^{(1+2\delta)sN^d}_{2\varepsilon^{1-\delta}}\left(\mathcal{X},\mathbf{d}_N, \varphi_N\right)
     \leq  \sum_{n=1}^\infty 2^{N^d}\, M(\tau_n)^{N^d/L_n}\left(2\tau_n^{1-\delta}\right)^{(1+2\delta)sN^d-\sup_{E_n}\varphi_N}. \]
Here the factor $2^{N^d}$ comes from the choice of $A\subset [N]^d$.
By $L_n = \tau_n^{-\delta}$ and (\ref{eq: choice of varepsilon_0 covering number})
\[  M(\tau_n)^{N^d/L_n} = \left(M(\tau_n)^{\tau_n^\delta}\right)^{N^d}  < 2^{N^d}. \]
Since $\varphi$ is a nonnegative function, 
\[ 2^{(1+2\delta)sN^d-\sup_{E_n}\varphi_N} \leq 2^{(1+2\delta)sN^d}. \]
Hence 
\[  \mathcal{H}^{(1+2\delta)sN^d}_{2\varepsilon^{1-\delta}}\left(\mathcal{X},\mathbf{d}_N, \varphi_N\right) \leq 
     \sum_{n=1}^\infty \left(2^{2+(1+2\delta)s}\right)^{N^d} \left(\tau_n^{1-\delta}\right)^{(1+2\delta)sN^d-\sup_{E_n}\varphi_N}. \]
We have
\begin{align*}
   \left(\tau_n^{1-\delta}\right)^{(1+2\delta)sN^d-\sup_{E_n}\varphi_N} &
   = \tau_n^{-\delta\left\{(1+2\delta)sN^d-\sup_{E_n}\varphi_N\right\}}\cdot 
    \tau_n^{(1+2\delta)sN^d-\sup_{E_n}\varphi_N}  \\
    &= \tau_n^{\delta\left\{(1-2\delta)sN^d+\sup_{E_n}\varphi_n\right\}}\cdot 
         \tau_n^{sN^d-\sup_{E_n}\varphi_N}.
\end{align*}         
Since $\varphi$ is nonnegative and $\tau_n < \varepsilon < \varepsilon_0< 1$
\[  \tau_n^{\delta\left\{(1-2\delta)sN^d+\sup_{E_n}\varphi_n\right\}} \leq \tau_n^{\delta(1-2\delta)sN^d}
      < \varepsilon_0^{\delta(1-2\delta)sN^d}. \]
Therefore 
\begin{align*}
    \mathcal{H}^{(1+2\delta)sN^d}_{2\varepsilon^{1-\delta}}\left(\mathcal{X},\mathbf{d}_N, \varphi_N\right) & \leq 
     \sum_{n=1}^\infty \underbrace{\left(2^{2+(1+2\delta)s} \cdot \varepsilon_0^{\delta(1-2\delta)s}
     \right)^{N^d}}_{\text{$<1$ by (\ref{eq: varepsilon_0 is small})}}
     \,  \tau_n^{sN^d-\sup_{E_n}\varphi_N} \\
     &\leq \sum_{n=1}^\infty  \tau_n^{sN^d-\sup_{E_n}\varphi_N} < 1 \quad 
     \text{by (\ref{eq: choice of tau_n})}.
\end{align*}     
Thus 
\[  \dimh(\mathcal{X},\mathbf{d}_N,\varphi_N,2\varepsilon^{1-\delta}) \leq (1+2\delta)sN^d. \]
This holds for infinitely many $N$. 
Hence
\[  \liminf_{N\to \infty} \frac{\dimh(\mathcal{X},\mathbf{d}_N,\varphi_N,2\varepsilon^{1-\delta})}{N^d} \leq 
     (1+2\delta)s. \]
Letting $\varepsilon \to 0$
\[  \lmdimh(\mathcal{X},T,\mathbf{d},\varphi) \leq  (1+2\delta)s. \]
Letting $\delta\to 0$ and $s\to \lmdimhl(\mathcal{X},T,\mathbf{d},\varphi)$, we get 
\[  \lmdimh(\mathcal{X},T,\mathbf{d},\varphi) \leq  \lmdimhl(\mathcal{X},T,\mathbf{d},\varphi). \]
\end{proof}

\subsection{$\mathbb{R}^d$-actions and $\mathbb{Z}^d$-actions}
\label{subsection: R^d-actions and Z^d-actions}

We naturally consider $\mathbb{Z}^d$ as a subgroup of $\mathbb{R}^d$.
Let $T\colon \mathbb{R}^d\times \mathcal{X}\to \mathcal{X}$ be a continuous action of $\mathbb{R}^d$ on a 
compact metrizable space $\mathcal{X}$.
We denote by $T|_{\mathbb{Z}^d}\colon \mathbb{Z}^d\times \mathcal{X}\to \mathcal{X}$ the restriction of $T$ to 
the subgroup $\mathbb{Z}^d$.
In this subsection we study relations between various mean dimensions of $T$ and $T|_{\mathbb{Z}^d}$.

Let $\mathbf{d}$ be a metric on $\mathcal{X}$ and $\varphi\colon \mathcal{X}\to \mathbb{R}$ a continuous function.

\begin{lemma} \label{lemma: mean dimension of R^d action and Z^d action}
We have:
\begin{align*}
  &  \mdim(\mathcal{X}, T, \varphi) = \mdim(\mathcal{X}, T|_{\mathbb{Z}^d}, \varphi^{\mathbb{R}}_1), \\
  & \umdimhl(\mathcal{X},T,\mathbf{d},\varphi)  =
     \umdimhl(\mathcal{X}, T|_{\mathbb{Z}^d}, \overline{\mathbf{d}}^{\mathbb{R}}_1, \varphi^{\mathbb{R}}_1), \\
  & \lmdimhl(\mathcal{X},T,\mathbf{d},\varphi)  =
     \lmdimhl(\mathcal{X}, T|_{\mathbb{Z}^d}, \overline{\mathbf{d}}^{\mathbb{R}}_1, \varphi^{\mathbb{R}}_1).
\end{align*}
Here $\overline{\mathbf{d}}^{\mathbb{R}}_1$ and $\varphi^{\mathbb{R}}_1$ are a metric and a function on $\mathcal{X}$ defined by 
\[  \overline{\mathbf{d}}^{\mathbb{R}}_1(x,y) = \int_{[0,1)^d} \mathbf{d}(T^u x, T^u y)\, du, \quad 
    \varphi^{\mathbb{R}}_1(x) = \int_{[0,1)^d}\varphi(T^u x)\, du. \]
We also have a similar result for mean Hausdorff dimension with potential by replacing $\overline{\mathbf{d}}_1^{\mathbb{R}}(x,y)$ by 
$\mathbf{d}_1^{\mathbb{R}}(x,y) = \sup_{u\in [0, 1)^d}\mathbf{d}(T^u x, T^u y)$.
\end{lemma}

\begin{proof}
Set $\rho= \overline{\mathbf{d}}^{\mathbb{R}}_1$ and $\psi = \varphi^{\mathbb{R}}_1$.
For a natural number $N$ we have 
\begin{align*}
   \overline{\rho}^{\mathbb{Z}}_N(x, y) & = \frac{1}{N^d}\sum_{u\in [N]^d}\rho(T^u x, T^u y) 
   = \frac{1}{N^d}\sum_{u\in [N]^d} \int_{v\in [0,1)^d}\mathbf{d}(T^{u+v}x, T^{u+v}y)\, dv \\
   & = \frac{1}{N^d} \int_{[0,N)^d}\mathbf{d}(T^v x, T^v y) \, dv = \overline{\mathbf{d}}^{\mathbb{R}}_N(x,y).
\end{align*}
Similarly 
\[ \psi^{\mathbb{Z}}_N(x) = \sum_{u\in [N]^d}\psi(T^u x)  = \int_{[0,N)^d} \varphi(T^v x) \, dv  = \varphi^{\mathbb{R}}_N. \]
By using Lemma \ref{lemma: mean Hausdorff dimension integer parameter}
\begin{align*}
   \lmdimhl(\mathcal{X}, T,\mathbf{d}, \varphi) & = 
   \lim_{\varepsilon\to 0} \left(\liminf_{\substack{N\in \mathbb{N} \\ N\to \infty}} 
   \frac{\dimh\left(\mathcal{X},\overline{\mathbf{d}}^{\mathbb{R}}_N, \varphi^{\mathbb{R}}_N, \varepsilon\right)}{N^d}\right)  \\
   & = \lim_{\varepsilon\to 0} \left(\liminf_{\substack{N\in \mathbb{N} \\ N\to \infty}} 
   \frac{\dimh\left(\mathcal{X},\overline{\rho}^{\mathbb{Z}}_N, \psi^{\mathbb{Z}}_N, \varepsilon\right)}{N^d}\right)  \\
   & = \lmdimhl(\mathcal{X}, T|_{\mathbb{Z}^d},\rho, \psi).
\end{align*}
We can prove the case of upper $L^1$-mean Hausdorff dimension with potential in the same way.
The case of (topological) mean dimension with potential can be also proved 
similarly by using $\left(\mathbf{d}^{\mathbb{R}}_1\right)^{\mathbb{Z}}_N = \mathbf{d}_N^{\mathbb{R}}$.
\end{proof}

\section{Mean dimension is bounded by mean Hausdorff dimension: proof of Theorem \ref{theorem: mdim is bounded by mdimh}}
\label{section: proof of mdim is bounded by mdimh}

In this section we prove Theorem \ref{theorem: mdim is bounded by mdimh}.

\subsection{A variation of the definition of mean dimension with potential}
\label{subsection: a variation of the definition of mean dimension with potential}

\begin{figure}[h]
    \centering
    \includegraphics[width=3.0in]{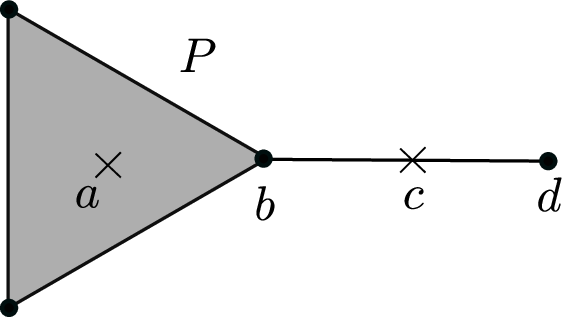}
    \caption{Here $P$ has four vertexes (denoted by dots), four $1$-dimensional simplexes and one $2$-dimensional simplex. 
    The points $b$ and $d$ are vertexes of $P$ wheres $a$ and $c$ are not. We have
    $\dim^\prime_a P =2$,  $\dim^\prime_b P =0$, $\dim^\prime_c P =1$ and $\dim^\prime_d P =0$.
    Recall $\dim_a P = \dim_b P =2$ and $\dim_c P = \dim_d P =1$.}   \label{figure: small local dimension}
\end{figure}

This subsection is a simple generalization of \cite[\S 3.2]{Tsukamoto potential}.
Here we introduce a variation of the definition of mean dimension with potential.
Let $P$ be a finite simplicial complex and $a\in P$.
We define \textbf{small local dimension} $\dim^\prime_a P$ as the minimum of $\dim \Delta$ where 
$\Delta$ is a simplex of $P$ containing $a$.
See Figure \ref{figure: small local dimension}.
(This is the same as \cite[Figure 2]{Tsukamoto potential}.)

Let $(\mathcal{X},\mathbf{d})$ be a compact metric space and $\varphi\colon \mathcal{X}\to \mathbb{R}$ a continuous function.
For $\varepsilon>0$ we set 
\begin{equation*} 
   \begin{split}
     & \widim^\prime_\varepsilon(\mathcal{X}, \mathbf{d}, \varphi)  \\ 
     & =  \inf\left\{\sup_{x\in \mathcal{X}} \left(\dim^\prime_{f(x)} P + \varphi(x)\right) \middle|
    \parbox{3in}{\centering $P$ is a finite simplicial complex and $f:\mathcal{X}\to P$ is an $\varepsilon$-embedding}\right\}.  
   \end{split}
\end{equation*}    
We also set 
\[  \mathrm{var}_\varepsilon(\varphi, \mathbf{d}) = 
    \sup\{|\varphi(x)-\varphi(y)| \, |\, \mathbf{d}(x,y) < \varepsilon\}. \]

The following lemma is given in \cite[Lemma 3.4]{Tsukamoto potential}.

\begin{lemma} \label{lemma: widim and small widim}
\[  \widim^\prime_\varepsilon(\mathcal{X}, \mathbf{d}, \varphi) 
     \leq \widim_\varepsilon(\mathcal{X}, \mathbf{d}, \varphi) 
     \leq \widim^\prime_\varepsilon(\mathcal{X}, \mathbf{d}, \varphi) + \mathrm{var}_\varepsilon(\varphi, \mathbf{d}). \]
\end{lemma}

The next lemma shows that we can use $\widim^\prime_\varepsilon(\mathcal{X}, \mathbf{d}, \varphi)$
instead of $\widim_\varepsilon(\mathcal{X}, \mathbf{d}, \varphi)$ in the definition of mean dimension with potential.

\begin{lemma}  \label{lemma: mean dimension by small local dimension}
   Let $T\colon \mathbb{Z}^d\times \mathcal{X}\to \mathcal{X}$ be a continuous action of $\mathbb{Z}^d$ on a compact metrizable space
   $\mathcal{X}$. Let $\mathbf{d}$ be a metric on $\mathcal{X}$ and $\varphi\colon \mathcal{X}\to \mathbb{R}$ a continuous function.
   Then 
  \begin{equation} \label{eq: definition of mean dimension by small local dimension}
      \mdim(\mathcal{X}, T, \varphi) = \lim_{\varepsilon\to 0} 
       \left(\lim_{N\to \infty} \frac{\widim^\prime_\varepsilon(\mathcal{X}, \mathbf{d}_N, \varphi_N)}{N^d}\right). 
   \end{equation}    
   Here the limits in the right-hand side exist as in \S \ref{subsection: mean dimension with potential of R^d-actions}.
\end{lemma}

\begin{proof}
By Lemma \ref{lemma: widim and small widim}, for any natural number $N$, we have 
\[  \widim^\prime_\varepsilon(\mathcal{X}, \mathbf{d}_N, \varphi_N) 
     \leq \widim_\varepsilon(\mathcal{X}, \mathbf{d}_N, \varphi_N) 
     \leq \widim^\prime_\varepsilon(\mathcal{X}, \mathbf{d}_N, \varphi_N) 
     + \mathrm{var}_\varepsilon(\varphi_N, \mathbf{d}_N). \]
We have 
\[  \mathrm{var}_\varepsilon(\varphi_N, \mathbf{d}_N) \leq  N^d  \mathrm{var}_\varepsilon(\varphi, \mathbf{d}). \]
Then 
\[  \lim_{\varepsilon\to 0} \left(\limsup_{N\to \infty} \frac{\mathrm{var}_\varepsilon(\varphi_N, \mathbf{d}_N)}{N^d}\right) 
      \leq  \lim_{\varepsilon \to 0} \mathrm{var}_\varepsilon(\varphi, \mathbf{d})  = 0. \]
\end{proof}

\subsection{Case of $\mathbb{Z}^d$-actions}   \label{subsection: case of Z^d-actions}

In this subsection we prove that, for $\mathbb{Z}^d$-actions, 
mean dimension with potential is bounded from above by lower mean Hausdorff dimension with potential.
A key ingredient of the proof is the following result on metric geometry.
This was proved in \cite[Lemma 3.8]{Tsukamoto potential}.

\begin{lemma}  \label{lemma: widim is bounded by using Hausdorff dimension}
Let $(\mathcal{X}, \mathbf{d})$ be a compact metric space and $\varphi\colon \mathcal{X}\to \mathbb{R}$ a continuous function.
Let $\varepsilon$ and $L$ be positive numbers, and let $s$ be a real number with $s>\max_{\mathcal{X}} \varphi$.
Suppose there exists a map $f\colon \mathcal{X}\to [0,1]^N$ such that 
   \begin{itemize}
      \item  $\norm{f(x)-f(y)}_\infty \leq L \, \mathbf{d}(x, y)$ for all $x, y\in \mathcal{X}$,
      \item  if $d(x, y) \geq \varepsilon$ then $\norm{f(x)-f(y)}_\infty =1$.
   \end{itemize}
Here $\norm{\cdot}_\infty$ is the $\ell^\infty$-norm.
Moreover we assume 
\[    4^N (L+1)^{1+s+\norm{\varphi}_\infty} \mathcal{H}_1^s\left(\mathcal{X},\mathbf{d},\varphi\right) < 1, \]
where $\norm{\varphi}_\infty = \max_{x\in \mathcal{X}} |\varphi(x)|$. 
Then we conclude that
\[    \widim^\prime_\varepsilon(\mathcal{X}, \mathbf{d},\varphi) \leq s + 1. \]
\end{lemma}

The following theorem is the main result of this subsection.

\begin{theorem}  \label{theorem: mdim is bounded by mdimh for Z^d actions}
   Let $T\colon \mathbb{Z}^d\times \mathcal{X}\to \mathcal{X}$ be a continuous action of $\mathbb{Z}^d$ on a compact metrizable space
   $\mathcal{X}$. Let $\mathbf{d}$ be a metric on $\mathcal{X}$ and $\varphi\colon \mathcal{X}\to \mathbb{R}$ a continuous function.  
   Then 
   \[  \mdim(\mathcal{X}, T, \varphi) \leq  \lmdimh(\mathcal{X}, T, \mathbf{d}, \varphi). \]
\end{theorem}

\begin{proof}
If $\lmdimh(\mathcal{X}, T, \mathbf{d}, \varphi)$ is infinite then the statement is trivial.
So we assume that it is finite.
Let $s$ be an arbitrary number larger than $\lmdimh(\mathcal{X}, T, \mathbf{d}, \varphi)$.
We prove that $\mdim(\mathcal{X}, T, \varphi) \leq s$.

Let $\varepsilon$ be a positive number.
There is a Lipschitz map $f\colon \mathcal{X} \to [0,1]^M$ such that\footnote{The construction of $f$ is as follows. 
Take a Lipschitz function $\psi\colon [0, \infty) \to [0,1]$ such that $\psi(t) = 1$ for $0\leq t \leq \varepsilon/4$ and 
$\psi(t) = 0$ for $t\geq \varepsilon/2$. Let $\{x_1, \dots, x_M\}$ be a $(\varepsilon/4)$-spanning set of $\mathcal{X}$.
Then we set $f(x) = \left(\psi(d(x, x_1)), \psi(d(x, x_2)), \dots, \psi(d(x, x_M))\right)$.}
if $\mathbf{d}(x, y) \geq  \varepsilon$ then $\norm{f(x)-f(y)}_\infty  =1$.  
Let $L$ be a Lipschitz constant of $f$. Namely we have $\norm{f(x)-f(y)}_\infty \leq L \, \mathbf{d}(x,y)$.
For each natural number $N$ we define $f_N\colon \mathcal{X}\to \left([0,1]^M\right)^{[N]^d}$ by 
\[  f_N(x) = \left(f(T^u x)\right)_{u \in [N]^d}. \]
Then we have
\begin{itemize}
   \item  $\norm{f_N(x)-f_N(y)}_\infty \leq L \, \mathbf{d}_N(x, y)$ for all $x, y \in \mathcal{X}$,
   \item  if $\mathbf{d}_N(x, y) \geq  \varepsilon$ then $\norm{f_N(x)-f_N(y)}_\infty =1$.
\end{itemize}

Let $\tau$ be an arbitrary positive number. We choose a positive number $\delta<1$ satisfying 
\begin{equation}   \label{eq: choice of delta in proof of mdim is bounded by mdimh}
    4^M (L+1)^{1+s+\tau+\norm{\varphi}_\infty} \, \delta^\tau < 1.
\end{equation}
From $\lmdimh(\mathcal{X}, T, \mathbf{d}, \varphi) < s$, there is a sequence of natural numbers 
$N_1<N_2<N_3<\dots$ such that 
\[  \dimh\left(\mathcal{X}, \mathbf{d}_{N_n}, \varphi_{N_n}, \delta\right) < sN_n^d  \quad \text{for all $n\geq 1$}. \]
Then $\mathcal{H}^{sN_n^d}_\delta\left(\mathcal{X}, \mathbf{d}_{N_n}, \varphi_{N_n}\right) < 1$ and hence 
\[  \mathcal{H}_\delta^{(s+\tau)N_n^d}\left(\mathcal{X}, \mathbf{d}_{N_n}, \varphi_{N_n}\right) \leq 
     \delta^{\tau N_n^d} \mathcal{H}^{sN_n^d}_\delta\left(\mathcal{X}, \mathbf{d}_{N_n}, \varphi_{N_n}\right) < \delta^{\tau N_n^d}. \]
Therefore 
\begin{align*}
   4^{MN_n^d}(L+1)^{1+(s+\tau)N_n^d+\norm{\varphi_{N_n}}_\infty} \, 
   \mathcal{H}_1^{(s+\tau)N_n^d}\left(\mathcal{X}, \mathbf{d}_{N_n}, \varphi_{N_n}\right) & <
   \left\{4^M (L+1)^{1+s+\tau+\norm{\varphi}_\infty} \, \delta^{\tau}\right\}^{N_n^d}  \\
   & < 1 \quad (\text{by (\ref{eq: choice of delta in proof of mdim is bounded by mdimh})}).
\end{align*}
Here we have used $\norm{\varphi_{N_n}}_\infty \leq N_n^{d} \norm{\varphi}_\infty$ in the first inequality.
Now we can use Lemma \ref{lemma: widim is bounded by using Hausdorff dimension} and conclude 
\[  \widim^\prime_{\varepsilon}\left(\mathcal{X}, \mathbf{d}_{N_n}, \varphi_{N_n}\right) \leq (s+\tau)N_n^d+1.  \]
Therefore 
\[  \lim_{N\to \infty} \frac{\widim^\prime_{\varepsilon}\left(\mathcal{X}, \mathbf{d}_{N}, \varphi_{N}\right)}{N^d} 
      = \lim_{n \to \infty} \frac{\widim^\prime_{\varepsilon}\left(\mathcal{X}, \mathbf{d}_{N_n}, \varphi_{N_n}\right)}{N_n^d}
        \leq s+\tau. \] 
Since $\tau>0$ is arbitrary, we have
\[  \lim_{N\to \infty} \frac{\widim^\prime_{\varepsilon}\left(\mathcal{X}, \mathbf{d}_{N}, \varphi_{N}\right)}{N^d}   \leq  s. \]
Thus, letting $\varepsilon \to 0$, we have $\mdim(\mathcal{X}, T, \varphi) \leq s$
by Lemma \ref{lemma: mean dimension by small local dimension}.
\end{proof}

\subsection{Proof of Theorem \ref{theorem: mdim is bounded by mdimh}}
\label{subsection: proof of Theorem mdim is bounded by mdimh}

Now we can prove Theorem \ref{theorem: mdim is bounded by mdimh}.
We write the statement again.

\begin{theorem}[$=$ Theorem \ref{theorem: mdim is bounded by mdimh}] 
Let $T\colon \mathbb{R}^d\times \mathcal{X}\to \mathcal{X}$ be a continuous action of $\mathbb{R}^d$ on a compact metrizable 
space $\mathcal{X}$. Let $\mathbf{d}$ be a metric on $\mathcal{X}$ and $\varphi\colon \mathcal{X}\to \mathbb{R}$ a continuous 
function. Then we have 
\[   \mdim(\mathcal{X}, T, \varphi) \leq   \lmdimhl\left(\mathcal{X}, T, \mathbf{d}, \varphi\right). \]
\end{theorem}

\begin{proof}
We use the notations in Lemma \ref{lemma: mean dimension of R^d action and Z^d action}.
Namely, we denote by $T|_{\mathbb{Z}^d}\colon \mathbb{Z}^d\times \mathcal{X} \to \mathcal{X}$ the restriction of $T$ to 
the subgroup $\mathbb{Z}^d$.
We define a metric $\overline{\mathbf{d}}_1^{\mathbb{R}}$ and a function $\varphi^{\mathbb{R}}_1$ on $\mathcal{X}$ by 
\[  \overline{\mathbf{d}}_1^{\mathbb{R}}(x, y) = \int_{[0,1)^d}\mathbf{d}(T^u x, T^u y)\, du, \quad 
     \varphi^{\mathbb{R}}_1(x) = \int_{[0,1)^d}\varphi(T^u x)\, du. \]
By Lemma \ref{lemma: mean dimension of R^d action and Z^d action}
\begin{align*}
   \mdim(\mathcal{X}, T, \varphi) & = \mdim\left(\mathcal{X}, T|_{\mathbb{Z}^d}, \varphi^{\mathbb{R}}_1\right) \\
   \lmdimhl\left(\mathcal{X}, T, \mathbf{d}, \varphi\right)  & 
   = \lmdimhl\left(\mathcal{X}, T|_{\mathbb{Z}^d}, \overline{\mathbf{d}}^{\mathbb{R}}_1, \varphi^{\mathbb{R}}_1\right).
\end{align*}

By Lemma \ref{lemma: tame growth of covering numbers} there exists a metric $\mathbf{d}^\prime$ on $\mathcal{X}$ such that 
$\mathbf{d}^\prime(x,y) \leq \overline{\mathbf{d}}^{\mathbb{R}}_1(x, y)$ for all $x, y\in \mathcal{X}$ and that
$(\mathcal{X}, \mathbf{d}^\prime)$ has the tame growth of covering numbers.
By $\mathbf{d}^\prime(x,y) \leq \overline{\mathbf{d}}^{\mathbb{R}}_1(x, y)$ we have 
\[  \lmdimhl\left(\mathcal{X}, T|_{\mathbb{Z}^d}, \mathbf{d}^\prime, \varphi^{\mathbb{R}}_1\right) 
     \leq   \lmdimhl\left(\mathcal{X}, T|_{\mathbb{Z}^d}, \overline{\mathbf{d}}^{\mathbb{R}}_1, \varphi^{\mathbb{R}}_1\right). \]
Since $(\mathcal{X}, \mathbf{d}^\prime)$ has the tame growth of covering numbers, 
we can apply Proposition \ref{proposition: L^1 mean Hausdorff dimension and mean Hausdorff dimension under tame growth} to it and 
get
\[  \lmdimhl\left(\mathcal{X}, T|_{\mathbb{Z}^d}, \mathbf{d}^\prime, \varphi^{\mathbb{R}}_1\right)
       =   \lmdimh \left(\mathcal{X}, T|_{\mathbb{Z}^d}, \mathbf{d}^\prime, \varphi^{\mathbb{R}}_1\right). \]
By Theorem \ref{theorem: mdim is bounded by mdimh for Z^d actions}
\[   \mdim\left(\mathcal{X}, T|_{\mathbb{Z}^d},  \varphi^{\mathbb{R}}_1\right) \leq 
      \lmdimh \left(\mathcal{X}, T|_{\mathbb{Z}^d}, \mathbf{d}^\prime, \varphi^{\mathbb{R}}_1\right). \]
Combining all the above inequalities, we conclude
\[  \mdim(\mathcal{X}, T, \varphi) \leq \lmdimhl\left(\mathcal{X}, T, \mathbf{d}, \varphi\right). \]
\end{proof}

\begin{remark}
Since $\lmdimhl\left(\mathcal{X}, T, \mathbf{d}, \varphi\right) \leq  \lmdimh\left(\mathcal{X}, T, \mathbf{d}, \varphi\right)$, 
we also have 
\[   \mdim(\mathcal{X}, T, \varphi) \leq  \lmdimh\left(\mathcal{X}, T, \mathbf{d}, \varphi\right). \]
We can directly prove this inequality (for $\mathbb{R}^d$-actions) without using $\mathbb{Z}^d$-actions.
The proof is almost the identical with the proof of Theorem \ref{theorem: mdim is bounded by mdimh for Z^d actions}.
However we have not found a direct proof Theorem \ref{theorem: mdim is bounded by mdimh} 
(the $L^1$ version) so far.
\end{remark}

\section{Mean Hausdorff dimension is bounded by rate distortion dimension: 
proof of Theorem \ref{theorem: dynamical Frostman lemma}}
\label{section: proof of dynamical Frostman lemma}

In this section we prove Theorem \ref{theorem: dynamical Frostman lemma} (dynamical Frostman’s lemma).
The proof is based on results on mutual information prepared in \S \ref{subsection: mutual information}.
Another key ingredient is the following version of Frostman’s lemma. 
This was proved in \cite[Corollary 4.4]{Lindenstrauss--Tsukamoto double}.

\begin{lemma}  \label{lemma: Frostman}
For any $0<c<1$ there exists $\delta_0 = \delta_0(c)\in (0,1)$ such that for any compact metric space 
$(\mathcal{X}, \mathbf{d})$ and any $0<\delta\leq \delta_0$ there exists a Borel probability measure $\nu$ on 
$\mathcal{X}$ satisfying 
\[  \nu(E) \leq \left(\diam \, E\right)^{c \cdot \dimh(\mathcal{X}, \mathbf{d}, \delta)} \quad 
     \text{for all Borel sets $E\subset \mathcal{X}$ with } \diam \, E <\frac{\delta}{6}. \] 
\end{lemma}

We also need the following elementary lemma.
This was proved in \cite[Appendix]{Lindenstrauss--Tsukamoto IEEE}.

\begin{lemma}  \label{lemma: construction of coupling}
Let $A$ be a finite set and $\{\mu_n\}$ a sequence of probability measures on $A$.
Suppose that $\mu_n$ converges to some probability measure $\mu$ in the weak$^*$ topology
(i.e. $\mu_n(a) \to \mu(a)$ for every $a\in A$).
Then there exists a sequence of probability measures $\pi_n$ on $A\times A$ such that 
    \begin{itemize}
       \item  $\pi_n$ is a coupling between $\mu_n$ and $\mu$, i.e., the first and second marginals of $\pi_n$ are 
                 $\mu_n$ and $\mu$ respectively,
       \item  $\pi_n$ converges to $(\mathrm{id}\times \mathrm{id})_*\mu$ in the weak$^*$ topology, namely 
        \[    \pi_n(a, b) \to \begin{cases}  \mu(a) & (\text{if $a=b$}) \\
                                                         0       &  (\text{if $a\neq  b$}) \end{cases}.  \]
    \end{itemize}
\end{lemma}

We write the statement of Theorem \ref{theorem: dynamical Frostman lemma} again.

\begin{theorem}[$=$ Theorem \ref{theorem: dynamical Frostman lemma}]
Let $T\colon \mathbb{R}^d\times \mathcal{X}\to \mathcal{X}$ be a continuous action of $\mathbb{R}^d$ on a compact metrizable space
$\mathcal{X}$. Let $\mathbf{d}$ be a metric on $\mathcal{X}$ and $\varphi\colon \mathcal{X}\to \mathbb{R}$ a continuous function.
Then we have 
\[  \umdimhl\left(\mathcal{X}, T, \mathbf{d}, \varphi\right) \leq \sup_{\mu \in \mathscr{M}^T(\mathcal{X})}
     \left(\lrdim(\mathcal{X}, T, \mathbf{d}, \mu) + \int_{\mathcal{X}} \varphi\, d\mu \right). \]
Here recall that $\mathscr{M}^T(\mathcal{X})$ is the set of all $T$-invariant Borel probability measures on $\mathcal{X}$. 
\end{theorem}

\begin{proof}
Let $c$ and $s$ be arbitrary real numbers with $0<c<1$ and $s< \umdimhl(\mathcal{X}, T, \mathbf{d}, \varphi)$.
We will construct $\mu\in \mathscr{M}^T(\mathcal{X})$ satisfying 
\begin{equation} \label{eq: requirement on mu in dynamical Frostman}
    \lrdim(\mathcal{X}, T, \mathbf{d}, \mu) + \int_{\mathcal{X}} \varphi\, d\mu \geq c s - (1-c)\norm{\varphi}_\infty. 
\end{equation}    
If this is proved then we get the claim of the theorem by letting $c\to 1$ and $s\to \umdimhl(\mathcal{X}, T, \mathbf{d}, \varphi)$.

Take $\eta>0$ with $\umdimhl(\mathcal{X}, T, \mathbf{d}, \varphi) > s+2\eta$.
Let $\delta_0 = \delta_0(c) \in (0, 1)$ be the constant introduced in Lemma \ref{lemma: Frostman}.
There are $\delta\in (0, \delta_0)$ and a sequence of positive numbers $L_1<L_2<L_3<\dots \to \infty$
satisfying 
\[  \dimh\left(\mathcal{X}, \overline{\mathbf{d}}_{L_n}, \varphi_{L_n}, \delta\right) > (s+2\eta) L_n^d \]
for all $n\geq 1$.

For a real number $t$ we set 
\[  \mathcal{X}_n(t) := \left(\frac{\varphi_{L_n}}{L_n^d}\right)^{-1}[t, t+\eta]  
      = \left\{x\in \mathcal{X}\middle|\, t\leq  \frac{\varphi_{L_n}(x)}{L_n^d} \leq t + \eta\right\}. \]

\begin{claim}  
  We can choose $t\in [-\norm{\varphi}_\infty, \norm{\varphi}_\infty]$ such that for infinitely many $n$ we have
  \[ \dimh\left(\mathcal{X}_n(t), \overline{\mathbf{d}}_{L_n}, \delta\right) \geq (s-t)L_n^d. \]
  Notice that, in particular, this inequality implies that $\mathcal{X}_n(t)$ is not empty because we assumed that $\dimh(\cdot)$ is $-\infty$ 
  for the empty set.
\end{claim}

\begin{proof}
We have\footnote{The quantity $\mathcal{H}^{(s+2\eta)L_n^d}_\delta\left(\mathcal{X}, \overline{\mathbf{d}}_{L_n}, \varphi_{L_n}\right)$
is defined only when $(s+2\eta)L_n^d >\max \varphi_{L_n}$. Therefore the following argument is problematic if we have 
$(s+2\eta)L_n^d \leq \max \varphi_{L_n}$ for all but finitely many $n$.
However, in this case, there is $t \in [-\norm{\varphi}_\infty, \norm{\varphi}_\infty]$ such that 
$t\geq s+\eta$ and $\mathcal{X}_n(t) \neq \emptyset$ for infinitely many $n$. Then we have 
$\dimh(\mathcal{X}_n(t), \overline{\mathbf{d}}_{L_n}, \delta) \geq 0 > c(s-t) L_n^d$ for infinitely many $n$ for this choice of $t$.}
$\mathcal{H}^{(s+2\eta)L_n^d}_\delta\left(\mathcal{X}, \overline{\mathbf{d}}_{L_n}, \varphi_{L_n}\right) \geq 1$.
Set $m = \lceil 2\norm{\varphi}_\infty/\eta\rceil$.
We have 
\[  \mathcal{X} = \bigcup_{\ell=0}^{m-1}\mathcal{X}_n\left(-\norm{\varphi}_\infty + \ell \eta\right). \]
Then there exists $t\in \{-\norm{\varphi}_\infty + \ell \eta\mid \ell = 0,1,\dots, m-1\}$ such that 
\[ \mathcal{H}^{(s+2\eta)L_n^d}_\delta \left(\mathcal{X}_n(t), \overline{\mathbf{d}}_{L_n}, \varphi_{L_n}\right)  \geq  \frac{1}{m} \quad
     \text{for infinitely many $n$}. \]
On the set $\mathcal{X}_n(t)$ we have 
\[ (s+2\eta)L_n^d - \varphi_{L_n} \geq (s+2\eta)L_n^d - (t+\eta)L_n^d = (s-t+\eta)L_n^d. \]     
Hence 
\begin{align*}
  \mathcal{H}^{(s+2\eta)L_n^d}_\delta \left(\mathcal{X}_n(t), \overline{\mathbf{d}}_{L_n}, \varphi_{L_n}\right) & \leq 
  \mathcal{H}^{(s-t+\eta)L_n^d}_\delta \left(\mathcal{X}_n(t), \overline{\mathbf{d}}_{L_n}\right)  \\
  & \leq \delta^{\eta L_n^d} \mathcal{H}^{(s-t)L_n^d}_\delta \left(\mathcal{X}_n(t), \overline{\mathbf{d}}_{L_n}\right).
\end{align*}
Therefore for infinitely many $n$ we have 
\[ \mathcal{H}^{(s-t)L_n^d}_\delta \left(\mathcal{X}_n(t), \overline{\mathbf{d}}_{L_n}\right)  
    \geq \frac{\delta^{-\eta L_n^d}}{m} \to \infty \quad 
    (n\to \infty). \]
Thus $\dimh\left(\mathcal{X}_n(t), \overline{\mathbf{d}}_{L_n}, \delta\right) \geq (s-t)L_n^d$ for infinitely many $n$.
\end{proof}

We fix $t\in [-\norm{\varphi}_\infty, \norm{\varphi}_\infty]$ satisfying the statement of this claim.
By choosing a subsequence  (also denoted by $L_n$) we can assume that 
\[  \dimh\left(\mathcal{X}_n(t), \overline{\mathbf{d}}_{L_n}, \delta\right) \geq (s-t)L_n^d \quad 
     \text{for all $n$}. \]

By a version of Frostman’s lemma (Lemma \ref{lemma: Frostman}), 
there is a Borel probability measure $\nu_n$ on $\mathcal{X}_n(t)$ such that 
\begin{equation}  \label{eq: choice of nu_n in the proof of dynamical Frostman}
  \nu_n(E) \leq  \left(\diam\left(E, \overline{\mathbf{d}}_{L_n}\right)\right)^{c(s-t)L_n^d} 
  \quad \text{for all Borel sets $E\subset \mathcal{X}$ with } \diam\left(E, \overline{\mathbf{d}}_{L_n}\right) < \frac{\delta}{6}. 
\end{equation}
We define a Borel probability measure $\mu_n$ on $\mathcal{X}$ by 
\[  \mu_n  = \frac{1}{L_n^d} \int_{[0, L_n)^d} T^u_*\nu_n \, du. \]
By choosing a subsequence (also denoted by $\mu_n$) we can assume that 
$\mu_n$ converges to $\mu\in \mathscr{M}^T(\mathcal{X})$ in the weak$^*$ topology.
We have 
\[
   \int_{\mathcal{X}} \varphi\, d\mu_n  =  \int_{\mathcal{X}} \frac{\varphi_{L_n}}{L_n^d}\, d\nu_n 
    = \int_{\mathcal{X}_n(t)} \frac{\varphi_{L_n}}{L_n^d}\, d\nu_n  \geq t. 
\]
Here we have used that $\nu_n$ is supported in $\mathcal{X}_n(t)$ in the second inequality and 
that $\varphi_{L_n}/L_n^d \geq t$ on $\mathcal{X}_n(t)$ in the last inequality.
Since $\mu_n \rightharpoonup  \mu$, we have 
\[   \int_{\mathcal{X}} \varphi \, d\mu  \geq  t. \]
If $t\geq  s$ then (\ref{eq: requirement on mu in dynamical Frostman}) trivially holds (recalling $0<c<1$):
\[ \lrdim(\mathcal{X},T,\mathbf{d},\mu) + \int_{\mathcal{X}} \varphi \, d\mu  \geq  t \geq cs - (1-c)\norm{\varphi}_\infty. \]
Therefore we assume $s> t$.

We will prove that for sufficiently small $\varepsilon>0$
\begin{equation} \label{eq: lower bound on rate distortion function in the proof of dynamical Frostman}
   R(\mathbf{d}, \mu, \varepsilon) \geq c(s-t) \log(1/\varepsilon) - Kc (s-t), 
\end{equation}   
where $R(\mathbf{d}, \mu, \varepsilon)$ is the rate distortion function and 
$K$ is the universal positive constant introduced in Proposition \ref{proposition: Kawabata--Dembo}.
Once this is proved, we have
\[  
  \lrdim(\mathcal{X}, T, \mathbf{d}, \mu)  = \liminf_{\varepsilon\to 0} \frac{R(\mathbf{d}, \mu, \varepsilon)}{\log(1/\varepsilon)} \geq c(s-t).
\]
Then we get (\ref{eq: requirement on mu in dynamical Frostman}) by
\[
   \lrdim(\mathcal{X}, T, \mathbf{d}, \mu) + \int_{\mathcal{X}} \varphi\, d\mu \geq c (s-t) + t
   = cs + (1-c) t \geq cs - (1-c) \norm{\varphi}_\infty.
\]
Here we have used $0<c<1$ and $t \geq -\norm{\varphi}_\infty$.
So the task is to prove (\ref{eq: lower bound on rate distortion function in the proof of dynamical Frostman}).

Let $\varepsilon$ be a positive number with $2\varepsilon \log (1/\varepsilon) < \delta/6$.
Let $M$ be a positive number, and let $X$ and $Y$ be random variables such that 
\begin{itemize}
   \item   $X$ takes values in $\mathcal{X}$ with $\mathrm{Law} X = \mu$,
   \item   $Y$ takes values in $L^1([0, M)^d, \mathcal{X})$ with 
              $\mathbb{E}\left( \int_{[0, M)^d}\mathbf{d}(T^v X, Y_v)\, dv\right) < \varepsilon\, M^d$.
\end{itemize}
We want to prove 
\[   \frac{1}{M^d} I(X;Y) \geq  c(s-t) \log(1/\varepsilon) - Kc(s-t). \]
If this is proved then we get (\ref{eq: lower bound on rate distortion function in the proof of dynamical Frostman}) and the proof is done.
We can assume that $Y$ takes only finitely many values
(Remark \ref{remark: we can assume that Y takes only finitely many values}).
We denote the set of values of $Y$ by $\mathcal{Y}$. This is a finite subset of $L^1([0, M)^d, \mathcal{X})$.

Take a positive number $\tau$ satisfying $\mathbb{E}\left( \int_{[0, M)^d}\mathbf{d}(T^v X, Y_v)\, dv\right) < (\varepsilon - 3\tau)M^d$.
We take a measurable partition 
\[  \mathcal{X} = P_1 \cup P_2 \cup \dots \cup P_\alpha \quad (\text{disjoint union}) \]
such that $\diam(P_i, \overline{\mathbf{d}}_M) < \tau$ and $\mu(\partial P_i) = 0$ for all $1\leq  i \leq \alpha$.
We pick a point $x_i \in P_i$ for each $i$ and set $A = \{x_1, \dots, x_\alpha\}$.
We define a map $\mathcal{P}\colon  \mathcal{X}\to  A$ by $\mathcal{P}(P_i) = \{x_i\}$.
Then we have 
\begin{equation}   \label{eq: average distance between P(X) and Y}
   \mathbb{E}\left( \frac{1}{M^d} \int_{[0, M)^d}\mathbf{d}\left(T^v \mathcal{P}(X), Y_v \right) \, dv\right) < \varepsilon - 2\tau.
\end{equation}
We consider the push-forward measures $\mathcal{P}_*\mu_n$ on $A$.
They converge to $\mathcal{P}_*\mu$ as $n\to \infty$ in the weak$^*$ topology by $\mu(\partial P_i) = 0$.

By Lemma \ref{lemma: construction of coupling}, we can construct random variables $X(n)$ coupled to $\mathcal{P}(X)$ such that 
$X(n)$ take values in $A$ with $\mathrm{Law} X(n) = \mathcal{P}_*\mu_n$ and 
\[ \mathbb{P}\left(X(n) = x_i, \mathcal{P}(X) = x_j\right) \to \delta_{ij} \mathbb{P}\left(\mathcal{P}(X)= x_j \right) \quad
  (n\to \infty).  \]
Then $\mathbb{E}\overline{\mathbf{d}}_M\left(X(n), \mathcal{P}(X)\right) \to 0$ as $n\to \infty$. 
We consider\footnote{This sentence is not rigorous. Strictly speaking, we can construct 
random variables $X(n), X^\prime, Y^\prime$ defined on a common probability space such that 
$\mathrm{Law}\left(X^\prime, Y^\prime\right) = \mathrm{Law}\left(\mathcal{P}(X), Y\right)$, 
\[ \mathbb{P}\left(X(n) = x_i, X^\prime = x_j\right) \to \delta_{ij} \mathbb{P}\left(X^\prime= x_j \right) \quad
  (n\to \infty),  \]
and 
\[  \mathbb{P}\left(X(n) = x_i, Y=y\middle|\, X^\prime = x_j \right) 
    =  \mathbb{P}\left(X(n) = x_i\middle|\, X^\prime = x_j \right) \cdot  \mathbb{P}\left(Y=y\middle|\, X^\prime = x_j \right). \]
    For simplicity we identify $X^\prime$ and $Y^\prime$ with $\mathcal{P}(X)$ and $Y$ respectively.}    
that $X(n)$ is coupled to $Y$ with the conditional distribution 
\[  \mathbb{P}\left(X(n) = x_i, Y=y\middle|\, \mathcal{P}(X) = x_j \right) 
    =  \mathbb{P}\left(X(n) = x_i\middle|\, \mathcal{P}(X) = x_j \right) \cdot  \mathbb{P}\left(Y=y\middle|\, \mathcal{P}(X) = x_j \right) \]
for $x_i, x_j\in A$ and $y\in \mathcal{Y}$.
Namely $X(n)$ and $Y$ are conditionally independent given $\mathcal{P}(X)$.
Then 
\begin{align*}
   \mathbb{P}\left(X(n)= x_i, Y=y\right) & = \sum_{j=1}^\alpha \mathbb{P}\left(X(n) = x_i, \mathcal{P}(X) = x_j \right) \cdot 
   \mathbb{P}\left(Y=y\middle|\, \mathcal{P}(X) = x_j \right)  \\ 
   & \to  \mathbb{P}\left(\mathcal{P}(X) = x_i, Y=y\right)  \quad   (n\to \infty). 
\end{align*}
By (\ref{eq: average distance between P(X) and Y})
\begin{equation}  \label{eq: average distance between X(n) and Y}
   \mathbb{E}\left(\frac{1}{M^d} \int_{[0, M)^d}\mathbf{d}\left(T^u X(n), Y_u\right) du \right)  <  \varepsilon - 2\tau
   \quad  \text{for large $n$}.
\end{equation}

Notice that $(X(n), Y)$ take values in a fixed finite set $A\times \mathcal{Y}$ and that their distributions converge to 
that of $\left(\mathcal{P}(X), Y\right)$.
Hence by Lemma \ref{lemma: convergence of mutual information}
\[  I\left(X(n); Y\right) \to I\left(\mathcal{P}(X); Y\right) \quad (n\to \infty). \]
We want to estimate $I\left(X(n); Y\right)$ from below.

Fix a point $x_0\in \mathcal{X}$. We will also denote by $x_0$ any constant function whose value is $x_0$.
For $x\in A$ and $y\in L^1([0, M)^d, \mathcal{X})$ we define a conditional probability mass function by 
\[  \rho_n(y|x) = \mathbb{P}\left(Y=y \middle|\, X(n) = x\right).  \]
This is nonzero only for $y\in \mathcal{Y}$.
(Here $\rho_n(\cdot|x)$ may be an arbitrary probability measure on $\mathcal{Y}$ if $\mathbb{P}\left(X(n) = x \right)=0$.)

We define $\Lambda\subset \mathbb{R}^d$ by 
\[  \Lambda = \left\{\left(M m_1, M m_2, \dots, M m_d\right) \middle|\, m_k\in \mathbb{Z}, \> 
                     0 \leq  m_k \leq \frac{L_n}{M}-2 \> (1\leq k \leq d) \right\}. \]
Let $v \in [0,M)^d$. 
We have
\[   \bigcup_{\lambda\in \Lambda} \left(v+\lambda+ [0, M)^d\right)  \subset  [0, L_n)^d. \]
Notice that the left-hand side is a disjoint union.
Here $v+\lambda+[0, M)^d = \{v+\lambda+w \mid  w\in [0,M)^d\}$.
Set 
\[  E_v = [0, L_n)^d \setminus  \bigcup_{\lambda\in \Lambda} \left(v+\lambda+ [0, M)^d\right). \]
See Figure \ref{figure: squares}.

\begin{figure}[h] 
    \centering
    \includegraphics[width=3.0in]{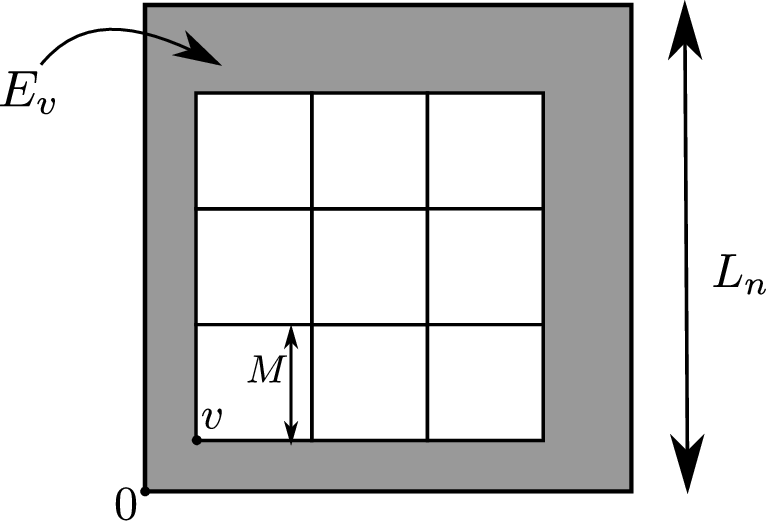}
    \caption{The big square is $[0, L_n)^d$ and small squares are $v+\lambda + [0, M)^d$
    ($\lambda \in \Lambda$).   The shadowed region is $E_v$.}    \label{figure: squares}
\end{figure}

For $x\in \mathcal{X}$ and $f\in L^1\left([0, L_n)^d, \mathcal{X}\right)$
we define a conditional probability mass function $\sigma_{n, v}(f|x)$ by 
\[  \sigma_{n, v}(f|x) =  \delta_{x_0} \left(f|_{E_v}\right) \cdot 
     \prod_{\lambda\in \Lambda} \rho_n\left(f|_{v+\lambda+[0,M)^d}\middle|\, \mathcal{P}(T^{v+\lambda}x)\right). \]
Here $f|_{E_v}$ is the restriction of $f$ to $E_v$ and $\delta_{x_0}$ is the delta probability measure concentrated at the 
constant function $x_0 \in L^1(E_v, \mathcal{X})$.
We naturally consider $f|_{v+\lambda+[0,M)^d}$ (the restriction of $f$ to $v+\lambda+[0,M)^d$) as an element of 
$L^1\left([0, M)^d, \mathcal{X}\right)$.

We define a transition probability $\sigma_n$ on $\mathcal{X}\times  L^1([0, L_n)^d, \mathcal{X})$ by 
\[  \sigma_n(B \mid x) =  \frac{1}{M^d}\int_{[0, M)^d} \sigma_{n, v}(B \mid x)\, dv \]
for $x\in \mathcal{X}$ and Borel subsets $B \subset L^1([0, L_n)^d, \mathcal{X})$.
Here $\sigma_{n, v}(B \mid  x) = \sum_{f\in B} \sigma_{n, v}(f|x)$.
(Notice that for each $x\in \mathcal{X}$ the function $\sigma_{n, v}(f|x)$ is nonzero only for finitely many $f$.)

Take random variables $Z$ and $W$ such that $Z$ takes values in $\mathcal{X}$ with $\mathrm{Law} Z = \nu_n$ 
and that $W$ takes values in $L^1([0, L_n)^d, \mathcal{X})$ with 
\[  \mathbb{P}\left(W \in B\middle|\, Z=x\right) = \sigma_n(B\mid x) \quad 
     (x\in \mathcal{X}, B\subset L^1([0, L_n)^d, \mathcal{X})). \]
Notice that $Z$ and $W$ depend on $n$. 
Rigorously speaking, we should denote them by $Z^{(n)}$ and $W^{(n)}$.
However we suppress their dependence on $n$ in the notations for simplicity.

We estimate $\mathbb{E}\left(\frac{1}{L_n^d}\int_{[0, L_n)^d}\mathbf{d}(T^u Z, W_u)\, du\right)$ and $I(Z;W)$.

\begin{claim}  \label{claim: distance between Z and W}
For all sufficiently large $n$ we have
\[   \mathbb{E}\left(\frac{1}{L_n^d}\int_{[0, L_n)^d}\mathbf{d}(T^u Z, W_u)\, du\right)  < \varepsilon. \]
\end{claim}

\begin{proof}
For each $v\in [0, M)^d$ we take a random variable $W(v)$ coupled to $Z$ such that 
$W(v)$ takes values in $L^1([0, L_n)^d, \mathcal{X})$ with 
$\mathbb{P}\left(W(v) = f\mid Z=x\right) = \sigma_{n, v}(f|x)$ for $x\in \mathcal{X}$ and 
$f\in  L^1([0, L_n)^d, \mathcal{X})$.
Then
\[   \mathbb{E}\int_{[0, L_n)^d} \mathbf{d}(T^u Z, W_u)\, du  =  
      \frac{1}{M^d} \int_{[0, M)^d} \mathbb{E}\left(\int_{[0, L_n)^d}\mathbf{d}(T^u Z, W(v)_u)\, du\right) \, dv.  \]
For each $v\in [0, M)^d$ we have 
\begin{align*}
   & \mathbb{E}\int_{[0, L_n)^d}\mathbf{d}(T^u Z, W(v)_u)\, du   \\ & = 
   \mathbb{E}\int_{E_v}\mathbf{d}(T^u Z, W(v)_u)\, du + \sum_{\lambda\in \Lambda}
   \mathbb{E}\int_{v+\lambda+[0, M)^d}\mathbf{d}(T^u Z, W(v)_u)\, du \\
   & \leq   C L_n^{d-1}  + \sum_{\lambda\in \Lambda}
   \mathbb{E}\int_{v+\lambda+[0, M)^d}\mathbf{d}(T^u Z, W(v)_u)\, du, 
\end{align*}
where $C$ is a positive constant independent of $v, L_n$.
In the last inequality we have used $\mathbf{m}(E_v) \leq \mathrm{const}\cdot L_n^{d-1}$.
Since $\overline{\mathbf{d}}_M(x, \mathcal{P}(x)) < \tau$ for every $x\in \mathcal{X}$, we have 
\begin{align*}
   & \mathbb{E}\int_{v+\lambda+[0, M)^d}\mathbf{d}(T^u Z, W(v)_u)\, du  \\
    &\leq M^d \tau  + \mathbb{E}\int_{[0, M)^d} \mathbf{d}\left(T^u \mathcal{P}(T^{v+\lambda}Z), W(v)_{v+\lambda+u}\right) du \\
    & \leq   M^d\tau  + \sum_{f\in \mathcal{Y}} 
    \int_{[0, M)^d}\left(\int_{\mathcal{X}}\mathbf{d}(T^u x, f_u)\rho_n(f|x) \, d\left(\mathcal{P}_*T^{v+\lambda}_*\nu_n(x)\right)\right)du.
\end{align*}
We sum up these estimates over $\lambda\in \Lambda$. 
Noting $M^d |\Lambda| \leq L_n^d$, we have 
\begin{align*}
    & \mathbb{E} \int_{[0, L_n)^d}\mathbf{d}(T^u Z, W(v)_u)\, du \\
    &  \leq  C L_n^{d-1} + \tau L_n^d 
        + \sum_{\substack{\lambda\in \Lambda \\ f\in \mathcal{Y}}} 
        \int_{[0, M)^d} \left(\int_{\mathcal{X}}\mathbf{d}(T^u x, f_u)\rho_n(f|x) d\left(\mathcal{P}_*T^{v+\lambda}_*\nu_n(x)\right)\right)du. 
\end{align*}        
We integrate this over $v\in [0, M)^d$.  Note that\footnote{For two measures $\mathbf{m}_1$ and $\mathbf{m}_2$ on $\mathcal{X}$ 
we write $\mathbf{m}_1\leq \mathbf{m}_2$ if we have $\mathbf{m}_1(B) \leq  \mathbf{m}_2(B)$ for all Borel subsets 
$B\subset \mathcal{X}$.}
\begin{align*}
     \int_{[0, M)^d} \left(\sum_{\lambda \in \Lambda} \mathcal{P}_*T^{v+\lambda}_*\nu_n\right)  dv  &=  
     \sum_{\lambda\in \Lambda} \int_{\lambda+[0,M)^d}  \mathcal{P}_*T^{v}_*\nu_n\, dv   \\ 
     & \leq 
     \int_{[0, L_n)^d}\mathcal{P}_*T^v_*\nu_n \, dv  = L_n^d \mathcal{P}_*\mu_n.
\end{align*}
Hence we have 
\begin{align*}
  & \frac{1}{M^d} \int_{[0, M)^d} \mathbb{E}\left(\int_{[0, L_n)^d}\mathbf{d}(T^u Z, W(v)_u)\, du\right) \, dv  \\
  & \leq  CL_n^{d-1} + \tau L_n^d + \frac{L_n^d}{M^d} \sum_{f\in \mathcal{Y}} 
   \int_{[0,M)^d}\left(\int_{\mathcal{X}} \mathbf{d}(T^u x, f_u)\rho_n(f|x) d\mathcal{P}_*\mu_n(x)\right) du \\
  & =  CL_n^{d-1} + \tau L_n^d  +  \frac{L_n^d}{M^d} \mathbb{E}\left(\int_{[0, M)^d}\mathbf{d}\left(T^u X(n), Y_u\right) du \right).
\end{align*}
In the last equality we have used $\mathrm{Law} X(n) = \mathcal{P}_*\mu_n$ and 
$\rho_n(f|x) = \mathbb{P}(Y=f\mid X(n) = x)$.
Therefore 
\[  \mathbb{E}\left(\frac{1}{L_n^d} \int_{[0, L_n)^d} \mathbf{d}(T^u Z, W_u)\, du \right)
     \leq  \frac{C}{L_n} + \tau 
    +  \mathbb{E}\left(\frac{1}{M^d} \int_{[0, M)^d}\mathbf{d}\left(T^u X(n), Y_u\right) du \right). \]
The third term in the right-hand side is smaller than $\varepsilon - 2\tau$ for large $n$ 
by (\ref{eq: average distance between X(n) and Y}).
Therefore we have
\[ \mathbb{E}\left(\frac{1}{L_n^d} \int_{[0, L_n)^d} \mathbf{d}(T^u Z, W_u)\, du \right) < \varepsilon \]
for all sufficiently large $n$.
\end{proof}

\begin{claim}  \label{claim: mutual information between Z and W}
  \[  \frac{1}{L_n^d}I(Z;W) \leq \frac{1}{M^d} I\left(X(n); Y\right). \]
\end{claim}

\begin{proof}
We have $I(Z;W) = I(\nu_n, \sigma_n)$.
Since $\sigma_n = (1/M^d) \int_{[0, M)^d} \sigma_{n, v} \, dv$, 
we apply to it Proposition \ref{proposition: concavity and convexity of I(mu, nu)} (2)
(the convexity of mutual information in transition probability):
\[  I(\nu_n, \sigma_n) \leq \frac{1}{M^d} \int_{[0, M)^d} I(\nu_n, \sigma_{n, v})\, dv. \]
By Lemma \ref{lemma: subadditivity of mutual information} (subadditivity of mutual information)\footnote{Let $W(v)$ be the random 
variable introduced in Claim \ref{claim: distance between Z and W}.
Then $I(\nu_n, \sigma_{n, v}) = I(Z;W(v))$. Consider the restrictions $W(v)|_{v+\lambda+[0, M)^d}$ $(\lambda\in \Lambda)$
and $W(v)|_{E_v}$. From the definition of the measure $\sigma_{n, v}$, they are conditionally independent given $Z$. 
By Lemma \ref{lemma: subadditivity of mutual information}
\[  I(Z;W) \leq  I\left(Z;W(v)|_{E_v}\right) + \sum_{\lambda\in \Lambda} I\left(Z;W(v)|_{v+\lambda+[0, M)^d}\right). \]
$I\left(Z;W(v)|_{E_v}\right) =0$ because $W(v)|_{E_v}$ is constantly equal to $x_0$. 
We have 
\[  I\left(Z;W(v)|_{v+\lambda+[0, M)^d}\right) = I\left(\mathcal{P}(T^{v+\lambda} Z); W(v)|_{v+\lambda+[0, M)^d}\right)
     = I\left(\mathcal{P}_*T^{v+\lambda}_*\nu_n, \rho_n\right). \]}
\[  I(\nu_n, \sigma_{n, v}) \leq \sum_{\lambda\in \Lambda} I\left(\mathcal{P}_*T^{v+\lambda}_*\nu_n, \rho_n\right). \]
Hence 
\begin{align*}
    I(\nu_n, \sigma_n) & \leq \frac{1}{M^d} \sum_{\lambda\in \Lambda} \int_{[0, M)^d}
     I\left(\mathcal{P}_*T^{\lambda+v}_*\nu_n, \rho_n\right) dv \\
     & =  \frac{1}{M^d}\int_{\cup_{\lambda\in \Lambda} \left(\lambda + [0, M)^d\right)}I\left(\mathcal{P}_*T^v_*\nu_n, \rho_n\right) dv \\
    & \leq \frac{1}{M^d}\int_{[0, L_n)^d}I\left(\mathcal{P}_*T^v_*\nu_n, \rho_n\right) dv
\end{align*}    
By Proposition \ref{proposition: concavity and convexity of I(mu, nu)} (1) 
(the concavity of mutual information in probability measure)
\begin{align*}
    \frac{1}{L_n^d} \int_{[0, L_n)^d}I\left(\mathcal{P}_*T^v_*\nu_n, \rho_n\right) dv  & \leq 
    I\left(\frac{1}{L_n^d}\int_{[0, L_n)^d}\mathcal{P}_*T^v_*\nu_n\, dv, \rho_n\right) \\
   & =  I\left(\mathcal{P}_*\mu_n, \rho_n\right) \\
   & = I(X(n); Y).
\end{align*}
Therefore we conclude 
\[  I(Z;W) =  I(\nu_n, \sigma_n) \leq  \frac{L_n^d}{M^d} I(X(n); Y). \]
\end{proof}

We define a metric $D_n$ on $L^1\left([0,L_n)^d, \mathcal{X}\right)$ by 
\[  D_n(f, g)  = \frac{1}{L_n^d} \int_{[0, L_n)^d} \mathbf{d}\left(f(u), g(u)\right) du. \]
Then the map 
\[  F_n\colon \left(\mathcal{X}, \overline{\mathbf{d}}_{L_n}\right) \ni x \mapsto   \left(T^t x\right)_{t\in [0, L_n)^d}
     \in \left(L^1\left([0,L_n)^d, \mathcal{X}\right), D_n\right) \]
is an isometric embedding.
Consider the push-forward measure ${F_n}_*\nu_n$ on $L^1\left([0,L_n)^d, \mathcal{X}\right)$.
It follows from (\ref{eq: choice of nu_n in the proof of dynamical Frostman}) that
\[  {F_n}_*\nu_n\left(E\right) \leq \left(\diam(E, D_n)\right)^{c(s-t)L_n^d}  \] 
for all Borel subsets $E\subset  L^1\left([0,L_n)^d, \mathcal{X}\right)$ with $\diam (E, D_n) < \delta/6$.

We have $\mathrm{Law} F_n(Z) = {F_n}_*\nu_n$, and by Claim \ref{claim: distance between Z and W}
\[  \mathbb{E}\left( D_n(F_n(Z), W) \right) < \varepsilon  \quad  \text{for large $n$}. \]
Since $2\varepsilon \log (1/\varepsilon) < \delta/6$, 
we apply Proposition \ref{proposition: Kawabata--Dembo} (Kawabata--Dembo estimate)
to $(F_n(Z), W)$ and get 
\[  I(Z;W) = I\left(F_n(Z); W\right) \geq c(s-t)L_n^d \log (1/\varepsilon) - K \left(c(s-t)L_n^d+1\right) \]
for large $n$. Here $K$ is a universal positive constant.
By Claim \ref{claim: mutual information between Z and W}
\[  \frac{1}{M^d} I\left(X(n); Y\right) \geq \frac{1}{L_n^d} I\left(Z;W\right) 
     \geq   c(s-t) \log(1/\varepsilon) - K\left(c(s-t) + \frac{1}{L_n^d}\right)  \]
for large $n$.    
Since $I\left(X(n); Y\right) \to  I\left(\mathcal{P}(X); Y\right)$ as $n\to \infty$, we get 
\[ \frac{1}{M^d}I\left(\mathcal{P}(X); Y\right) \geq  c(s-t) \log (1/\varepsilon) - K c (s-t). \]
Then we have 
\[ \frac{1}{M^d}I(X; Y) \geq  \frac{1}{M^d}I\left(\mathcal{P}(X); Y\right) \geq  c(s-t) \log (1/\varepsilon) - K c (s-t). \]
This is what we want to prove.
\end{proof}

Now we have proved Theorems \ref{theorem: mdim is bounded by mdimh} and \ref{theorem: dynamical Frostman lemma}.
These two theorems implies Theorem \ref{main theorem} (Main theorem) as we already explained in 
\S \ref{section: metric mean dimension with potential and mean Hausdorff dimension with potential}.
Therefore we have completed the proof of Theorem \ref{main theorem}.

\section{Local nature of metric mean dimension with potential}
\label{section: local nature of metric mean dimension with potential}

This section is independent of the proof of Theorem \ref{main theorem} (Main Theorem).
It can be read independently of Sections 2, 4, 5 and 6.
Here we present a formula expressing metric mean dimension with potential by using a certain \textit{local} quantity.
We plan to use it in a future study of geometric examples of dynamical systems \cite{Gromov, Tsukamoto Brody curves}.

\subsection{A formula of metric mean dimension with potential}  \label{subsection: a formula of metric mean dimension with potential}

Let $T\colon \mathbb{R}^d\times \mathcal{X}\to \mathcal{X}$ be a continuous action of $\mathbb{R}^d$ on a compact metrizable space 
$\mathcal{X}$.
Let $\mathbf{d}$ be a metric on $\mathcal{X}$ and $\varphi\colon \mathcal{X}\to \mathbb{R}$ a continuous function.
For a subset $E\subset \mathcal{X}$ and $\varepsilon>0$ we set 
\[ P_T(E, \mathbf{d}, \varphi, \varepsilon) = \liminf_{L\to \infty} \frac{\log \#\left(E, \mathbf{d}_L, \varphi_L, \varepsilon\right)}{L^d}. \]
Here recall that
\[  \#\left(E, \mathbf{d}_L, \varphi_L, \varepsilon\right)  =
    \inf\left\{\sum_{i=1}^n (1/\varepsilon)^{\sup_{U_i}\varphi} \middle|\, 
    \parbox{3in}{\centering  $E \subset  U_1\cup \dots \cup U_n$. Each $U_i$ is an open set of $\mathcal{X}$
    with $\diam(U_i, \mathbf{d}_L) < \varepsilon$.} \right\}. \]
Also recall that the upper and lower metric mean dimensions with potential are defined by 
\begin{align*}
     \umdimm(\mathcal{X}, T, \mathbf{d}, \varphi) & = \limsup_{\varepsilon\to 0} \frac{P_T(\mathcal{X}, \mathbf{d}, \varphi, \varepsilon)}{\log (1/\varepsilon)},  \\
     \lmdimm(\mathcal{X}, T, \mathbf{d}, \varphi)   & = \liminf_{\varepsilon\to 0} \frac{P_T(\mathcal{X}, \mathbf{d}, \varphi, \varepsilon)}{\log (1/\varepsilon)}.   
\end{align*}

For a (not necessarily bounded) subset $A$ of $\mathbb{R}^d$ we define a metric $\mathbf{d}_A$ on $\mathcal{X}$ by 
\[  \mathbf{d}_A(x, y) = \sup_{a\in A} \mathbf{d}\left(T^a x, T^a y\right). \]
(If $A$ is unbounded, this metric is not compatible with the given topology of $\mathcal{X}$ in general.)
For $x\in  \mathcal{X}$ and $\delta>0$ we define $B_\delta(x, \mathbf{d}_{\mathbb{R}^d})$ as the closed $\delta$-ball with respect to $\mathbf{d}_{\mathbb{R}^d}$:
\[ B_\delta(x, \mathbf{d}_{\mathbb{R}^d}) = \{y\in \mathcal{X}\mid  \mathbf{d}_{\mathbb{R}^d}(x, y) \leq \delta\}
    = \{y\in \mathcal{X}\mid  \mathbf{d}(T^u x, T^u y) \leq \delta  \> (\forall  u\in \mathbb{R}^d)\}. \]

The following is the main result of this section.

\begin{theorem}  \label{theorem: local formula of metric mean dimension with potential}
For any $\delta>0$ we have 
      \begin{align*}
       \umdimm(\mathcal{X}, T, \mathbf{d}, \varphi) & = \limsup_{\varepsilon\to 0} 
      \frac{\sup_{x\in \mathcal{X}} P_T\left(B_\delta(x, \mathbf{d}_{\mathbb{R}^d}), \mathbf{d}, \varphi, \varepsilon\right)}{\log (1/\varepsilon)}, \\
      \lmdimm(\mathcal{X}, T, \mathbf{d}, \varphi) & = \liminf_{\varepsilon\to 0} 
      \frac{\sup_{x\in \mathcal{X}} P_T\left(B_\delta(x, \mathbf{d}_{\mathbb{R}^d}), \mathbf{d}, \varphi, \varepsilon\right)}{\log (1/\varepsilon)}.
      \end{align*}
\end{theorem}

Notice that $B_\delta(x, \mathbf{d}_{\mathbb{R}^d})$ is not a neighborhood of $x$ with respect to the original metric $\mathbf{d}$ in general.
Nevertheless we can calculate the metric mean dimension with potential by gathering such information.

In the case that $\varphi$ is identically zero, Theorem \ref{theorem: local formula of metric mean dimension with potential} was proved in 
\cite{Tsukamoto local nature}.
The proof of Theorem \ref{theorem: local formula of metric mean dimension with potential} follows the argument of \cite{Tsukamoto local nature}, 
which is in turn based on the paper of Bowen \cite{Bowen entropy expansive}.

\subsection{Tiling argument}  \label{subsection: tiling argument}

Here we prepare a technical lemma (Lemma \ref{lemma: tiles} below).
For $x = (x_1, \dots, x_d)\in \mathbb{R}^d$ we set $\norm{x}_\infty = \max_{1\leq i \leq d} |x_i|$.
A \textbf{cube} of $\mathbb{R}^d$ is a set $\Lambda$ of the form 
\[ \Lambda = u + [0, L]^d = \{u+v\mid v\in [0, L]^d\}, \]
where $u\in \mathbb{R}^d$ and $L>0$. We set $\ell(\Lambda)  = L$.
For $r>0$ and $A\subset \mathbb{R}^d$ we define 
\[ \partial (A, r) = \left\{x \in \mathbb{R}^d\mid \exists y\in A, \exists z\in \mathbb{R}^d\setminus A: \norm{x-y}_\infty \leq r \text{ and } \norm{x-z}_\infty \leq r \right\}, \]
\[ B_r(A)  = A\cup \partial(A, r) = \{x\in \mathbb{R}^d\mid \exists y\in A: \norm{x-y}_\infty \leq  r\}. \]

For a finite set $\mathcal{C} = \{\Lambda_1, \dots, \Lambda_n\}$ of cubes of $\mathbb{R}^d$ we set 
\[  \ell_{\min}(\mathcal{C}) = \min_{1\leq i \leq n} \ell(\Lambda_i), \quad 
     \ell_{\max}(\mathcal{C})  = \max_{1\leq  i \leq n} \ell(\Lambda_i). \]

The following lemma was proved in \cite[Proposition 3.4]{Tsukamoto local nature}.

\begin{lemma} \label{lemma: tiles}
For any $\eta>0$ there exists a natural number $k_0 = k_0(\eta)>0$ for which the following statement holds.
Let $A$ be a bounded Borel subset of $\mathbb{R}^d$.
Let $\mathcal{C}_k$ $(1\leq  k  \leq k_0)$ be finite sets of cubes of $\mathbb{R}^d$ such that  
   \begin{enumerate}
    \item  $\ell_{\max}(\mathcal{C}_1) \geq 1$ and $\ell_{\min}(\mathcal{C}_{k+1}) \geq k_0\cdot \ell_{\max}(\mathcal{C}_k)$ 
    for all $1\leq  k \leq  k_0-1$,
    \item  $\mathbf{m}\left(\partial (A, \ell_{\max}(\mathcal{C}_{k_0}))\right) < (\eta/3)\cdot \mathbf{m}(A)$, 
    \item  $A \subset \bigcup_{\Lambda\in \mathcal{C}_k} \Lambda$ for every $1\leq  k  \leq  k_0$.
   \end{enumerate}
 Then there is a disjoint subfamily $\mathcal{C}^\prime \subset \mathcal{C}_1\cup \dots \cup \mathcal{C}_{k_0}$ satisfying 
  \[  \bigcup_{\Lambda\in \mathcal{C}^\prime}\Lambda \subset A, \quad   
       \mathbf{m}\left(B_1\left(A \setminus \bigcup_{\Lambda\in \mathcal{C}^\prime} \Lambda\right)\right) < \eta \cdot \mathbf{m}(A). \]
Here ``disjoint" means that for any two distinct $\Lambda_1, \Lambda_2\in \mathcal{C}^\prime$ we have $\Lambda_1\cap \Lambda_2 = \emptyset$.
\end{lemma}

This is a rather technical statement.
The assumption (1) means that some cube of $\mathcal{C}_1$ is not so small and that
every cube in $\mathcal{C}_{k+1}$ is much larger than cubes in $\mathcal{C}_k$.
The assumption (2) means that $A$ is much larger than all the given cubes.
The assumption (3) means that each $\mathcal{C}_k$ covers $A$.
The conclusion means that we can find a disjoint subfamily $\mathcal{C}^\prime$ which covers a substantial portion of $A$.

\subsection{The case that $\varphi$ is nonnegative} \label{subsection: the case that varphi is nonnegative}

This subsection is also a preparation for the proof of Theorem \ref{theorem: local formula of metric mean dimension with potential}.
Let $T\colon \mathbb{R}^d\times \mathcal{X}\to \mathcal{X}$ be a continuous action of $\mathbb{R}^d$ on a compact metrizable space
$\mathcal{X}$. Let $\mathbf{d}$ be a metric on $\mathcal{X}$ and $\varphi\colon \mathcal{X} \to \mathbb{R}$ a continuous function.
Throughout this subsection, we assume that $\varphi$ is a nonnegative function.

Recall that for a bounded Borel subset $A\subset \mathbb{R}^d$ a new function $\varphi_A\colon \mathcal{X}\to \mathbb{R}$ is defined by 
\[  \varphi_A(x) = \int_A \varphi(T^u x) \, du. \]

\begin{lemma} \label{lemma: block coding}
Let $0<\varepsilon <1$ and $E\subset \mathcal{X}$. Let $A, A_1, A_2, \dots, A_n$ be bounded Borel subsets of $\mathbb{R}^d$.
If $A\subset A_1 \cup A_2 \cup \dots \cup A_n$ then 
\[   \#\left(E, \mathbf{d}_A, \varphi_A, \varepsilon\right) \leq  \prod_{k=1}^n \#\left(E, \mathbf{d}_{A_k}, \varphi_{A_k}, \varepsilon\right). \]
\end{lemma}

\begin{proof}
Suppose we are given an open cover $E\subset U_{k1}\cup \dots \cup U_{k m_k}$ with $\diam (U_{k j}, \mathbf{d}_{A_k}) < \varepsilon$
for each $1\leq k \leq n$.
Then 
\[  E\subset  \bigcup
     \left\{ U_{1 j_1}\cap U_{2 j_2} \cap \dots \cap U_{n j_n}  \middle|\, 1\leq j_1 \leq m_1, 1 \leq j_2 \leq m_2, \dots, 1\leq j_n \leq m_n \right\}.   \]
From $A\subset A_1\cup \dots \cup A_n$, the diameter of $U_{1 j_1}\cap U_{2 j_2} \cap \dots \cap U_{n j_n}$ is smaller than 
$\varepsilon$ with respect to the metric $\mathbf{d}_A$.
Since $\varphi$ is nonnegative (here we use this assumption), we have 
\[  \varphi_A \leq \varphi_{A_1} + \varphi_{A_2}+ \dots + \varphi_{A_n} \]
and hence 
\[  \sup_{U_{1 j_1} \cap \dots \cap U_{n j_n}} \varphi_{A} \leq  \sum_{k=1}^n \sup_{U_{k j_k}} \varphi_{A_k}. \]
Therefore we have
\[   \sum_{\substack{1\leq j_1 \leq m_1 \\ \svdots \\ 1\leq j_n \leq m_n}}
      \left(\frac{1}{\varepsilon}\right)^{\sup_{U_{1 j_1} \cap \dots \cap U_{n j_n}} \varphi_{A}} 
      \leq  \prod_{k=1}^n \left(\sum_{j=1}^{m_k} \left(\frac{1}{\varepsilon}\right)^{\sup_{U_{k j}} \varphi_{A_k}}\right). \]
This proves the claim of the lemma.
\end{proof}

\begin{lemma}  \label{lemma: crude estimate of covering number}
For $0< \varepsilon < 1$ and a bounded Borel subset $A\subset \mathbb{R}^d$, we have
\[  \#\left(\mathcal{X}, \mathbf{d}_A, \varphi_A, \varepsilon\right) \leq 
      \{\#\left(\mathcal{X}, \mathbf{d}_{[0, 1]^d}, \varphi_{[0, 1]^d}, \varepsilon\right)\}^{\mathbf{m}\left(B_1(A)\right)}. \]
Notice that $\#\left(\mathcal{X}, \mathbf{d}_{[0, 1]^d}, \varphi_{[0, 1]^d}, \varepsilon\right) \geq 1$
because $0<\varepsilon <1$ and $\varphi$ is nonnegative.
\end{lemma}

\begin{proof}
Let $\Omega$ be the set of $u\in \mathbb{Z}^d$ with $\left(u+[0, 1]^d\right)\cap A \neq \emptyset$.
We have 
\[  A \subset \bigcup_{u\in \Omega} \left(u+[0,1]^d\right) \subset B_1(A). \]
In particular the cardinality of $\Omega$ is bounded from above by $\mathbf{m}\left(B_1(A)\right)$.
Then by Lemma \ref{lemma: block coding}
\begin{align*}
    \#\left(\mathcal{X}, \mathbf{d}_A, \varphi_A, \varepsilon\right) & \leq 
    \prod_{u\in \Omega} \#\left(\mathcal{X}, \mathbf{d}_{u+[0,1]^d}, \varphi_{u+[0,1]^d}, \varepsilon\right) \\
    & = \prod_{u\in \Omega} \#\left(\mathcal{X}, \mathbf{d}_{[0,1]^d}, \varphi_{[0,1]^d}, \varepsilon\right)  \\
    & \leq        \{\#\left(\mathcal{X}, \mathbf{d}_{[0, 1]^d}, \varphi_{[0, 1]^d}, \varepsilon\right)\}^{\mathbf{m}\left(B_1(A)\right)}. 
\end{align*}
\end{proof}

The following is the main result of this subsection.
This is a modification of a classical result of Bowen \cite[Proposition 2.2]{Bowen entropy expansive}.
Here recall that we have assumed that $\varphi$ is nonnegative.

\begin{proposition} \label{prop: modification of Bowen paper}
For positive numbers $\delta, \beta$ and $0<\varepsilon<1$ there is a positive number $D = D(\delta, \beta, \varepsilon)$
for which the following statement holds.
Set 
\[  a = \frac{\sup_{x\in \mathcal{X}} 
          P_T\left(B_\delta(x, \mathbf{d}_{\mathbb{R}^d}), \mathbf{d}, \varphi, \varepsilon\right)}{\log (1/\varepsilon)}. \]
Then for all sufficiently large $L$ we have 
\[   \sup_{x\in \mathcal{X}} \#\left(B_\delta(x, \mathbf{d}_{[-D, L+D]^d}), \mathbf{d}_L, \varphi_L, \varepsilon\right)
       \leq  \left(\frac{1}{\varepsilon}\right)^{(a+\beta)L^d}. \]
Here $B_\delta(x, \mathbf{d}_{[-D, L+D]^d}) = \{y\in \mathcal{X}\mid \mathbf{d}_{[-D, L+D]^d}(x, y) \leq \delta\}$.
\end{proposition}

\begin{proof}
Choose a positive number $\eta$ satisfying 
\begin{equation}  \label{eq: choice of eta in the proof of local nature}
  \left(\#(\mathcal{X}, \mathbf{d}_{[0, 1]^d}, \varphi_{[0, 1]^d}, \varepsilon)\right)^\eta 
  < \left(\frac{1}{\varepsilon}\right)^{\frac{\beta}{2}}.
\end{equation}
Let $k_0 = k_0(\eta)$ be the natural number introduced in Lemma \ref{lemma: tiles}.

We will construct the following data inductively on $k =1, 2, \dots, k_0$.
\begin{itemize}
    \item   A finite set $Y_k \subset \mathcal{X}$.
    \item   Positive numbers $L_k(y)$ and $M_k(y)$ for each $y\in Y_k$.
    \item   Open neighborhoods $V_k(y)$ and $U_k(y)$ of $y$ in $\mathcal{X}$ with $V_k(y) \subset U_k(y)$ for each $y\in Y_k$.
\end{itemize}
We assume the following conditions.
\begin{enumerate}
   \item  $L_1(y) > 1$ for all $y\in Y_1$.
   \item  $L_k(y) > k_0  L_{k-1}(z)$ for all $y\in Y_k$, $z\in Y_{k-1}$ and $2\leq k \leq k_0$.
   \item  $\#\left(U_k(y), \mathbf{d}_{L_k(y)}, \varphi_{L_k(y)}, \varepsilon\right) < (1/\varepsilon)^{(a+\frac{\beta}{2})L_k(y)^d}$
             for all $y\in Y_k$.
   \item  $B_\delta(v, \mathbf{d}_{[-M_k(y), M_k(y)]^d}) \subset U_k(y)$ for all $y\in Y_k$ and $v\in V_k(y)$.
   \item  $X = \bigcup_{y\in Y_k} V_k(y)$ for every $1\leq k \leq k_0$.         
\end{enumerate}

The construction of these data go as follows.
Suppose that the data of $(k-1)$-th step (i.e. $Y_{k-1}, L_{k-1}(y), M_{k-1}(y), V_{k-1}(y), U_{k-1}(y)$) have been constructed.
We consider the $k$-th step. (The case of $k=1$ is similar.)

Take an arbitrary $y\in \mathcal{X}$. Since we have
\[  \frac{P_T\left(B_\delta(y, \mathbf{d}_{\mathbb{R}^d}), \mathbf{d}, \varphi, \varepsilon\right)}{\log (1/\varepsilon)} \leq 
     a < a + \frac{\beta}{2}, \]
there is a positive number $L_k(y)$ larger than $k_0 \max_{z\in Y_{k-1}} L_{k-1}(z)$ (we assume $L_1(y)>1$ in the case of $k=1$) satisfying
\[  \frac{1}{L_k(y)^d} \log \#\left(B_\delta(y, \mathbf{d}_{\mathbb{R}^d}), \mathbf{d}_{L_k(y)}, \varphi_{L_k(y)}, \varepsilon\right)
      < \left(a+\frac{\beta}{2}\right) \log (1/\varepsilon). \]
Then there is an open set $U_k(y) \supset B_\delta(y, \mathbf{d}_{\mathbb{R}^d})$ such that 
\[ \frac{1}{L_k(y)^d} \log \#\left(U_k(y), \mathbf{d}_{L_k(y)}, \varphi_{L_k(y)}, \varepsilon\right)
      < \left(a+\frac{\beta}{2}\right) \log (1/\varepsilon). \]
Namely we have
\[  \#\left(U_k(y), \mathbf{d}_{L_k(y)}, \varphi_{L_k(y)}, \varepsilon\right) < \left(\frac{1}{\varepsilon}\right)^{(a+\frac{\beta}{2})L_k(y)^d}. \]
Since $B_\delta(y, \mathbf{d}_{\mathbb{R}^d}) = \bigcap_{M=1}^\infty B_{\delta}(y, \mathbf{d}_{[-M, M]^d}) \subset U_k(y)$, 
there is a positive number $M_k(y)$ with $B_\delta(y, \mathbf{d}_{[-M_k(y), M_k(y)]^d}) \subset U_k(y)$.
There is an open neighborhood $V_k(y)$ of $y$ such that for every $v\in V_k(y)$ we have $B_\delta(v, \mathbf{d}_{[-M_k(y), M_k(y)]^d}) \subset U_k(y)$.
Since $\mathcal{X}$ is compact we can find a finite set $Y_k\subset \mathcal{X}$ satisfying $\mathcal{X} = \bigcup_{y\in Y_k} V_k(y)$.
The construction of the $k$-th step has been finished.

Take $D>1$ with $D>\max\{M_k(y)\mid 1\leq k \leq k_0, y\in Y_k\}$.
Let $L$ be a sufficiently large number so that the cube $A :=[0, L]^d$ satisfies 
$\mathbf{m}\left(\partial(A, \max_{y\in Y_{k_0}} L_{k_0}(y))\right) < \frac{\eta}{3} L^d$.
We will show that
$\sup_{x\in \mathcal{X}} \#\left(B_\delta(x, \mathbf{d}_{[-D, L+D]^d}), \mathbf{d}_L, \varphi_L, \varepsilon\right) \leq  \left(\frac{1}{\varepsilon}\right)^{(a+\beta)L^d}$.

Take an arbitrary point $x\in \mathcal{X}$.
For each $1 \leq  k \leq k_0$ and $t\in A\cap \mathbb{Z}^d$ we pick $y\in Y_k$ with $T^t x\in V_k(y)$.
Set $\Lambda_{k, t} = t + [0, L_k(y)]^d$.
It follows from the choice of $D$ that
\[  T^t\left(B_\delta(x, \mathbf{d}_{[-D, L+D]^d})\right) \subset  B_\delta(T^t x, \mathbf{d}_{[-M_k(y), M_k(y)]^d}) \subset  U_k(y). \]
Hence 
\[ \#\left(B_\delta(x, \mathbf{d}_{[-D, L+D]^d}), \mathbf{d}_{\Lambda_{k, t}}, \varphi_{\Lambda_{k,t}}, \varepsilon\right) \leq 
    \#\left(U_k(y), \mathbf{d}_{L_k(y)}, \varphi_{L_k(y)}, \varepsilon\right)  
< \left(\frac{1}{\varepsilon}\right)^{(a+\frac{\beta}{2})L_k(y)^d}. \]
Set $\mathcal{C}_k = \{\Lambda_{k, t}\mid  t\in A\cap \mathbb{Z}^d\}$. This is a finite family of cubes covering $A= [0, L]^d$.
Notice that $\mathcal{C}_k$ depends on the choice of $x$; we suppress its dependence on $x$ in the notation for simplicity.

By Lemma \ref{lemma: tiles} there is a disjoint subfamily $\mathcal{C}^\prime\subset \mathcal{C}_1\cup \dots \cup \mathcal{C}_{k_0}$ such that 
\[  \bigcup_{\Lambda\in \mathcal{C}^\prime} \Lambda \subset A, \quad 
    \mathbf{m}\left(B_1\left(A\setminus \bigcup_{\Lambda\in \mathcal{C}^\prime} \Lambda \right)\right) < \eta \, \mathbf{m}(A). \]
Set 
\[   A^\prime = A\setminus \bigcup_{\Lambda\in \mathcal{C}^\prime} \Lambda.  \]
For every $\Lambda\in \mathcal{C}^\prime$ we have
\[   \#\left(B_\delta(x, \mathbf{d}_{[-D, L+D]^d}), \mathbf{d}_{\Lambda}, \varphi_{\Lambda}, \varepsilon\right) < 
      \left(\frac{1}{\varepsilon}\right)^{(a+\frac{\beta}{2})\mathbf{m}(\Lambda)}. \]
By Lemma \ref{lemma: crude estimate of covering number}
\begin{align*}
  \#\left(B_\delta(x, \mathbf{d}_{[-D, L+D]^d}), \mathbf{d}_{A^\prime}, \varphi_{A^\prime}, \varepsilon\right)
   & \leq  \#\left(\mathcal{X}, \mathbf{d}_{A^\prime}, \varphi_{A^\prime}, \varepsilon\right)  \\
   &  \leq  \left\{\#\left(\mathcal{X}, \mathbf{d}_{[0, 1]^d}, \varphi_{[0, 1]^d}, \varepsilon\right)\right\}^{\mathbf{m}(B_1(A^\prime))} \\
   &  < \left\{\#\left(\mathcal{X}, \mathbf{d}_{[0, 1]^d}, \varphi_{[0,1]^d}, \varepsilon\right)\right\}^{\eta \, \mathbf{m}(A)}  \\
  &  <  \left(\frac{1}{\varepsilon}\right)^{\frac{\beta}{2} \mathbf{m}(A)}.
\end{align*}
In the last inequality we have used (\ref{eq: choice of eta in the proof of local nature}).
From Lemma \ref{lemma: block coding}
\begin{align*}
   & \#\left(B_\delta(x, \mathbf{d}_{[-D, L+D]^d}), \mathbf{d}_A, \varphi_A, \varepsilon\right)  \\
   &  \leq   \#\left(B_\delta(x, \mathbf{d}_{[-D, L+D]^d}), \mathbf{d}_{A^\prime}, \varphi_{A^\prime}, \varepsilon\right) 
             \prod_{\Lambda\in \mathcal{C}^\prime} \#\left(B_\delta(x, \mathbf{d}_{[-D, L+D]^d}), \mathbf{d}_{\Lambda}, \varphi_{\Lambda}, \varepsilon\right)  \\
    & <  \left(\frac{1}{\varepsilon}\right)^{(a+\beta)\mathbf{m}(A)}.
\end{align*}
This holds for every point $x\in \mathcal{X}$.
Thus we have proved the claim of the proposition.
\end{proof}

\subsection{Proof of Theorem \ref{theorem: local formula of metric mean dimension with potential}}
\label{subsection: proof of the local formula}

Here we prove Theorem \ref{theorem: local formula of metric mean dimension with potential}.
Let $T\colon \mathbb{R}^d\times \mathcal{X}\to \mathcal{X}$ be a continuous action of $\mathbb{R}^d$ on a compact metrizable 
space $\mathcal{X}$.
Let $\mathbf{d}$ be a metric on $\mathcal{X}$ and $\varphi\colon \mathcal{X}\to \mathbb{R}$ a continuous function.
We do not assume that $\varphi$ is nonnegative.

The next proposition looks the same as Proposition \ref{prop: modification of Bowen paper}.
The point is that we do not assume the nonnegativity of $\varphi$ here whereas we assumed it in 
Proposition \ref{prop: modification of Bowen paper}.

\begin{proposition} \label{prop: modification of Bowen paper revisited}
Let $\delta, \beta, \varepsilon$ be positive numbers with $0<\varepsilon <1$.
Set 
\[ a = \sup_{x\in \mathcal{X}} \frac{P_T\left(B_\delta(x, \mathbf{d}_{\mathbb{R}^d}), \mathbf{d},\varphi, \varepsilon\right)}{\log (1/\varepsilon)}. \]
Then for all sufficiently large $L$ we have 
\[ \sup_{x\in \mathcal{X}}\#\left(B_\delta(x, \mathbf{d}_{[-D, L+D]^d}), \mathbf{d}_L, \varphi_L, \varepsilon\right) 
    \leq  \left(\frac{1}{\varepsilon}\right)^{(a+\beta)L^d}. \]
Here $D = D(\delta, \beta, \varepsilon)$ is the positive constant\footnote{Strictly speaking, the constant $D$ depends on not only 
$\delta, \beta, \varepsilon$ but also $(\mathcal{X},T, \mathbf{d}, \psi)$ where $\psi := \varphi - \min_{\mathcal{X}}\varphi$.} 
introduced in Proposition \ref{prop: modification of Bowen paper}.
\end{proposition}

\begin{proof}
Set $c = \min_{x \in \mathcal{X}}\varphi(x)$ and $\psi(x) = \varphi(x)-c$. 
We have $\psi(x) \geq 0$. For any positive number $L$ we have
\[ \psi_L(x) = \varphi_L(x) - cL^d. \]
For any subset $E\subset \mathcal{X}$
\[ \#\left(E, \mathbf{d}_L, \psi_L, \varepsilon\right) = (1/\varepsilon)^{-c L^d} \#\left(E, \mathbf{d}_L, \varphi_L, \varepsilon\right). \]
Hence 
\[  P_T(E, \mathbf{d}, \psi, \varepsilon) = P_T(E,\mathbf{d},\varphi, \varepsilon) - c\log (1/\varepsilon),   \]
\[   \sup_{x\in \mathcal{X}} \frac{P_T\left(B_\delta(x, \mathbf{d}_{\mathbb{R}^d}), \mathbf{d},\psi,\varepsilon\right)}{\log (1/\varepsilon)}   
       = a-c. \]
Since $\psi$ is a nonnegative function, we apply Proposition \ref{prop: modification of Bowen paper} and get 
\[ \sup_{x\in \mathcal{X}}\#\left(B_\delta(x, \mathbf{d}_{[-D,L+D]^d}), \mathbf{d}_L, \psi_L, \varepsilon\right) 
    \leq    \left(\frac{1}{\varepsilon}\right)^{(a-c+\beta)L^d} \]
for sufficiently large $L$. 
We have 
\[  \#\left(B_\delta(x, \mathbf{d}_{[-D,L+D]^d}), \mathbf{d}_L, \psi_L, \varepsilon\right) 
     = \left(\frac{1}{\varepsilon}\right)^{-c L^d} \#\left(B_\delta(x, \mathbf{d}_{[-D,L+D]^d}), \mathbf{d}_L, \varphi_L, \varepsilon\right). \]
Therefore 
\[ \sup_{x\in \mathcal{X}} \#\left(B_\delta(x, \mathbf{d}_{[-D,L+D]^d}), \mathbf{d}_L, \varphi_L, \varepsilon\right)
   \leq \left(\frac{1}{\varepsilon}\right)^{(a+\beta)L^d}. \]
\end{proof}

Now we prove Theorem \ref{theorem: local formula of metric mean dimension with potential}.
We write the statement again.

\begin{theorem}[$=$ Theorem \ref{theorem: local formula of metric mean dimension with potential}]
For any positive number $\delta$
  \begin{equation}  \label{eq: local formula of metric mean dimension with potential}
   \begin{split}
   \umdimm(\mathcal{X}, T, \mathbf{d}, \varphi) = \limsup_{\varepsilon\to 0} 
   \frac{\sup_{x\in \mathcal{X}} P_T\left(B_\delta(x, \mathbf{d}_{\mathbb{R}^d}), \mathbf{d},\varphi,\varepsilon\right)}{\log (1/\varepsilon)}, \\
    \lmdimm(\mathcal{X}, T, \mathbf{d}, \varphi) = \liminf_{\varepsilon\to 0} 
   \frac{\sup_{x\in \mathcal{X}} P_T\left(B_\delta(x, \mathbf{d}_{\mathbb{R}^d}), \mathbf{d},\varphi,\varepsilon\right)}{\log (1/\varepsilon)}.
   \end{split}
  \end{equation}
\end{theorem}

\begin{proof}
It is obvious that the left-hand sides of (\ref{eq: local formula of metric mean dimension with potential})
are greater than or equal to the right-hand sides.
So it is enough to prove the reverse inequalities.
Let $\beta$ and $\varepsilon$ be arbitrary positive numbers with $0<\varepsilon<1$.
Let $D= D(\delta, \beta, \varepsilon)$ be the positive constant introduced in Proposition \ref{prop: modification of Bowen paper}.
Set 
\[ a = \sup_{x\in \mathcal{X}} \frac{P_T\left(B_\delta(x, \mathbf{d}_{\mathbb{R}^d}), \mathbf{d},\varphi, \varepsilon\right)}{\log (1/\varepsilon)}. \]
For any positive number $L$ we can take points $x_1, \dots, x_M\in \mathcal{X}$ such that
\[   \mathcal{X} = \bigcup_{m=1}^M  B_\delta(x_m, \mathbf{d}_{[-D, L+D]^d}),  \]
\[   M\leq  \#\left(\mathcal{X}, \mathbf{d}_{[-D,L+D]^d}, \delta\right) = \#\left(\mathcal{X}, \mathbf{d}_{[0, L+2D]^d},\delta\right).  \]
Then we have 
\begin{align*}
  \#\left(\mathcal{X}, \mathbf{d}_L, \varphi_L, \varepsilon\right) & \leq \sum_{m=1}^M 
  \#\left(B_\delta(x_m, \mathbf{d}_{[-D, L+D]^d}), \mathbf{d}_L, \varphi_L, \varepsilon\right) \\
  & \leq  M\, \sup_{x\in \mathcal{X}} \#\left(B_\delta(x, \mathbf{d}_{[-D, L+D]^d}), \mathbf{d}_L, \varphi_L, \varepsilon\right) \\
  & \leq  M \left(\frac{1}{\varepsilon}\right)^{(a+\beta)L^d}.
\end{align*}
The last inequality holds for all sufficiently large $L$ by Proposition \ref{prop: modification of Bowen paper revisited}.
Therefore 
\[ \log \#\left(\mathcal{X}, \mathbf{d}_L, \varphi_L, \varepsilon\right)  \leq  \log \#\left(\mathcal{X}, \mathbf{d}_{L+2D}, \delta\right)
    + (a+\beta) L^d \log (1/\varepsilon). \]
Dividing this by $L^d$ and letting $L\to \infty$, we have 
\begin{align*}
   P_T(\mathcal{X}, \mathbf{d}, \varphi, \varepsilon) & \leq  \lim_{L\to \infty} \frac{\log \#\left(\mathcal{X}, \mathbf{d}_L, \delta\right)}{L^d}
                                                                                 + (a+\beta) \log (1/\varepsilon)  \\
  & \leq  \log \#\left(\mathcal{X}, \mathbf{d}_1, \delta\right)  + (a+\beta) \log (1/\varepsilon).
\end{align*}
We can let $\beta\to 0$ and get
\begin{align*}
   P_T(\mathcal{X}, \mathbf{d},\varphi, \varepsilon) & \leq 
   \log \#\left(\mathcal{X}, \mathbf{d}_1, \delta\right)   + a \log (1/\varepsilon) \\
   & =
   \log \#\left(\mathcal{X}, \mathbf{d}_1, \delta\right)  
    + \sup_{x\in \mathcal{X}} P_T\left(B_\delta(x, \mathbf{d}_{\mathbb{R}^d}), \mathbf{d}, \varphi, \varepsilon\right). 
\end{align*}    
We divide this by $\log (1/\varepsilon)$ and let $\varepsilon \to 0$.
Then we conclude that the left-hand sides of (\ref{eq: local formula of metric mean dimension with potential}) are less than or equal to the 
right-hand sides.
\end{proof}

\end{document}